\documentclass[a4paper,oneside,11pt]{article}
 \usepackage[utf8]{inputenc} 
\usepackage{amsfonts}
\usepackage{amsthm}
\newenvironment{customproof}{\noindent\emph{Proof.}}{\qed}
\newenvironment{customproofsketch}{\noindent\emph{Proof} (sketch).}{\qed}
\usepackage{graphicx}
\usepackage{pgfplots}
\usepackage{subfigure}
\usepackage{tikz}
\usepackage{hyperref}
\usepackage{dsfont}
\usepackage[toc,page]{appendix}
\usepackage{enumitem}
\usepackage{bm,amsmath,tabularx}
\newcolumntype{L}[1]{>{\raggedright\arraybackslash}p{#1}}
\newcolumntype{C}[1]{>{\centering\arraybackslash}p{#1}}
\newcolumntype{R}[1]{>{\raggedleft\arraybackslash}p{#1}}

\newcommand{\Prob}{\mathbb{P}}

\usepackage{rotating}
\usepackage{multirow}
\usepackage{float}

\newcommand{\eqdef}{\mathrel{\mathop:}=}

\definecolor{dgreen}{rgb}{0,0.7,0}
\definecolor{dred}{rgb}{0.8,0,0}
\definecolor{dblue}{rgb}{0,0,0.8}

\newcommand{\1}{\mathds{1}}
\usepackage{bm}
\usepackage{mathtools}
\usepackage{empheq}
\usepackage{tikz}
\usepackage{subcaption}
\usepackage{booktabs}
\newcommand{\bP}{\mathbb{P}}
\newcommand{\cP}{\mathcal{P}}

\newcommand{\bE}{\mathbb{E}}

\usepackage{dsfont}
\usepackage{stackengine}
\usepackage{accents}

\newcommand{\usigma}{\bar{\sigma}}
\newcommand{\lsigma}{\underaccent{\bar}{\sigma}}
\newcommand{\uupsilon}{\bar{\upsilon}}
\newcommand{\lupsilon}{\underaccent{\bar}{\upsilon}}

\usepackage{amsthm}
\newtheorem{theorem}{Theorem}[section]  % This defines a theorem environment
\newtheorem{lemma}{Lemma}[section]      % This defines a lemma environment
\newtheorem{proposition}{Proposition}[section] % This defines a proposition environment
\newtheorem{corollary}{Corollary}[section]

\usepackage{natbib}
 \bibpunct[, ]{(}{)}{,}{a}{}{,}%
\usepackage{tikz}
\usetikzlibrary{shapes,arrows,positioning}
\usetikzlibrary{automata}
\usepackage{bbm}
\makeatletter
\newcommand{\thickbar}{\mathpalette\@thickbar}
\newcommand{\@thickbar}[2]{{#1\mkern1.5mu\vbox{
  \sbox\z@{$#1\mkern-1.5mu#2\mkern-1.5mu$}%
  \sbox\tw@{$#1\overline{#2}$}%
  \dimen@=\dimexpr\ht\tw@-\ht\z@-.8\p@\relax
  \hrule\@height.8\p@ % adjust for the desired rule thickness
  \vskip\dimen@
  \box\z@}\mkern1.5mu}
}
\makeatother
\usetikzlibrary{arrows,intersections}

\usepackage{a4wide}
\begin{document}

\thispagestyle{empty}

\title{Distributionally robust monopoly pricing:\\
  Switching from low to high prices in volatile markets
  }

\author{Tim S.G. van Eck, Pieter Kleer and Johan S.H.~van Leeuwaarden}
  \maketitle

\begin{center}
\textit{Manufacturing \& Service Operations Management}\\
Published version:
\href{https://pubsonline.informs.org/doi/10.1287/msom.2024.0952}{https://pubsonline.informs.org/doi/10.1287/msom.2024.0952}
\end{center}

\begin{abstract}
\noindent {\bf Problem definition:} 
Traditional monopoly pricing assumes sellers have full information about consumer valuations. 
We consider monopoly pricing under limited information, where a seller only knows the mean, variance and support
of the valuation distribution. The objective is to maximize expected revenue by selecting the optimal fixed price. 
{\bf Methodology/results:} 
We adopt a distributionally robust framework, where the seller considers all valuation distributions that comply with the limited information. We formulate a maximin problem which seeks to maximize expected revenue for the worst-case valuation distribution. The minimization problem that identifies the worst-case valuation distribution is solved using 
 primal-dual methods, and in turn leads to an   explicitly solvable  maximization problem. This yields a closed-form optimal pricing policy and a new fundamental principle prescribing when to use low and high robust prices.\\ 
{\bf Managerial implications:} 
 We show that the optimal policy switches from low to high prices when variance becomes sufficiently large, yielding significant performance gains compared with existing robust prices that generally decay with market uncertainty. This presents guidelines for when the seller should switch from targeting mass markets to niche markets. Similar guidelines are obtained for delay-prone services with rational utility-maximizing customers, underlining the universality and wide applicability of the novel pricing policy. 
 \end{abstract}

\section{Introduction}
Optimizing revenue in contemporary markets represents a pivotal challenge, extensively explored within theoretical pricing models spanning disciplines such as operations research, computer science, and economics. Prevailing in the theoretical pricing literature is the assumption that sellers possess complete knowledge of how the market responds to price variations. This assumption relies on a full understanding of the demand function, serving as a prior in Bayesian analysis, and significantly influencing the derived optimal price. However, the theoretical advantage of assuming full knowledge is compromised by the scarcity or unreliability of market data. A principle known as Wilson's doctrine calls for less detailed yet parsimonious models, and the development of more resilient pricing theories, capable of accommodating imprecise market details \citep{wilson1987,carroll2019robustness}.
Robust pricing theories seek strategies that perform well across diverse markets, reducing sensitivity to detailed market knowledge. The interplay between a seller's pricing strategy and market knowledge then becomes crucial. A complete lack of knowledge may result in an overly pessimistic worst-case scenario, rendering a robust price excessively conservative and impractical. A more realistic approach is to assume partial knowledge. This paper addresses this scenario in the context of a monopolist selling a single type of good to a diverse, independent consumer base.

The monopolist aims to determine the revenue-maximizing price by balancing aggressive (low) pricing for a larger market share against higher pricing for a smaller market share. Each consumer's valuation of the good is an independent draw from a given valuation distribution with the purchase occurring when this valuation exceeds the price. The expected revenue is defined as the product of the price times the probability that the valuation exceeds the price.
While the optimal price is easily determined when the valuation distribution is known \citep{riley1983optimal, myerson1981optimal}, the challenge arises when only partial market knowledge is available. 
This study focuses on the case when only the mean, variance, and an upper limit of the valuation distribution are known. Mean-variance information proves suitable for robust decision-making, capturing both the anticipated future scenario (mean) and the risk associated with less likely scenarios (variance). This choice is motivated by the often scarce or non-existent data faced by sellers, where the mean, variance, and support upper bound can be based on expert opinions—intuitive statistics that are easier to estimate than the entire valuation distribution; see e.g.~\cite{Scarf1958,azar2012optimal}.

The partially informed seller needs to decide on the price in the face of uncertain market demand, much in the same way as an investor who only knows the risk and return of assets. To do so, the seller employs a maximin decision framework. The maximin price, or robust price, is determined by first establishing the worst-case valuation distribution that minimizes expected revenue for all prices and then maximizing expected revenue by setting the optimal price. This involves solving two optimization problems—a semi-infinite linear program and finding the global optimum of a bimodal function that depends delicately on the price and information parameters. The comprehensive nature of this approach reveals the intricate interplay between partial knowledge and robust pricing, surpassing earlier studies \citep{azar2012optimal,kos2015selling,suzdaltsev2018distributionally,chen2022distribution,chen2021screening} and highlighting a novel robust pricing principle that emerges only when knowledge includes both the support upper bound and variance.

\begin{figure}[t]
\centering
\subfigure[Optimal price]{
\label{}
\begin{tikzpicture}
    \node at (0,0) {\includegraphics[width=0.45\linewidth]{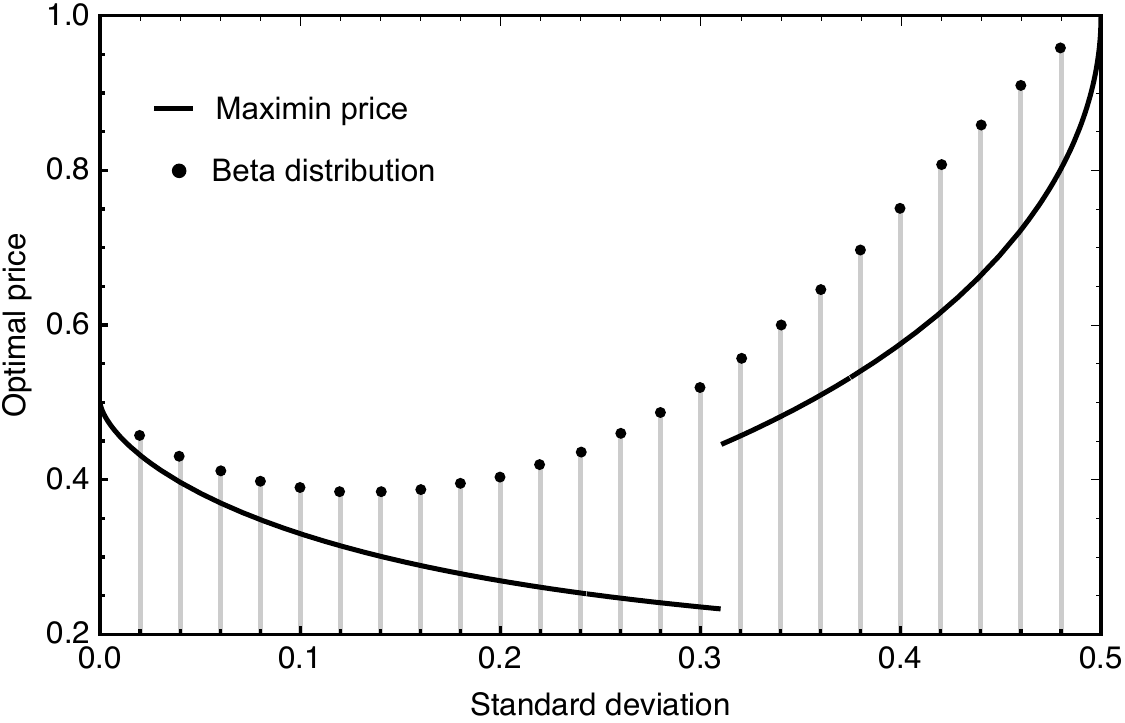}};
    \end{tikzpicture}
}
\subfigure[Maximal revenue]{
\label{}
\begin{tikzpicture}
    \node at (0,0) {\includegraphics[width=0.45\linewidth]{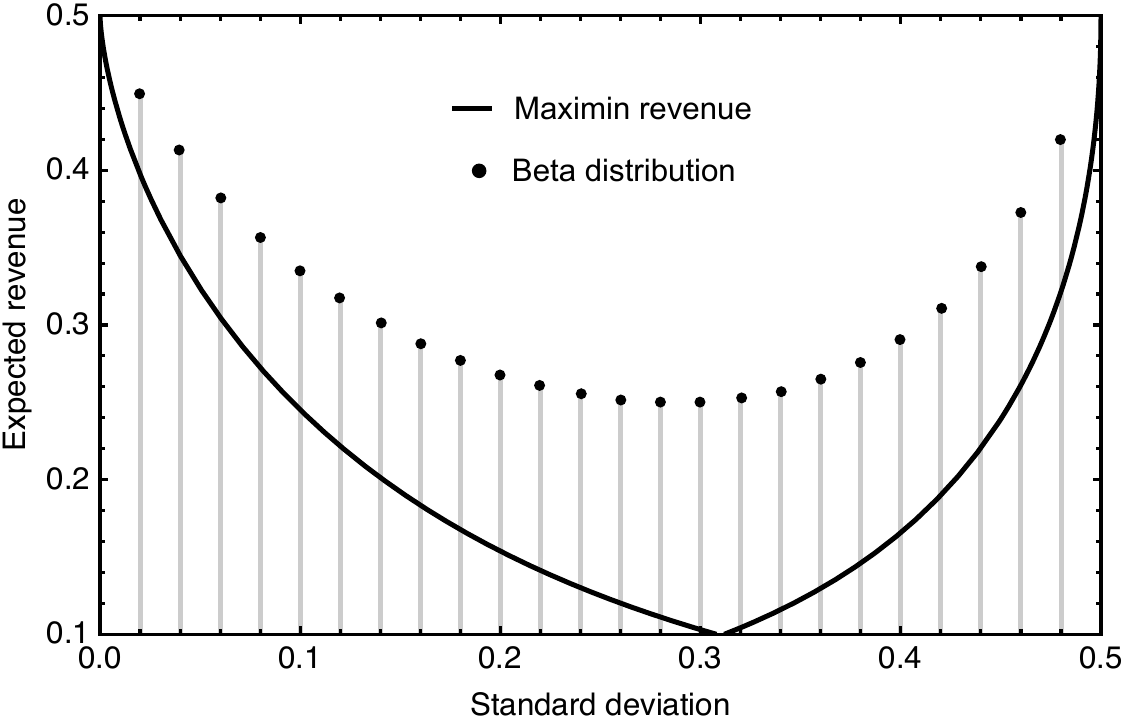}};
    \end{tikzpicture}
}
\caption{
Optimal price and revenue for full information setting when the seller knows the valuation follows a beta distribution on $[0,1]$ with mean $1/2$ and standard deviation ranging from $0$ to $0.5$. The black curves represent the robust prices and revenue revealed in this paper.}\label{fig:beta}
\end{figure}

\vspace{.3cm}
\noindent {\bf Classical versus robust pricing.}
Before we explain this general principle, consider classical monopoly pricing with a fully informed monopolist who knows that valuations  follow some cumulative distribution function $F$ with probability density function $f$. For any given price $p$, expected revenue is $p(1-F(p))$ and the monopolist thus solves the maximization problem
$\sup_p p\left(1-F(p)\right)$. First-order conditions then prescribe that the optimal price $p^*$ solves
\begin{flalign}\label{foc}
    &&p^* = \frac{1-F(p^*)}{f(p^*)}.&&
\end{flalign}
Figure \ref{fig:beta} shows this optimal price for valuations that follow a beta distribution on the unit interval. 
Observe that the optimal price is $\mu$ for zero variance, a feature that holds irrespective of the specific valuation distribution. Another feature that can be observed which provided a strong motivation for our work is the following: Viewed as function of $\sigma$, the optimal price $p^*$ first decreases to values below the mean and then increases to values above the mean and eventually to the maximal valuation.
 
The ``U-shaped" nature (or quasi-convexity) of revenue as function of demand dispersion -- as displayed in the right panel of Figure~\ref{fig:beta} -- is a known phenomenon in monopolistic markets \citep{johnson2006simple,sun2012does}. Rigorous underpinning is provided in \cite{johnson2006simple} by an explanatory analysis that shows the U-shape occurs naturally by imposing only mild conditions on the valuation distribution. Consequently, when there is sufficiently low dispersion, indicative of relatively uniform consumer valuations, sellers are inclined to cater to ``mass markets" through competitive pricing strategies. Conversely, as dispersion increases, the average consumer's willingness to pay diminishes. The increased heterogeneity among consumers prompts firms to adopt a strategy of restricting sales to a more specialized ``niche market characterized by higher prices. We add a theory aligned with the above intuition, showing that for a rich class of potential consumer valuations, a robust seller indeed switches from low pricing with large market share to high pricing and small market share in response to increased consumer valuation dispersion.

To build some first intuition for the robust pricing in this paper, assume that the seller only knows the mean, variance, and support upper bound, and chooses to fit to these summary statistics a beta distribution. Since the monopolist now only pretends to know the valuation distribution, the optimal price found from assuming a beta distribution is only an approximation, as the ground truth valuation is likely not to be beta distributed. In fact, the beta distribution is only one of infinitely many distributions with the same mean, variance and support, and only serves to calculate the expected revenue. This paper presents an alternative pricing strategy, where the monopolist does not commit to a specific valuation distribution, and instead first determines the worst-case valuation distribution that minimizes demand, and only then solves for the revenue-maximizing price. We prove that this robust maximin price corresponds to the black solid line in Figure~\ref{fig:beta}. As this paper's title says, the monopolist switches from low to high prices when the variance exceeds a certain threshold. Before that threshold, the maximin price decreases with variance, at the threshold the price jumps to a higher value, and continues to increase after the threshold. 

\vspace{.3cm}
\noindent {\bf Worst-case demand implies segmentation.}
It turns out that the structure of the worst-case valuation distribution identified in this study suggests implicit market segmentation. Specifically, this worst-case market can be interpreted as comprising a low-end segment of customers willing to pay only up to a certain amount and a high-end segment willing to pay significantly more. While detailed knowledge of these segments might lead a seller to consider third-degree price discrimination by setting distinct prices for each group, the seller in this setting only has limited information available. Consequently, price differentiation is not addressed here; the analysis focuses on determining the optimal uniform price for all customers. Due to customers' perceptions of fairness, however, firms often adopt uniform pricing even in markets with an underlying segmented structure. \cite{dellavigna2019uniform} provide substantial evidence that food and mass-merchandise chains often charge nearly uniform prices across stores, despite significant variations in consumer demographics.  Beyond fairness perceptions, uniform pricing can perform well in terms of profitability under specific conditions in segmented markets. When the revenue functions for each segment are concave and share a common support (i.e., the same price range applies across all segments), \cite{bergemann2022third} demonstrate that uniform pricing captures at least 50\% of the maximum profits achievable under price discrimination. This result holds irrespective of the number of segments or the specific discriminatory prices. However, this profit guarantee collapses when the assumption of common support is violated, underscoring that the degree of market segmentation—particularly differences in price support—has a critical impact on the effectiveness of uniform pricing.

\vspace{.3cm}
\noindent {\bf Jumps in sticky prices.} Our analysis focuses on a firm selling to a large market without detailed knowledge of the demand function or its potential segmentation. Adopting a robust approach, the firm determines an optimal yet "sticky" uniform price based on the mean, variance, and support of demand. In the worst-case scenario, demand corresponds to a segmented market where segments differ in price support, making the revenue function fall outside the scope of \cite{bergemann2022third}. These differences in support provide strong incentives to adjust the uniform price to better target certain segments. However, such adjustments must be carefully balanced against the risks of customer anger and brand equity concerns when revising prices. For instance, in 2023 and 2024, Tesla implemented significant price reductions across its flagship models, marking a strategic shift from targeting niche markets of early adopters to a broader mass-market audience \citep{reuters_us_2024}.
 By lowering entry barriers, Tesla aimed to attract middle-income buyers previously excluded from the electric vehicle market. However, these abrupt price cuts provoked anger among existing Tesla owners, who felt cheated and faced reduced resale values for their vehicles.  Despite this backlash, the strategic objective of expanding market reach was prioritized. In our analysis, which solely focuses on revenue maximization and disregards factors such as customer anger and brand reputation, we also identify the triggers for when significant price reductions become profitable, particularly in targeting lower-end segments; in Figure~\ref{fig:beta} this trigger is the threshold value for the variance, where the optimal price jumps for low to high or vice versa.

\vspace{.3cm}
\noindent {\bf Contributions.}
We reveal a new phenomenon of abruptly switching from low to high pricing when the market segmentation becomes sufficiently large. In this way, we contribute to advancing our understanding of uniform pricing in the presence of limited market information.
\cite{azar2012optimal} pioneered this direction for knowledge of the mean and variance, demonstrating that the optimal uniform price maximizing expected revenue in the worst case can be expressed as an explicit function that decays with the variance. Our more general results align with theirs, and their pricing function appears prominently in several of our theorems for low-risk markets. Another relevant study by \cite{kos2015selling} considers information about the mean and bound on support, so the seller lacks knowledge about the variance. \cite{kos2015selling} also derive an explicit pricing function which appears in our study when we consider high-risk markets. We summarize the existing works and our new results in Table~\ref{table_2}.

\begin{table}[t]
\begin{center}
    \begin{tabular}{  p{2.5cm} | p{4cm} | p{7cm} }
    \hline
    {\bf Information} & {\bf Pioneering works} & {\bf Pricing insights} \\ \hline
    \vtop{\hbox{\strut Mean}
    \hbox{\strut Variance}}
     & \cite{azar2012optimal} \cite{carrasco2018optimal}  & {\small Azar and Micali's pioneering work revealed a robust price that decreases with variance. Carrasco et al.~considered an upper bound on variance, yielding the same price.}  \\ \hline
     \vtop{\hbox{\strut Mean}\hbox{\strut Support}} & \cite{kos2015selling}  & {\small Kos en Messner found a robust price that increases with the valuation support and is (by definition) insensitive to variance.} 
 \\ \hline
     \vtop{\hbox{\strut Mean}
    \hbox{\strut Variance}\hbox{\strut Support}}& 
    \vtop{\hbox{\strut \cite{suzdaltsev2018distributionally}}
    \hbox{\strut This paper:}\hbox{\strut -Theorems~\ref{3pthm} and \ref{general_jump_in_words}}\hbox{\strut -Proposition~\ref{theorem_s} and \ref{theorem_h}}
    \hbox{\strut -Theorems~\ref{th:queue}, \ref{THMswithQ}, \ref{th:queue2}, \ref{THMswithQ2}}}& {\small This paper and unpublished work of Suzdaltsev present robust prices that first decrease and then increase with variance, causing the switch from low to high pricing strategies in response to increased market volatility. The high pricing strategy targets niche markets.} \\ \hline
    \end{tabular}
\end{center}
\vspace{.1cm}
\caption{Pioneering papers and contributions in this paper.
}
\label{table_2}
\end{table}

From a more theoretical perspective, we rigorously establish several propositions and theorems (see Table~\ref{table_2}) for both the classical pricing problem and pricing in an extended setting with delay-prone services. These theoretical results together provide a firm underpinning of this principle of low and high robust pricing. In all cases, we derive the optimal uniform prices for when the seller has full information of the mean, a support upper bound, and a variance range of the valuation distribution, generalizing all earlier works including \cite{azar2012optimal,kos2015selling}. 
We obtain explicit expressions for the optimal prices by solving a maximin problem with distributionally robust optimization methods. We solve the minimization problem with a classical primal-dual method, see e.g.~\cite{rogosinski1958moments, shapiro2001duality}, and show that the solution renders a bimodal revenue function with up to three candidate prices that achieve local maxima. Solving the maximization problem then requires
showing that the global optimum is always one of these three local optima, ruling out boundary points; 
the robust seller then always chooses one of these three candidate prices, a crucial insight that generalizes all existing contributions listed in Table~\ref{table_2}.

\vspace{.3cm}
\noindent {\bf Outline.}
The paper is structured as follows. Section~\ref{sec:lit} presents an overview of related literature on monopoly pricing and distributionally robust optimization. 
Section~\ref{sec:pre} obtains the worst-case market in terms of tight tail bounds for the valuation distribution, and presents the main result for low-high pricing by showing that the partially informed seller always selects one of three closed-form pricing functions. 
This result is proven in Section~\ref{sec:twoprices}, starting with first presenting results on robust pricing and the low-high pricing phenomenon for special cases. 
Section~\ref{sec:queue} treats pricing in delay-prone services and shows that similar robust pricing strategies are optimal. 
Section~\ref{sec:con} concludes and reflects on future work.

\section{Related literature}\label{sec:lit}
Our study adds to the stream of research on robust pricing decisions under limited information, where sellers must commit to fixed prices without possessing complete knowledge of market valuation and the corresponding demand function.
Our approach focuses on mean-variance information, a widely studied information set in applications such as investment theory \citep{natarajan2017reduced}, option pricing \citep{de2007distribution}, insurances \citep{roos2019chebyshev}, service operation management \citep{vaneekelen2021mad,taylor2018demand} and inventory management \citep{Scarf1958,gallego1992minmax}. In the management science and operations literature, this information setting is commonly referred to as `distribution-free', highlighting that the decision-making under uncertainty does not require detailed or full distributional knowledge. 

\vspace{.2cm}
\noindent{\bf Distributionally robust optimization.}
Such distribution-free analysis combined with optimization falls in the more general domain of distributionally robust optimization (DRO), a mathematical framework that addresses uncertainty in optimization problems by considering an ambiguity set of possible probability distributions for the involved random variables, rather than relying on a single assumed distribution. The ambiguity sets considered in this paper are based on moment information, as in the pioneering work of \cite{Scarf1958} on robust inventory management, and later generalized to a powerful method with many applications \citep{Delage2010,wiesemann2014distributionally}.  Applied to monopoly pricing, this moment-based DRO approach helps create more robust and conservative pricing strategies that perform well under a range of scenarios. Instead of making strong assumptions about the specific valuation distribution, DRO considers a set of possible distributions, forming a certain ambiguity set, and resulting in pricing strategies that are robust across different scenarios, without relying on a single, possibly incorrect, assumed distribution. 
This aligns well with Occam's razor and Wilson's doctrine \citep{carroll2019robustness}, two principles of parsimony that aim for explanatory models able to capture the essence without making unnecessary assumptions or requiring overly strong parametric assumptions.

\vspace{.2cm}
\noindent{\bf Robust monopoly pricing.}
Azar and Micali's foundational result for the information set with mean and variance (and no bound on the support) 
has been rediscovered and rederived in various papers \citep{carrasco2018optimal,chen2022distribution,suzdaltsev2018distributionally}.
The information set consisting of the triplet mean, variance, and bound on support has remained largely unexplored. A notable exception is \cite{suzdaltsev2018distributionally}, where this setting was studied as a special single-bidder case of an auction. In fact, one of the results in this paper (Proposition~\ref{theorem_h}) was mentioned without proof in an early draft of \cite{suzdaltsev2018distributionally} and was removed from the final published version. 
We thus focus on maximizing expected revenue with information about the first two moments, just as in \cite{azar2012optimal,chen2022distribution,suzdaltsev2018distributionally}. Alternative information sets that were studied in the context of robust pricing include knowledge of higher moments \citep{carrasco2018optimal}, percentiles \citep{eren2010monopoly}, mean absolute deviation \citep{roos2021distributionally, elmachtoub2021value},
mean preserving contraction and spread \citep{condorelli2020information,roesler2017buyer},
or knowing that the valuation distribution is within the proximity of a given reference distribution \citep{bergemann2011robust,chen2021screening}. Instead of maximin expected revenue, alternative objectives studied include minimax regret \citep{bergemann2008pricing,bergemann2011robust,elmachtoub2021value, chen2021screening} and the competitive or approximation ratio \citep{eren2010monopoly,giannakopoulos2019robust,elmachtoub2021value,wang2024minimax}.
Let us also mention that our work  connects to 
 Bayesian persuasion \citep{kamenica2011bayesian,kamenica2019bayesian}, a research line that explores the optimal type or level of information for optimizing individual utility or social welfare.

\vspace{.2cm}

\noindent {\bf Regularity assumptions.}
The lion's share of the pricing literature assumes strong regularity conditions on the valuation distribution that guarantee existence and uniqueness of an optimal price. In particular, assuming that the valuation has a monotone hazard rate, or some closely related condition, can guarantee that the revenue function is concave and has a unique optimizer \citep{ziya2004relationships,schweizer2019performance}. 
Such regularity assumptions greatly simplify the mathematical framework, but true valuation distributions are more often irregular than not \citep{morganti2022non,roughgarden2016ironing}. 
\cite{loertscher2022monopoly} write:  ``The standard assumption, which is almost universally maintained in economics, is that expected revenue is concave. The typical justification for this assumption, other than it being standard, is that it is deemed an analytic simplification that permits one to focus on the key economic insights without cluttering the analysis with case distinctions and multiplicity of local maxima. We have never seen it justified on the basis of empirical evidence, and we will not impose it.''  See \cite{celis2014buy} for empirical evidence that the value distribution need not be regular. See also the literature overview in \cite{loertscher2022monopoly} for more references on non-concave revenue and a proof that whenever two or more separate markets merge, the revenue is no longer concave. Even in robust settings, regularity conditions are often maintained for mathematical tractability,  see e.g.~\cite{maskin1984monopoly}, \cite{eren2010monopoly}, \cite{allouah2022pricing}. 
For robust pricing as in this paper, one should let go of regularity completely, as the worst-case distribution that plays a pivotal role in the maximin analysis typically fails this assumption. In fact, if you would add to the minimization problem the condition that the worst-case valuation needs to be regular, we would lose the key take-aways and novel economic insights of this paper. 
\vspace{.2cm}

\noindent {\bf Pricing in service operations.}
Monopoly pricing in the context of delay-sensitive service systems has been extensively examined. A well-established body of literature, as detailed by \cite{hassin2003queue}, delves into mathematical models that capture customer decision-making regarding joining delay-sensitive services based on individual utility considerations. The pioneering work by \citet{naor1969regulation} laid the groundwork for this branch of operations management. Naor's seminal model views the system as an M/M/1 queue and models customer utility as linearly dependent on the price, value of service, and the anticipated waiting time. Naor revealed that rational individuals have a tendency to overcrowd the system and suggested the introduction of a fixed fee for service that optimizes social welfare. The determination of this socially optimal fee involves solving first-order conditions, akin to classical monopoly pricing. Expanding upon the foundational work by Naor, numerous studies investigated pricing and control considerations for a variety of extended queueing scenarios, with often the focus on determining the price that maximizes revenue or welfare; see
e.g.~\cite{dewan1990user,chen2004monopoly,ha1998incentive,taylor2018demand,baron2023omnichannel,chen2022food}. \cite{stidham1992pricing} made a significant contribution by introducing the concept of random valuation for individual customers and assuming that the arrival rate of customers with a valuation surpassing a specified threshold. In this framework, valuation is often referred to as willingness-to-pay, and as in classical monopoly pricing, is a random variable with some given distribution.
\cite{cachon2011pricing,haviv2014pricing} build on Stidham's framework and assume that rational customers join when their random valuation exceeds the sum of the price and the expected delay costs. This requires specifying the valuation distribution in order to determine the optimal price. We will perform a distribution-free analysis of this rational queueing model and show that again the optimal pricing strategy switches from low to high pricing when the variance crosses a certain threshold. Another stream of research on rational queues assumes that customers are heterogeneous in their delay sensitivity rather than valuation  \citep{Afeche2013, Ata2018, Yu2013, Aksin2013, Nazerzadeh2018}.
We will also establish the low-high pricing phenomenon for this setting with a partially known delay-sensitivity distribution.

\section{Robust pricing framework}\label{sec:pre}
The optimal Bayesian price $p^*$ in \eqref{foc} strongly depends on the properties of the given valuation distribution. The robust pricing developed in this paper, however, only uses summary statistics and renders a revenue-guarantee accomplished across all possible valuation distributions, without referring to any specific distribution.
Therefore, the basic problem for robust monopoly pricing considered in this paper is
\begin{flalign}\label{rmpp}
&&\sup_{p}R(p,\cP),&&
\end{flalign} where
\begin{flalign}\label{rr}
    &&R(p,\cP) \eqdef \inf_{\bP\in\cP}p \bP(X\geq p)&&
\end{flalign}
with $\bP$ a valuation distribution from $\cP$, the ambiguity set with all valuation distributions that comply with the partial knowledge that the seller has. Specifically, we consider the ambiguity set (with $\mathcal{B}$ the Borel sigma-algebra over $\mathbb{R_+}$)
\begin{flalign}\label{ambset}
    &&\cP(\mu,\lsigma,\usigma,\beta) = \left\{\bP:\mathcal{B}\to[0,1]\mid \bP(X\in[0,\beta])=1,\, \bE(X) = \mu,\bE(X^2) \in [\lsigma^2+\mu^2, \usigma^2 + \mu^2]\right\}&&
\end{flalign}
containing all distributions with mean $\mu$, standard deviation in the interval $[\lsigma, \usigma]$ and support upper bounded by $\beta$. So instead of knowing an upper bound or exact value, the seller knows that the standard deviation is contained in an interval $[\lsigma,\usigma]$. The width of this interval thus reflects how confident the seller is about the actual value $\sigma$, or equivalently, how accurate the seller can estimate $\sigma$. 

\subsection{Worst-case market}\label{sec:wcmm}
For the ambiguity set \eqref{ambset}, we determine tight bounds for the tail probabilities, solving the minimization part \eqref{rr} of the central maximin problem \eqref{rmpp}. To this end, we introduce
\begin{flalign}\label{thresholdvalues}
    &&\uupsilon_1 \eqdef \mu - \frac{\usigma^2}{\beta-\mu}, \lupsilon_1 \eqdef \mu - \frac{\lsigma^2}{\beta-\mu}, \text{ and } \lupsilon_2 \eqdef \mu + \frac{\lsigma^2}{\mu},&&
\end{flalign}
which serve as threshold values for Lemma \ref{pdinf} and play a significant role in the subsequent analysis.

\begin{lemma}[Tight tail bounds]\label{pdinf}
Consider ambiguity set $\cP(\mu,\lsigma,\usigma, \beta)$. Then, 
    \begin{flalign}\label{ttb}
    &&\inf_{\bP \in \cP(\mu,\lsigma,\usigma, \beta)}\bP(X \geq p)=\inf_{\bP \in \cP(\mu,\lsigma,\usigma, \beta)}\bP(X > p)=
    \begin{cases}
     \frac{(\mu-p)^2}{(\mu-p)^2+\usigma^2}, \quad& p\in(0,\uupsilon_1], \\
     \frac{\mu-p}{\beta-p}, \quad& p\in[\uupsilon_1,\lupsilon_1], \\
     \frac{\mu^2+\lsigma^2-\mu  p}{\beta(\beta-p)},\quad& p\in[\lupsilon_1,\lupsilon_2], \\
     0, \quad& p\in[\lupsilon_2,\beta].
    \end{cases}&&
\end{flalign}
\end{lemma}
\begin{customproofsketch}
Lemma~\ref{pdinf} generalizes the one-sided Chebyshev bound, also known as Cantelli's inequality, for mean-variance information \citep{feller1991introduction1}, and also extends results in \cite{de1995general} for precise knowledge of the variance.
Let us sketch the proof of Lemma~\ref{pdinf}.  The key step is to formulate the tight bound as the solution to a semi-infinite linear program, see e.g.~\cite{shapiro2001duality}. Writing the tail probability $\bP(X > p)$ as the expectation of an indicator function, 
$\int_x \1{\{x> p\}}{\rm d} \mathbb{P}(x)$, we see that the tight bound follows from solving the primal problem (with $\mathcal{M}$ be the set of probability measures defined on the measurable space $([0,\beta],\mathcal{B})$) 
\begin{equation}\label{primal1x}
\begin{aligned}
&\! \inf_{\mathbb{P}\in \mathcal{M}} &  &\int_x \1{\{x> p\}}{\rm d} \mathbb{P}(x),\\
&\text{s.t.} &      &  \int_x {\rm d}\mathbb{P}(x)=1, \ \int_x x{\rm d}\mathbb{P}(x)=\mu,\ \int_x x^2{\rm d}\mathbb{P}(x)\geq\lsigma^2+\mu^2,\ \int_x x^2{\rm d}\mathbb{P}(x)\leq\usigma^2+\mu^2.
%&               &      & \int p(x){\rm d}x=1\\
%&               &      & p(x)\geq 0.\\
\end{aligned}
\end{equation}
The proof then continues with solving the corresponding dual problem 
\begin{equation}\label{dual1x}
\begin{aligned}
&\sup_{\lambda_0,\lambda_1,\lambda_2,\lambda_3} &  &\lambda_0 + \lambda_1 \mu+\lambda_2 (\lsigma^2+\mu^2) + \lambda_3 (\usigma^2+\mu^2),\\
&\text{s.t.} &      & \1{\{x> p\}}\geq \lambda_0  +\lambda_1 x+(\lambda_2+\lambda_3) x^2 =: F(x) \ \forall x\in[0,\beta], \lambda_2 \geq 0, \lambda_3 \leq 0,
\end{aligned}
\end{equation}
by considering all possible shapes of the minorant $F(x)$. For each shape we solve the dual and use complementary slackness to obtain a primal solution, showing that the primal and dual solutions agree. The full proof is presented in Section \ref{pdp}.
\end{customproofsketch}

From Lemma~\ref{pdinf} we see that 
the revenue maximization problem \eqref{rmpp} becomes
    $\sup_{p}\Tilde{R}(p)$
with 
\begin{flalign}\label{rev_g}
    &&\Tilde{R}(p) \eqdef R(p,\cP(\mu,\lsigma,\usigma,\beta)) = \begin{cases}
        R_1(p)=\frac{p(\mu-p)^2}{(\mu-p)^2+\usigma^2}, \quad & p \in (0,\lupsilon_1], \\
        R_2(p)=\frac{p(\mu-p)}{\beta-p}, \quad & p \in [\uupsilon_1,\lupsilon_1], \\
     R_3(p)=\frac{p(\mu(\mu-p)+\lsigma^2)}{\beta(\beta-p)}, \quad & p \in [\lupsilon_1,\lupsilon_2], \\
     0, \quad & p \in [\lupsilon_2,\beta].
    \end{cases}&&
\end{flalign}
\begin{figure}[t]
\centering
\subfigure[Worst-case demand]{
\label{}
{\includegraphics[width=0.45\linewidth]{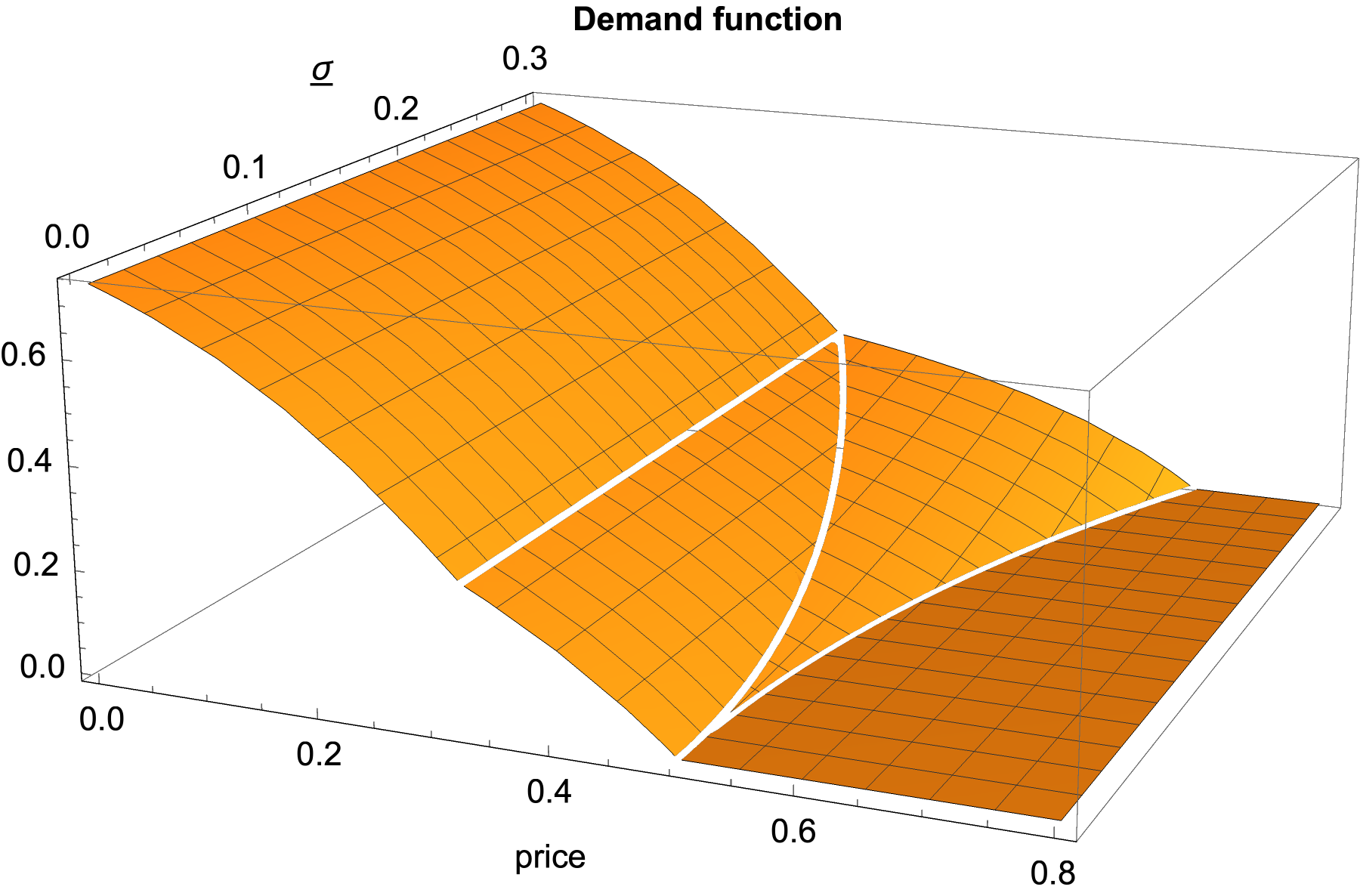}}
}
\subfigure[Worst-case revenue]{
\label{}
{\includegraphics[width=0.45\linewidth]{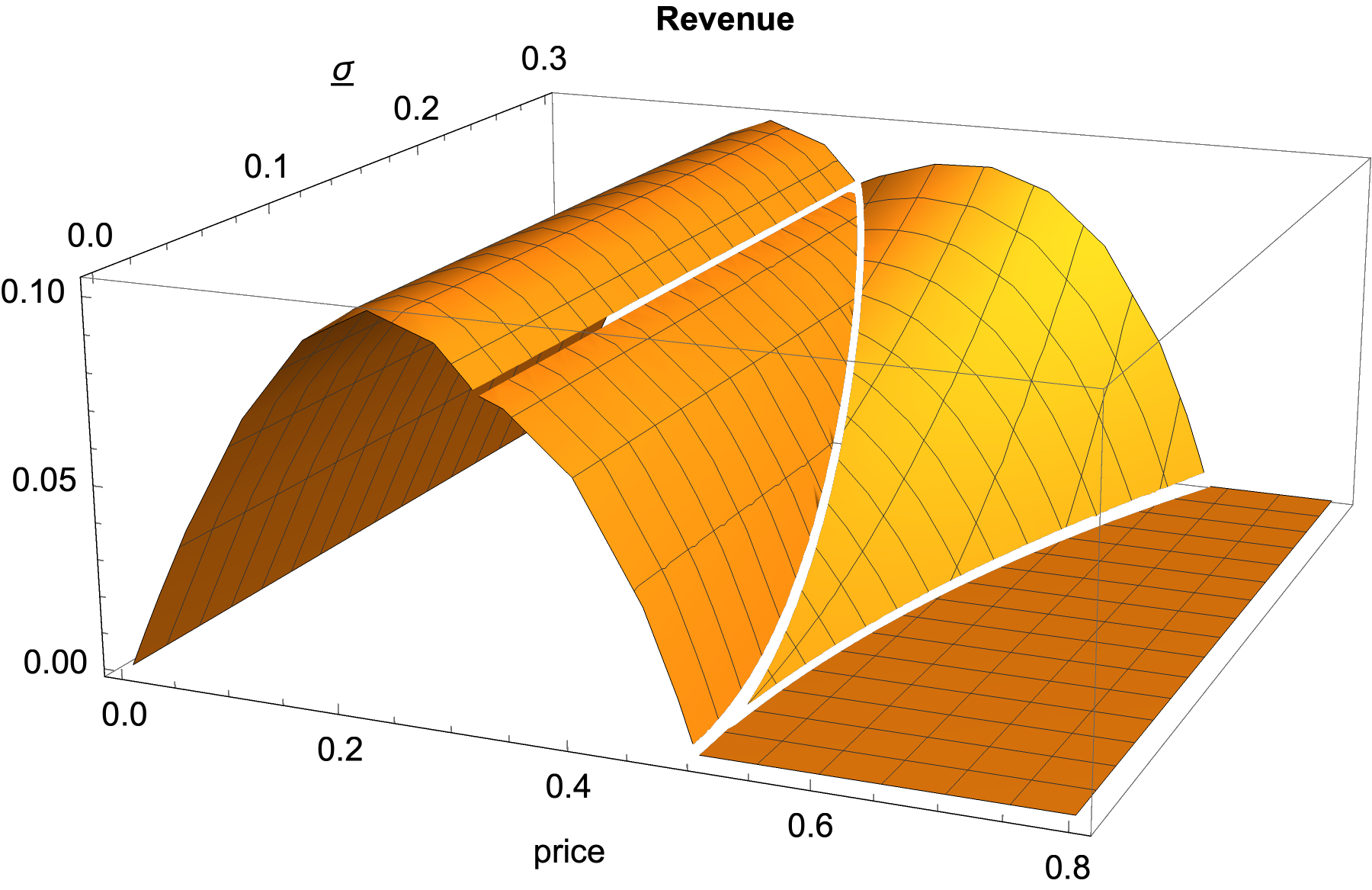}}
}
\caption{Worst-case demand function  \eqref{ttb}
    and revenue \eqref{rev_g} for $\mu=0.5$, $\beta=1$ and $\usigma=0.3$. 
    }\label{fig:1demrev}
\end{figure}
\noindent Let us now examine the structural properties of the worst-case demand function \eqref{ttb} and the corresponding worst-case revenue function \eqref{rev_g}. 
The worst-case demand function faced by the robust seller is derived in \eqref{ttb} by selecting, for every fixed price \( p \), the valuation distribution that minimizes expected demand. 
Although this worst-case valuation distribution is a two- or three-point distribution for each fixed \( p \), the resulting worst-case demand function is more intricate. 
To formalize this, we write  
\begin{flalign}
    &&
\inf_{\bP \in \cP(\mu, \lsigma, \usigma, \beta)} \bP(X \geq p) =: 1 - F^*(p), \quad \forall p \in [0, \beta].
&&
\end{flalign}
This representation effectively transforms the distributionally robust setting into the standard Bayesian setting with a known valuation distribution $F^*(p)$ defined via the tail bound \eqref{ttb}. 

\vspace{.2cm}
\noindent{\bf Lack of smoothness.}
One might then consider solving the classic optimization problem \( \sup_p p(1 - F^*(p)) \) using first-order conditions as in \eqref{foc}. 
However, this approach fails because the distribution \( F^*(p) \) does not satisfy the regularity conditions required for a convex revenue function or the existence of a unique optimal price.
Figure~\ref{fig:1demrev} illustrates \( 1 - F^*(p) \) and the resulting revenue function \( p(1 - F^*(p)) \) for \( \usigma = 0.3 \) across various pairs of \( (\mu,\lsigma, \usigma,\beta) \). 
The left panel reveals that the tail cumulative distribution function (CDF) is non-smooth at certain transition points. 
At these points, the distribution becomes non-differentiable, leading to non-convexity in the revenue function, as shown in the right panel. 
Consequently, the minimization problem \eqref{rr} results in a maximization problem \eqref{rmpp} that is non-convex and admits multiple local optima. 

\vspace{.2cm}
\noindent{\bf Segmented market.} Despite the lack of smoothness, the worst-case demand function exhibits a structural shape that aligns intuitively with robust pricing. 
A relatively large \( \usigma \) may arise from demand functions associated with heterogeneous customer segments, such as those differing in income, age, or geography. 
However, the seller possesses only coarse information about market uncertainty in the form of \( (\lsigma, \usigma) \). 
If the seller had precise knowledge of distinct market segments and could implement segment-specific pricing, third-degree price discrimination might be an option. 
Here, we do not consider that possibility; instead, the seller must determine the optimal uniform price.
Interestingly, the structure of the worst-case demand function inferred from our analysis suggests a segmented market. 
The left panel of Figure~\ref{fig:1demrev} and Lemma~\ref{pdinf} indicate that the worst-case market can be interpreted as consisting of three or four segments: consumers willing to pay nothing or prices within intervals \( (0, \uupsilon_1) \), \( (\uupsilon_1, \lupsilon_1) \), and \( (\uupsilon_1, \lupsilon_2) \). 
The non-smoothness of \( 1 - F^*(p) \) at points like \( p = \uupsilon_1 \) suggests that as the price surpasses this threshold, a specific segment ceases to purchase the product.

\vspace{.2cm}
\noindent{\bf Rise of the niche segment.} It is crucial to emphasize that we did not assume a priori the existence of multiple customer groups or a segmented market. 
Rather, the segmentation emerges as a natural interpretation of the worst-case demand function, with the degree of segmentation increasing as market risk, represented by \( \usigma \), grows. 
This phenomenon is evident in Figure~\ref{fig:1demrev}, where increasing \( \lsigma \) amplifies the dominance of the third segment \( (\uupsilon_1, \lupsilon_2) \). 
Notably, for \( \lsigma = \usigma \), the revenue function exhibits two nearly comparable local maxima. 
As market risk intensifies, the maximum from the higher segment \( (\uupsilon_1, \lupsilon_2) \) approaches or even surpasses the maximum from the lower segment \( (0, \uupsilon_1) \), ultimately delivering the global maximum. 
In Appendix~\ref{9plotsexamples}, we present figures for \( \usigma \) ranging from 0 to 0.45. 
These illustrate that for \( \usigma \geq 0.35 \), the higher segment's local maximum overtakes the lower segment's maximum. 
This initial insight—that the robust seller may shift focus from the lower segment to the higher segment in response to increased market risk—will be further substantiated below.

\subsection{Main results on robust pricing}
We will now present this paper's main result for robust monopoly pricing. To do so, we introduce three price functions:
\begin{flalign}\label{price_h}
      &&p^*_{h}(\lsigma,\beta) &= \beta - \sqrt{\beta\left(\beta-\mu-\frac{\lsigma^2}{\mu}\right)},&&\\ \label{price_m}
     &&p^*_{m}(\beta) &= \beta-\sqrt{\beta(\beta-\mu)},&&\\ \label{price_l}
      &&p^*_{l}(\usigma) &= \mu - \usigma \kappa(\mu,\usigma),&&
\end{flalign}
with
\begin{flalign}\label{kappa}
&&\kappa(\mu,\usigma) &\eqdef \left(\tfrac{\mu}{\usigma}+\sqrt{1+(\tfrac{\mu}{\usigma})^2}\right)^{\frac{1}{3}} + \left(\tfrac{\mu}{\usigma}-\sqrt{1+(\tfrac{\mu}{\usigma})^2}\right)^{\frac{1}{3}}.&&
\end{flalign}
The arguments of these price functions are chosen deliberately to express the crucial fact that the price  $p^*_{l}(\usigma)$ responds to the risk level $\usigma$, while the price $p^*_{m}(\beta)$ does not. While all three price functions depend on $\mu$, it is always considered fixed and thus omitted as an argument. Notice that unlike $p^*_{l}(\usigma)$ and $p^*_{m}(\beta)$ the price function $p^*_{h}(\lsigma,\beta)$ depends on both $\lsigma$ and $\beta$. This price function turns out pivotal in high-risk settings.

\begin{figure}[t!]
\begin{center}
\begin{tikzpicture}
    \node at (0,0) {\includegraphics[width=0.5\linewidth]{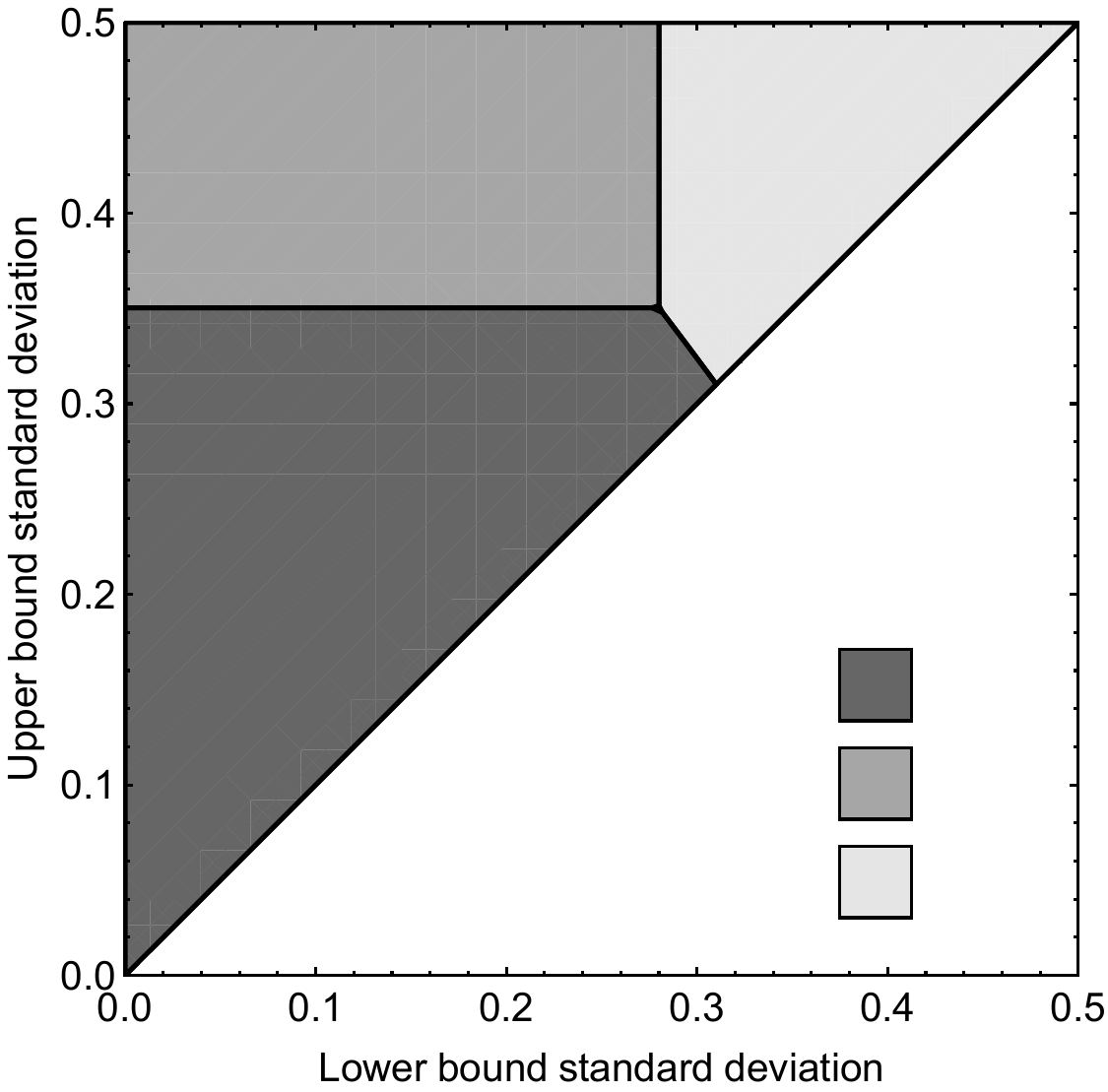}};
     \node at (2.75, -1)    
    {$\Sigma_l$};
    \node at (2.82, -1.67)    
    {$\Sigma_m$};
    \node at (2.78, -2.32)    
    {$\Sigma_h$};
    \end{tikzpicture}
\caption{Three prices theorem for $\mu=0.5$ and $\beta=1$.}\label{Figmain3price}
\end{center}
\end{figure}

\begin{theorem}[Three-Prices Theorem]\label{3pthm}
When the seller knows $\mu$, $\beta$, and that $\sigma\in [\lsigma, \usigma]$, the optimal robust price $p^* = \arg\sup_p\Tilde{R}(p)$ is contained in the set $\{p^*_l(\usigma),p^*_m(\beta),p^*_h(\lsigma,\beta)\}$. 
\begin{flalign}\label{p^*_gmain}
    &&p^* = 
    \begin{cases}
    p^*_{l}(\usigma) \quad &  {\rm if} \ (\lsigma,\usigma)   \in \Sigma_l, \\
    p^*_{m}(\beta) \quad &  {\rm if} \ (\lsigma,\usigma) \in \Sigma_m,
    \\
    p^*_{h}(\lsigma,\beta) \quad &  {\rm if} \ (\lsigma,\usigma) \in \Sigma_h,
    \end{cases}&&
\end{flalign}
with $\Sigma_l$, $\Sigma_m$ and $\Sigma_h$ speficied in 
\eqref{sets} and such that
$\Sigma_l\cup\Sigma_m\cup\Sigma_h$ spanning all pairs $(\lsigma,\usigma)$ with  $\lsigma \leq \usigma$.
\end{theorem} 
See Figure \ref{Figmain3price} for a visual representation of \eqref{sets}.

\vspace{.2cm}
\noindent {\bf Straight-line boundaries.} Figure~\ref{Figmain3price} visualizes Theorem~\ref{3pthm}. The 2D map, divided by straight lines, represents regions corresponding to optimal prices, determined by the market risk bounds—lower and upper—which influence the worst-case demand function \eqref{ttb}. For low prices, the demand depends only on the upper bound $\usigma$, while for high prices, it is governed by the lower bound $\lsigma$. For medium prices, the demand becomes independent of both bounds, making neither significant for determining the optimal price. The boundaries between these regions reflect transitions in how the demand depends on the bounds, whether shifting from one to the other or to neither. This resulting piecewise linear partitioning of the 2D space creates a clear and interpretable map for understanding how market risk bounds influence robust revenue. The oblique line separating the regions corresponding to low and high prices is neither horizontal nor vertical, as each region is determined by a different bound of the standard deviation.

\vspace{.2cm}
\noindent {\bf General price jump.}
When either $\lsigma$ or $\usigma$ is held constant and does not coincide with the intersection of $\Sigma_l$ and $\Sigma_m$, varying the other parameter will result in at most two relevant prices, but never three. The high pricing strategy $p^*_h(\lsigma,\beta)$ becomes optimal when $\lsigma$ and $\usigma$ are both sufficiently large. This insight can be described more formally as follows:

\begin{theorem}[Price-Jump Theorem]\label{general_jump_in_words}
When the upper bound $\usigma$ increases from the lower bound $\lsigma$ to the maximal value, the optimal price jumps from the low price $p^*_l(\usigma)$ to a higher price, which is either $p^*_m(\beta)$ or $p^*_h(\lsigma,\beta)$, given that $\lsigma$ is sufficiently small.
\end{theorem}
Theorem~\ref{general_jump_in_words} not only identifies the three candidate prices corresponding to the three distinct variance regions but also reveals a critical insight: as demand uncertainty increases, the robust seller strategically introduces discontinuous price jumps to adapt to changing conditions. This insight into the occurrence of discontinuous price adjustments forms the basis of what we term the low-high pricing phenomenon, which serves as the central theme of this paper.

The proofs of Theorems~\ref{3pthm}
and \ref{general_jump_in_words} are provided in Section~\ref{sec:twoprices}. 
These proofs build on results for two special cases: knowing the precise value of the variance and knowing the upper bound on the support (but not both). 
Knowing the true variance corresponds with the diagional in Figure~\ref{Figmain3price}, where the lower and upper bound on $\sigma$ agree. The discontinuous price jump for this case can also be seen in Figure~\ref{fig:beta}. Knowing only the upper bound corresponds to the vertical axis in Figure~\ref{Figmain3price} with a lower bound of zero. Both cases are analyzed in great detail in Section~\ref{sec:twoprices}.

\newpage
\noindent  {\bf Ground truth markets.}
A natural benchmark for assessing robust pricing is comparison with strategies based on full empirical knowledge of the valuation distribution. If abundant and reliable data are available, empirical optimization will, by definition, achieve higher performance than any method that only uses partial summaries. However, in practice, such data are rarely available, and fitted classical distributions are often used as surrogates. When we examine truncated normal, exponential, lognormal, gamma, and beta distributions (see Appendix \ref{app:gtms} for some examples), we find that they all yield smooth revenue functions. Truncated classical families produce unimodal prices that decrease gradually with variance, while beta distributions may exhibit degenerate bimodality under high variance, but without meaningful segmentation of consumers. In short, while actual empirical markets may display segmentation if it truly exists, fitting classical smooth distributions enforces unimodality and cannot capture the sharp transitions or the low–high pricing behavior that robust models reveal.

\vspace{.2cm}
\noindent {\bf Added value of robustness.}
The strength of robustness lies in its ability to work with limited information while avoiding overly pessimistic prescriptions. When only mean and variance are known, robust models produce unimodal and conservative prices. Once range information is included, however, the worst-case distribution becomes genuinely bimodal, splitting consumers into low- and high-valuation segments. This latent segmentation creates a bimodal revenue function and explains the discontinuous switching between low and high optimal prices. Crucially, because the combination of mean, variance, and range already provides substantial information, the resulting solutions are not excessively conservative. Instead, they uncover the appropriate pricing mechanism that fitted smooth distributions obscure. Robustness therefore serves not only as protection against missing or noisy data, but also as a constructive modeling tool that highlights underlying market structures.

\section{Proofs of robust pricing with two and three prices}\label{sec:twoprices}
This section serves a dual purpose: proving the main pricing theorems and developing further robust pricing insights. We begin by examining two scenarios: one in which the seller knows the exact value (Section~\ref{secupiv}) and another where the seller is aware only of the maximum value of the standard deviation (Section~\ref{secupiv2}). In both cases, we show that the optimal robust pricing strategy involves selecting one of two prices, as formalized in Propositions \ref{theorem_s} and \ref{theorem_h}. These results are instrumental in extending the analysis to the more general setting where the standard deviation lies within a specified interval. The optimal strategy involves three candidate prices, as established in Theorems~\ref{3pthm} and \ref{general_jump_in_words}, proven in Sections~\ref{ProofT1} and \ref{ProofT2}.

\subsection{Upper bound on standard deviation}\label{secupiv}
 First assume that the seller knows the mean valuation $\mu$, the maximal valuation $\beta$, and that the standard deviation of the valuation distribution is at most $\usigma$. The robust price should then solve the maximin problem 
%$\begin{equation}\label{rmpp_s}
    %\textup{REV}(\mu,\usigma^2,\beta) %\eqdef %\sup_{p}\inf_{\bP\in\cP_{(\mu,\usigma,\beta)}}\bP(X \geq p)=
    $\sup_{p}\Bar{R}(p)$, 
 where the function $\Bar{R}(p)$ follows from Lemma~\ref{pdinf} for $\lsigma=0$:
\begin{flalign}\label{rev_s}
&&\Bar{R}(p) \eqdef R(p,\cP(\mu,\usigma,\beta)) = \begin{cases}
        R_1(p)\eqdef\frac{p(\mu-p)^2}{(\mu-p)^2+\usigma^2}, \quad & p \in (0,\uupsilon_1], \\
    R_2(p)\eqdef\frac{p(\mu-p)}{\beta-p}, \quad & p \in [\uupsilon_1,\mu], \\
     0, \quad & p \in [\mu,\beta],
    \end{cases}&&
\end{flalign}
with $\uupsilon_1$ as in \eqref{thresholdvalues} and $\cP(\mu,\usigma,\beta)$ the set of all valuation distribution with mean $\mu$, upper bound on valuation $\beta$, and upper bound on standard deviation $\usigma$.
In this case, the worst-case revenue $\Bar{R}(p)$ can thus be viewed as the maximum of  $R_1(p)$, $R_2(p)$ and $0$, so the optimal robust price will follow from either $R_1(p)$ or $R_2(p)$. This is illustrated in Figure~\ref{figmain1}. Observe that $\Bar{R}(p)$ indeed consists of the two concatenated functions $R_1(p)$, $R_2(p)$, and that both functions might create a local maximum. For $\bar\sigma=0.36$ there are two local maxima, corresponding to two prices, and the global maximum follows from the largest price. As it turns out, the local maximum delivered by $R_1(p)$, associated with a relatively low price, is globally optimal when the risk level $\bar{\sigma}$ is sufficiently small. And vice versa, when $\bar{\sigma}$ is sufficiently large, the optimal price is relatively high and corresponds to the argmax of $R_2(p)$.

\begin{figure}[ht]
\centering
\subfigure[$\bar{R}(p)$ with $\mu=0.5$ and $\beta=1$.]{
\label{figmain1}
\begin{tikzpicture}
    \node at (0,0) {\includegraphics[width=0.46\linewidth]{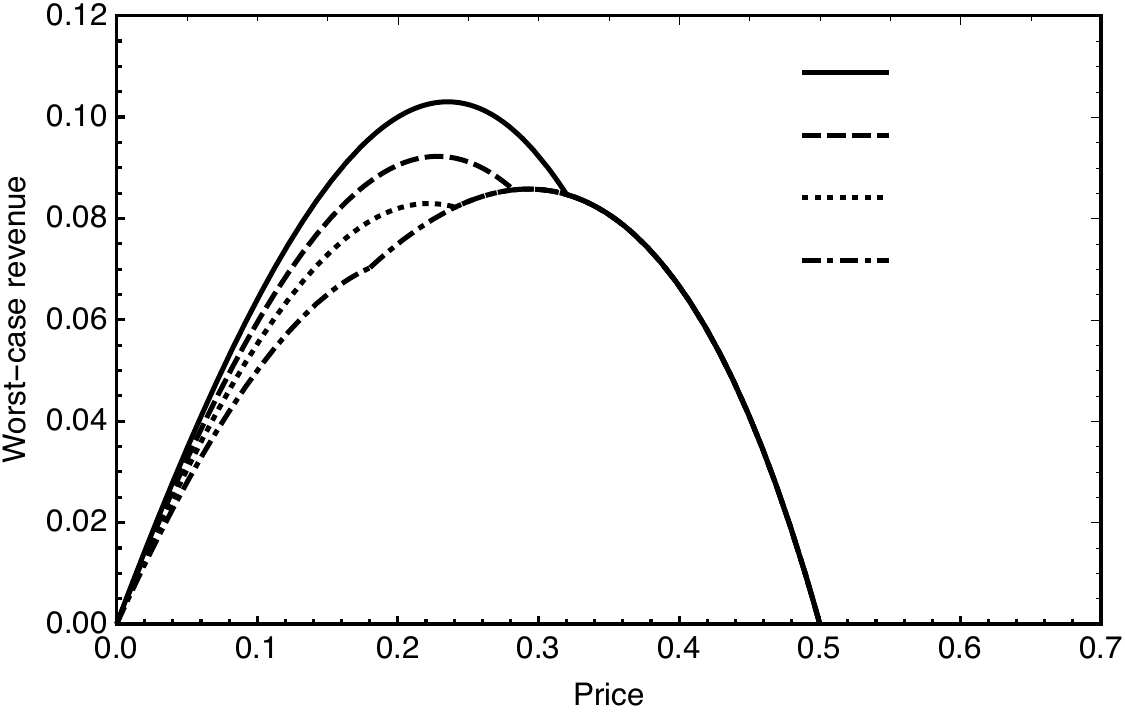}};
    \node at (2.75, 1.48)   {\begin{scriptsize}$\bar \sigma=0.33$\end{scriptsize}};
     \node at (2.75, 1.88)   {\begin{scriptsize}$\bar \sigma=0.30$\end{scriptsize}};
     \node at (2.75, 1.06)   {\begin{scriptsize}$\bar \sigma=0.36$\end{scriptsize}};
     \node at (2.75, 0.65)   {\begin{scriptsize}$\bar \sigma=0.40$\end{scriptsize}};
    \end{tikzpicture}
}
\subfigure[Proposition~\ref{theorem_s} with $\mu=0.5$ and $\beta=1$.]{
\label{figmain1}
\begin{tikzpicture}
    \node at (0,0) {\includegraphics[width=0.46\linewidth]{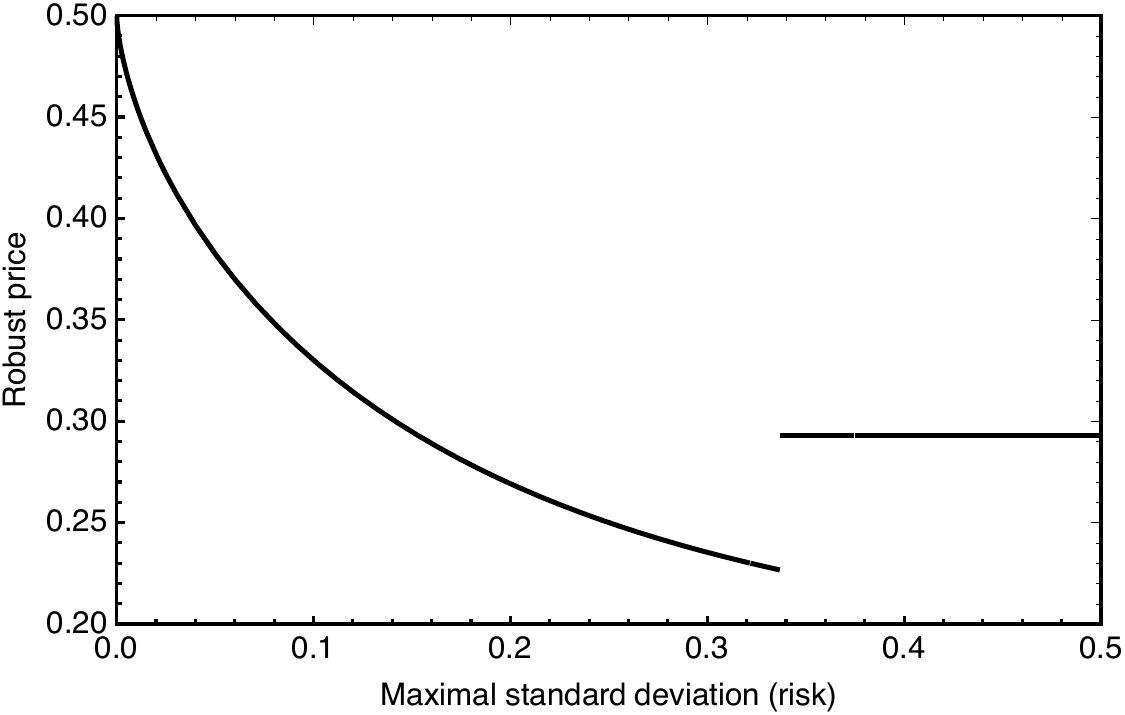}};
     \node at (2.3, 0.4)   {\begin{scriptsize}risk-insensitive\end{scriptsize}};
     \node at (2.3, 0)   {\begin{scriptsize} price $p^*_{m}(\beta)$\end{scriptsize}};
      \node at (-1, 0.8)   {\begin{scriptsize}risk-sensitive\end{scriptsize}};
     \node at (-1, 0.4)   {\begin{scriptsize} price $p^*_{l}(\bar\sigma)$\end{scriptsize}};
    \end{tikzpicture}
}
\caption{Illustration of the first main pricing theorem. }\label{Figmain}
\end{figure}

Here is then a pricing result, saying that for a large enough maximal risk-level $\bar{\sigma}$, the optimal price becomes risk-insensitive:
\begin{proposition}[Risk-sensitive and risk-insensitive pricing]\label{theorem_s}
The optimal robust price $\bar p^*=\arg\sup_{p}\Bar{R}(p)$ with $\Bar{R}(p)$ in  \eqref{rev_s}, and $
\bar{\sigma}^*=\left(\tfrac{32}{27}(\beta-\mu)\left(\sqrt{\beta(\beta-\mu)}-(\beta-\mu)\right)\right)^{\frac{1}{2}}$, is
\begin{flalign}
    &&\bar{p}^* = 
    \begin{cases}
    p^*_{l}(\bar\sigma) \quad &  {\rm if} \ \bar{\sigma} \leq \bar{\sigma}^*, \\
    p^*_{m}(\beta) \quad &  {\rm if} \ \bar{\sigma} \geq \bar{\sigma}^*,
    \end{cases}&&
\end{flalign}
\end{proposition}

\begin{customproof}  
The derivative of $R_1(p)$ equals zero at $p^*_l(\usigma)$ and $\mu$, but only when $p^*_l(\usigma)$ lies in the interval $(0,\uupsilon_1]$ is this the maximum ($\mu$ is never the maximum). Since $R_1(p)$ increases in  $[0,p^*_l(\usigma)]$, the maximum is attained in $\uupsilon_1$ when $\uupsilon_1=\min\{p^*_l(\usigma),\uupsilon_1\}$. Similarly, $R_2(p)$ is concave and attains 
a unique maximum at $p^*_m(\beta)\in [0,\mu]$. Since $R_2(p)$ increases in $[0,p^*_m(\beta)]$, the maximum of $R_2(p)$ in $[\uupsilon_1,\mu]$ is attained in $\uupsilon_1$ when $\uupsilon_1 = \max\{ p^*_m(\beta),\uupsilon_1\}$. 

We next show that, despite $\uupsilon_1$ being a candidate local maximizer, it can never be the global maximizer for the function $\Bar{R}(p)$. Observe that
$R(p,\cP(\mu,\usigma,\beta)) = R(p,\cP(\mu,\usigma) \cap \cP(\mu,\beta))
        \geq \max\{R(p,\cP(\mu,\usigma)),R(p,\cP(\mu,\beta))\}$
and
\begin{align*}
    R(p,\cP(\mu,\usigma,\beta)) &= R(p,\cP(\mu,\usigma))\1\{p\in (0,\uupsilon_1]\}+R(p,\cP(\mu,\beta))\1\{p\in (\uupsilon_1,\beta]\}\\
    &\leq \max\{R(p,\cP(\mu,\usigma)),R(p,\cP(\mu,\beta))\}.
\end{align*}
Hence, $R(p,\cP(\mu,\usigma,\beta)) = \max\{R(p,\cP(\mu,\sigma)),R(p,\cP(\mu,\beta))\}$. Consequently, since the maximization problem is equivalent to the maximization problem without either $\usigma$ or $\beta$, we have $\bar{p}^* \in \{p^*_l(\usigma),p^*_m(\beta)\}$. 
Since $R_1$  is continuously differentiable and attains its maximum in $p^*_l(\usigma)$, we can employ the envelope theorem
    to show that
    \begin{flalign}
        &&\frac{\partial R_{1}(p^*_l(\usigma))}{\partial\usigma^2} = \frac{\partial R_{1}(p)}{\partial\usigma^2}\biggr\rvert_{p = p^*_l(\usigma)} = -\frac{p^*_l(\usigma)(\mu-p^*_l(\usigma))^2}{((\mu-p^*_l(\usigma))^2+\usigma^2)^2} \leq 0,&&
    \end{flalign}
    since $p^*_{l}(\usigma) \geq 0$. Hence, $R_1(p^*_l(\usigma))$ is non-increasing in $\usigma$, so that $\usigma^*$ is the unique value of $\usigma$ that satisfies $R_1(p^*_l(\usigma)) = R_2(p^*_m(\beta))$, which proves the claim.
    \end{customproof}

One important conclusion from Proposition~\ref{theorem_s} is that pricing becomes risk-insensitive when the upper bound on standard deviation is sufficiently large; see Figure~\ref{figmain1}.
However, this will change drastically when the seller has better or more precise standard deviation information. 

\subsection{Precise knowledge of standard deviation}\label{secupiv2}
\begin{figure}[ht!]
\centering
\subfigure[$\hat{R}(p)$ with $\mu=0.5$ and $\beta=1$]{
\label{figmain1h}
\begin{tikzpicture}
    \node at (0,0) {\includegraphics[width=0.46\linewidth]{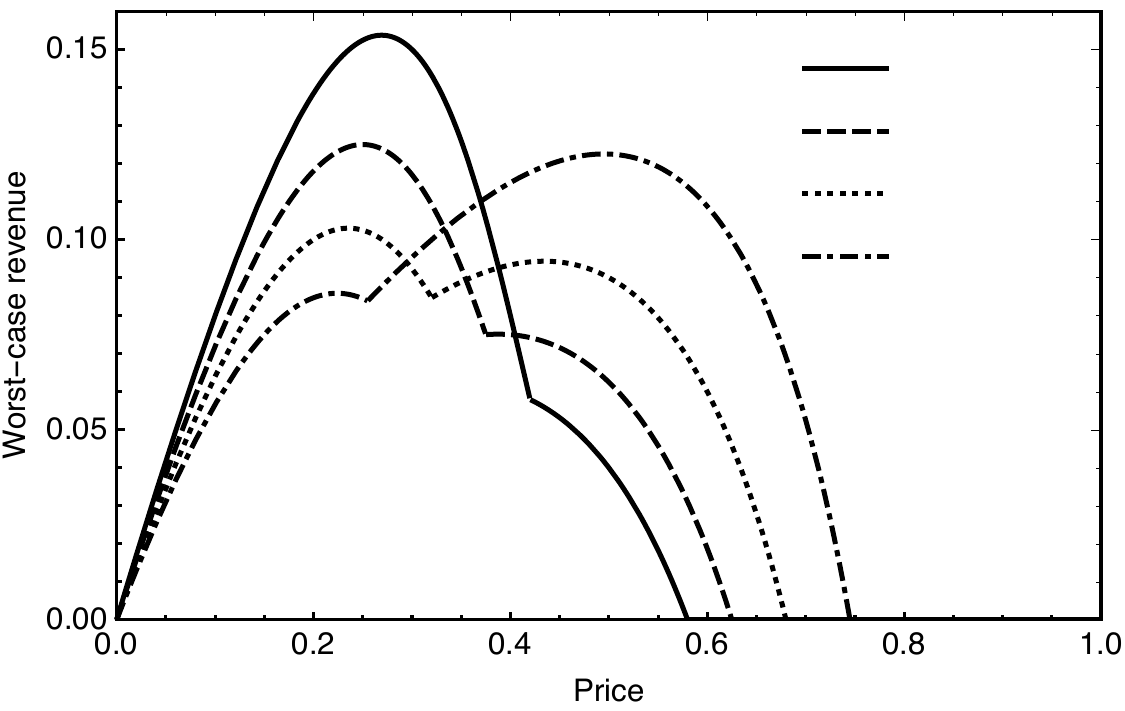}};
    \node at (2.75, 1.49)   {\begin{scriptsize}$ \sigma=0.20$\end{scriptsize}};
     \node at (2.75, 1.88)   {\begin{scriptsize}$ \sigma=0.25$\end{scriptsize}};
     \node at (2.75, 1.08)   {\begin{scriptsize}$ \sigma=0.30$\end{scriptsize}};
     \node at (2.75, 0.66)   {\begin{scriptsize}$ \sigma=0.35$\end{scriptsize}};
    \end{tikzpicture}
}
\subfigure[Proposition~\ref{theorem_h} with $\mu=0.5$ and $\beta=1$.]{
\label{figmain1hopt}
\begin{tikzpicture}
    \node at (0,0) {\includegraphics[width=0.46\linewidth]{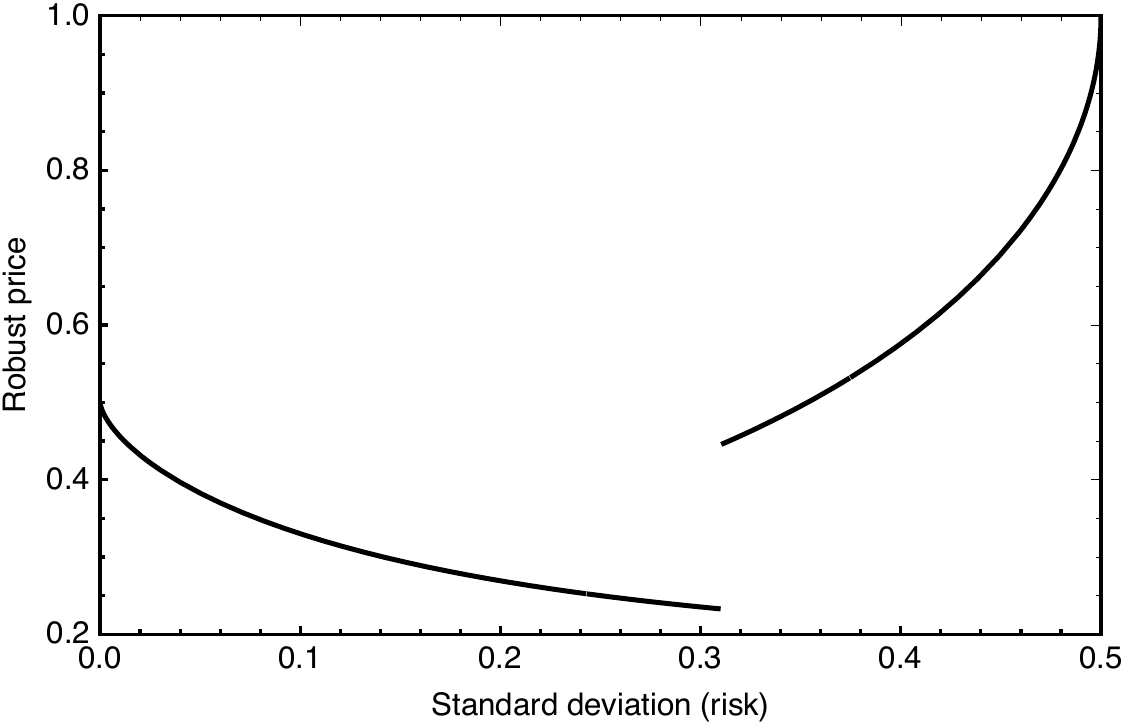}};
     \node at (2, 1.1)   {\begin{scriptsize}high price\end{scriptsize}};
     \node at (2, 0.7)   {\begin{scriptsize} $p^*_{h}(\sigma,\beta)$\end{scriptsize}};
      \node at (-1, -0.4)   {\begin{scriptsize}low price\end{scriptsize}};
     \node at (-1, -0.8)   {\begin{scriptsize} $p^*_{l}(\sigma)$\end{scriptsize}};
    \end{tikzpicture}
}
\caption{Illustration of the second main pricing theorem. }\label{Figmainh}
\end{figure}

Next, consider the setting where the seller knows the precise standard deviation which means $\lsigma = \usigma$. For notational convenience, we will drop the bar notation and write $\sigma$ instead. We solve
$\sup_{p}\hat{R}(p)$, where it follows from Lemma~\ref{pdinf} for $\lsigma = \usigma$ that
\begin{flalign}\label{rev_h}
    &&\hat{R}(p) \eqdef R(p,\cP(\mu,\sigma,\beta)) = \begin{cases}
        R_1(p)=\frac{p(\mu-p)^2}{(\mu-p)^2+\sigma^2}, \quad & p \in (0,\lupsilon_1], \\
     R_3(p)\eqdef\frac{p(\mu(\mu-p)+\sigma^2)}{\beta(\beta-p)}, \quad & p \in [\lupsilon_1,\lupsilon_2], \\
     0, \quad & p \in [\lupsilon_2,\beta],
    \end{cases}&&
\end{flalign}\\ 
with $\lupsilon_1$ and $\lupsilon_2$ as in \eqref{thresholdvalues} and $\cP(\mu,\sigma,\beta)$ the set of all valuation distribution with mean $\mu$, upper bound on valuation $\beta$, and standard deviation $\sigma$. Figure~\ref{figmain1h} displays the worst-case revenue for several risk levels. Observe that $\hat{R}(p)$ still consists of the two quadratic functions, pasted together, and that both functions present a candidate for the global maximum. However, compared with Figure~\ref{figmain1}, there are some notable differences. First, the price range with non-degenerate prices increases with $\sigma$, now allowing robust prices that exceed the mean $\mu$. Second, the risk level from which the high price starts dominating the low price seems to be lower than in Figure~\ref{figmain1}. And third, the plots in Figure~\ref{figmain1h} show that for high prices, the revenue increases significantly with $\sigma$, which leads to premium pricing targeted at users with the highest valuations.

Taken together, this gives the following pricing result:
\begin{proposition}[Low and high pricing]\label{theorem_h}
    The optimal robust price $ \hat{p}^*=\arg\sup_{p}{\hat{R}}(p)$ with ${\hat{R}}(p)$ in \eqref{rev_h} is
    \begin{flalign}
        &&\hat{p}^* = \begin{cases}
            p^*_l(\sigma) \quad & {\rm if}\  \sigma \leq \sigma^*, \\
            p^*_h(\sigma,\beta) \quad & {\rm if}\    \sigma \geq \sigma^*,
        \end{cases}&&
    \end{flalign}
    with $\sigma^*$ defined implicitly as the solution to $R_1(p^*_{l}(\sigma)) = R_3(p^*_{h}(\sigma,\beta))$.
\end{proposition}

Figure~\ref{figmain1hopt}
illustrates Proposition~\ref{theorem_h}. In response to increased risk, the robust seller initially lowers the price, aligned with classical pricing theory. When the risk surpasses some threshold, however, the seller changes its strategy drastically and switches to a high price that increases as a function of the risk. This high price function typically exceeds the mean valuation and converges to the maximal valuation when the risk approaches its maximal value. 

Proving  Proposition~\ref{theorem_h} is considerably more challenging than proving Proposition~\ref{theorem_s}.   
Our proof starts with the following result:
\begin{lemma}\label{lemma_h_1}
    Consider $ \hat{p}^*=\arg\sup_{p}{\hat{R}}(p)$ with ${\hat{R}}(p)$ in \eqref{rev_h}. The optimal price $\hat{p}^*$ is $p^*_l(\sigma)$, $p^*_h(\sigma,\beta)$ or $\lupsilon_1$. 
\end{lemma}
\begin{customproof}
From Proposition \ref{theorem_s} we know that the maximum of $R_1(p)$ is either attained in $p^*_l(\sigma)$ or in $\lupsilon_1$. Also, as $\sigma^2$ is bounded from above by $\mu(\beta-\mu)$, 
\begin{flalign}
    &&R_3''(p) = -\frac{2(\mu(\beta-\mu)-\sigma^2)}{(\beta-p)^3} \leq 0.&&
\end{flalign}
 Hence, $R_3(p)$ is concave and attains a unique maximum at $p^*_h(\sigma,\beta) \in [\lupsilon_1,\lupsilon_2]$. Since $R_3(p)$ decreases in $[p^*_h(\sigma,\beta),\lupsilon_2]$, the maximum of $R_3(p)$ in $[\lupsilon_1,\lupsilon_2]$ is attained in $\lupsilon_1$ when $\lupsilon_1 = \max\{p^*_h(\sigma,\beta),\lupsilon_1\}$.
\end{customproof}

As in the proof of Proposition~\ref{theorem_s}, we want to rule out $\lupsilon_1$ as candidate optimal price, but this now turns out more challenging. The next lemma states that indeed $\lupsilon_1$ can never be the optimal price, and the proof using some advanced analysis is deferred to the appendix. 
\begin{lemma}\label{middlepointnot}
Consider $ \hat{p}^*=\arg\sup_{p}{\hat{R}}(p)$ with ${\hat{R}}(p)$ in \eqref{rev_h}. The optimal price $\hat{p}^*$  can never be $\lupsilon_1$. 
\end{lemma}
Lemmas~\ref{lemma_h_1} and \ref{middlepointnot} together show that the optimal price is either $p^*_l(\sigma)$ or $p^*_h(\sigma,\beta)$. The remaining step is to show that the low (high) price is optimal for low (high) standard deviation. 
\begin{lemma}\label{lemma_h_2}
Consider $\sigma^*$ defined in {\rm Proposition~\ref{theorem_h}}. Then, $\hat{R}(p^*_l(\sigma)) > \hat{R}(p^*_h(\sigma,\beta))$ for $\sigma < \sigma^*$,    $\hat{R}(p^*_l(\sigma)) < \hat{R}(p^*_h(\sigma,\beta))$ for $\sigma > \sigma^*$  and  $\hat{R}(p^*_l(\sigma)) = \hat{R}(p^*_h(\sigma,\beta))$ for $\sigma = \sigma^*$. 
\end{lemma}
\begin{customproof}
    {\noindent Since $R_1(p)$ is maximal and continuously differentiable in $p^*_l(\sigma)$, we can employ the envelope theorem to show that
    \begin{flalign}
        &&\frac{\partial R_{1}(p^*_l(\sigma))}{\partial\sigma^2} = \frac{\partial R_{1}(p)}{\partial\sigma^2}\biggr\rvert_{p = p^*_l(\sigma)} = -\frac{p^*_l(\sigma)(\mu-p^*_l(\sigma))^2}{((\mu-p^*_l(\sigma))^2+\sigma^2)^2} \leq 0,&&
    \end{flalign}
    since $p^*_{l}(\sigma) \geq 0$. Similarly, since $R_3(p)$ is maximal and continuously differentiable in $p^*_h(\sigma,\beta)$, the envelope theorem yields
    \begin{flalign}
        &&\frac{\partial R_{3}(p^*_h(\sigma,\beta))}{\partial\sigma^2} = \frac{\partial R_{3}(p)}{\partial\sigma^2}\biggr\rvert_{p = p^*_h(\sigma,\beta)} = \frac{p^*_h(\sigma,\beta)}{\beta(\beta-p^*_h(\sigma,\beta))} \geq 0,&&
    \end{flalign}
    since $p^*_h(\sigma,\beta) \geq 0$. Together these two statements prove that the threshold value $\sigma^*$ exists and is unique.}
    
    Next, consider $\sigma \xrightarrow{}0$ and assume $p^*_h(\sigma,\beta) \in [\lupsilon_1,\lupsilon_2]$. Then $p^*_h(\sigma,\beta) = \mu$, which is a contradiction, and $\hat{p}^* = p^*_l(\sigma)$. Furthermore, consider $\sigma \xrightarrow{}\sigma_{\textup{max}} = \sqrt{\mu(\beta-\mu)}$ and assume $p^*_l(\sigma) \in [0,\lupsilon_1]$. Then $p^*_l(\sigma) = 0$, which cannot be optimal, so $\hat{p}^* = p^*_h(\sigma,\beta)$.
\end{customproof}

Proposition~\ref{theorem_h} then follows from Lemmas~\ref{lemma_h_1},  \ref{middlepointnot} and \ref{lemma_h_2}.

\subsection{Proof of Theorem~\ref{3pthm}}\label{ProofT1}
%\begin{customproof}
    When $\lsigma = \usigma$, we know from Proposition~\ref{theorem_h} that $p^* \in \{p^*_l(\usigma),p^*_h(\lsigma,\beta)\}$. For the remainder of the proof, assume $\lsigma < \usigma$. From Proposition \ref{theorem_s} and \ref{theorem_h}, we know that $p^*_l(\usigma),p^*_m(\beta),p^*_h(\lsigma,\beta)$ are the optimal prices of the unconstrained functions $R_1(p),R_2(p),R_3(p)$ respectively. What remains is to exclude the potential boundary optimal prices $\uupsilon_1$ and $\lupsilon_1$. 

 We will now show that $\uupsilon_1$ cannot be the unique optimal price. Since $\uupsilon_1 \in (0,\lupsilon_1)$, we can assume without loss of generality that $p^* < \lupsilon_1$ (since $p^* \geq \lupsilon_1$ implies $p^*\neq \uupsilon_1$). The following holds due to Proposition \ref{theorem_s}:
    \begin{flalign}
        &&R(p,\cP(\mu,\lsigma,\usigma,\beta)) &\geq R(p,\cP(\mu,\usigma,\beta))
        = \max\{R(p,\cP(\mu,\usigma)),R(p,\cP(\mu,\beta))\}.&&
    \end{flalign}
    In addition, it holds that
    \begin{flalign}
&&R(p,\cP(\mu,\lsigma,\usigma,\beta)) &=
        R(p,\cP(\mu,\usigma))\1\{p\in (0,\uupsilon_1]\}+R(p,\cP(\mu,\beta))\1\{p\in (\uupsilon_1,\lupsilon_1)\},&&\nonumber \\
        && &\leq \max\{R(p,\cP(\mu,\usigma)),R(p,\cP(\mu,\beta))\}.&&
    \end{flalign}
    \noindent Hence, $R(p,\cP(\mu,\lsigma,\usigma,\beta)) = \max\{R(p,\cP(\mu,\usigma)),R(p,\cP(\mu,\beta))\}$ and therefore, $p^*$ is either $p^*_l(\usigma)$ or $p^*_m(\beta)$ whenever $p^* \in [0,\lupsilon_1)$, since the maximization problem is equivalent to the maximization problem without either $\usigma$ or $\beta$. This means we can exclude $\uupsilon_1$ as an optimal price candidate.

    Next, we will show that $\lupsilon_1$ cannot be the unique optimal price. Notice that for $\lupsilon_1$ to be the unique global maximum, the ordering $p^*_h(\lsigma,\beta) < \lupsilon_1 < p^*_m(\beta)$ has to hold. However, $p^*_h(\lsigma,\beta) \geq p^*_m(\beta)$, as both prices exhibit the form of $\beta - \sqrt{\beta(\beta-c)}$ with $c$ being $\mu$ for $p^*_m(\beta)$ and $\mu + {\lsigma^2}/{\mu}$ for $p^*_h(\lsigma,\beta)$. Hence, we can exclude $\lupsilon_1$ as an optimal price candidate as well. 

We next build on the insight that there are three prices, and the extreme cases covered in  Theorems~\ref{theorem_s} and \ref{theorem_h}, to obtain a more detailed characterization of the optimal pricing strategy based on three pairwise comparisons, for which we shall introduce three functions. 
The function $f$ is the value of $\usigma$ for which $R_1(p^*_l(\usigma)) = R_2(p^*_m(\beta))$ and corresponds with the threshold value $\bar{\sigma}$ that appeared in Proposition~\ref{theorem_s}: 
\begin{flalign}
    &&f(\mu,\beta) = \left(\tfrac{32}{27}(\beta-\mu)\left(\sqrt{\beta(\beta-\mu)}-(\beta-\mu)\right)\right)^{\frac{1}{2}}.&&
\end{flalign} 
The function $g$ is the value of $\lsigma$ for which $R_2(p^*_m(\beta)) = R_3(p^*_h(\lsigma,\beta))$. Since $R_2(p^*_m(\beta))$ does not depend on $\lsigma$, we write
\begin{flalign}
    &&g(\mu,\beta) = \left(2\beta\sqrt{R_2(p^*_m(\beta))\mu}-\beta R_2(p^*_m(\beta))-\mu^2\right)^{\frac{1}{2}}.&&
\end{flalign}
The function $h$ is the value of $\lsigma$ for which $R_1(p^*_l(\usigma)) = R_3(p^*_h(\lsigma,\beta))$. Since $R_1(p^*_l(\usigma))$ does not depend on $\lsigma$, we write %we can leave it implicit and write
\begin{flalign}    
    &&h(\usigma,\mu,\beta) = \left(2\beta\sqrt{R_1(p^*_l(\usigma))\mu}-\beta R_1(p^*_l(\usigma))-\mu^2\right)^{\frac{1}{2}}.&&
\end{flalign}
We can now characterize the optimal price in terms of risk.
%\begin{theorem}\label{cor_g}
The optimal robust price $p^* = \arg\sup_p\Tilde{R}(p)$ is then as in \eqref{p^*_gmain}
%\begin{flalign}\label{p^*_g2}
 %   &&p^* = 
  %  \begin{cases}
   % p^*_{l}(\usigma) \quad &  {\rm if} \ (\lsigma,\usigma)   \in \Sigma_l, \\
  %  p^*_{m}(\beta) \quad &  {\rm if} \ (\lsigma,\usigma) \in \Sigma_m,
  %  \\
   % p^*_{h}(\lsigma,\beta) \quad &  {\rm if} \ (\lsigma,\usigma) \in \Sigma_h,
    %\end{cases}&&
%\end{flalign}
with the risk threshold sets given by
\begin{flalign}\label{sets}
    && &\Sigma_{l} = \{(\lsigma,\usigma) \in \mathbf{R}_+^2|\lsigma \leq h(\usigma,\mu,\beta),\usigma \leq f(\mu,\beta), \lsigma \leq \usigma\},&& \nonumber \\
    && &\Sigma_{m} = \{(\lsigma,\usigma) \in \mathbf{R}_+^2|\lsigma \leq g(\mu,\beta),\usigma \geq f(\mu,\beta), \lsigma \leq \usigma\},&& \\
    && &\Sigma_{h} = \{(\lsigma,\usigma) \in \mathbf{R}_+^2|\lsigma \geq g(\mu,\beta),\lsigma \geq h(\usigma,\mu,\beta), \lsigma \leq \usigma\}.&&\nonumber
\end{flalign}

\subsection{Proof of Theorem~\ref{general_jump_in_words}}\label{ProofT2}
We now show that as the bounds of the standard deviation interval increase, a discontinuous jump eventually occurs in which we transition from the low price $p^*_l(\usigma)$ to a higher price, either $p^*_m(\beta)$ or $p^*_h(\lsigma,\beta)$. 

We begin by reformulating Theorem~\ref{general_jump_in_words}, originally stated in fairly loose terms, into a more quantifiable and precise form. We allow $\lsigma$ and $\usigma$ to increase in a general manner by letting $\usigma$ increase while imposing $\lsigma = \Phi(\usigma)$, where $\Phi$ is any non-decreasing function of $\usigma$ constrained by $\Phi(\usigma) \leq \usigma$. Special cases include the standard deviation as an upper bound when $\Phi(\usigma) = 0$ and as exact when $\Phi(\usigma) = \usigma$.

\begin{theorem}\label{general_jump}
When $\usigma < f(\mu,\beta)$ and $\Phi(\usigma) < h(\usigma,\mu,\beta)$, then $$\tilde{R}(p^*_l(\usigma)) < \max\{\tilde{R}(p^*_m(\beta)),\tilde{R}(p^*_h(\Phi(\usigma),\beta))\},$$ and when $\usigma > f(\mu,\beta)$ or $\Phi(\usigma) > h(\usigma,\mu,\beta)$, then $\tilde{R}(p^*_l(\usigma)) > \max\{\tilde{R}(p^*_m(\beta)),\tilde{R}(p^*_h(\Phi(\usigma),\beta))\}$.
\end{theorem}
\begin{customproof}
\noindent First, since $R_1(p)$ is maximal and continuously differentiable in $p^*_l(\usigma)$, we can employ the envelope theorem to show that
    \begin{flalign}
        &&\frac{\partial R_{1}(p^*_l(\usigma))}{\partial\usigma^2} = \frac{\partial R_{1}(p)}{\partial\usigma^2}\biggr\rvert_{p = p^*_l(\usigma)} = -\frac{p^*_l(\usigma)(\mu-p^*_l(\usigma))^2}{((\mu-p^*_l(\usigma))^2+\usigma^2)^2} \leq 0,&&
    \end{flalign}
    since $p^*_{l}(\usigma) \geq 0$. Second,
    \begin{flalign}
        &&\frac{\partial R_{2}(p^*_m(\beta))}{\partial\usigma^2} = 0,&&
    \end{flalign}
    as $R_2(p^*_m(\beta))$ does not depend on $\usigma$. Finally, since $R_3(p)$ is maximal and continuously differentiable in $p^*_h(\Phi(\usigma),\beta)$, the envelope theorem yields
    \begin{flalign}
        &&\frac{\partial R_{3}(p^*_h(\Phi(\usigma),\beta))}{\partial\Phi(\usigma)^2} = \frac{\partial R_{3}(p)}{\partial\Phi(\usigma)^2}\biggr\rvert_{p = p^*_h(\Phi(\usigma),\beta)} = \frac{p^*_h(\Phi(\usigma),\beta)}{\beta(\beta-p^*_h(\Phi(\usigma),\beta))} \geq 0,&&
    \end{flalign}
    since $p^*_h(\Phi(\usigma),\beta) \geq 0$. Hence, since $\Phi(\usigma)$ is increasing in $\usigma$, we also have that $R_{3}(p^*_h(\Phi(\usigma),\beta))$ is increasing in $\usigma$.
    
    Next, consider $(\Phi(\usigma),\usigma) = (0,0)$. Then $\uupsilon_1 = \lupsilon_1 = \lupsilon_2 = \mu$, thus $p^*=p^*_l(\usigma)$. Furthermore, consider $(\Phi(\usigma),\usigma) = (\Phi(\sigma_{\textup{max}}),\sigma_{\textup{max}}) = (\Phi(\sqrt{\mu(\beta-\mu)}),\sqrt{\mu(\beta-\mu)})$. Then $\lupsilon_1 = 0$, thus $p^* \in \{p^*_m(\beta),p^*_h(\Phi(\usigma),\beta)\}$.
\end{customproof}

Theorem \ref{general_jump} displays the robustness of our main insight: even when the standard deviation lies within an interval $[\lsigma,\usigma]$, an abrupt shift from a lower to a higher price still occurs. However, in contrast to the specific cases where the standard deviation is an upper bound or exact, the general scenario may exhibit two discontinuous jumps depending on the choice of $\Phi$. Notably, these shifts always follow a clear progression from the low price to the middle price, and then to the high price. This key insight follows from the proof of Theorem \ref{general_jump}, which reveals that $R_2(p^*_m(\beta))$ is independent of $\usigma$, while $R_3(p^*_h(\Phi(\usigma),\beta))$ is increasing in $\usigma$.

\section{Low and high pricing in service operations management}\label{sec:queue}
Monopoly pricing as considered so far considers a firm seeking to maximize the revenue function $p\bP(X\geq p)$ where $\bP(X\geq p)$ represents the demand function, i.e., the probability that a consumer is willing to pay a price $p$ for some product. When transitioning from products to services, the analysis becomes more intricate. In service systems, the demand function is not fixed but emerges from the decisions of rational customers who evaluate both the price and the expected waiting time. Customers will opt into the service only if the utility they derive exceeds the combined cost of the price and the waiting time. This creates a demand function that is endogenously determined through equilibrium analysis, shaped by the interplay of individual decisions and the operational characteristics of the system. This added complexity makes pricing for services a more challenging problem, as it requires integrating pricing strategies with queueing models to account for congestion and customer behavior.
We address robust pricing for services by applying the same max-min framework, focusing on information sets defined by the mean, variance, and range of market characteristics. As will become clear, the core research thread remains consistent: we observe shifts from lower to higher pricing strategies as the market becomes more volatile or heterogeneous. While the resulting price functions and thresholds are less explicit and more involved than in standard monopoly pricing, they underline the universality of the main insight: a robust seller should shift from high to low pricing as the market becomes sufficiently predictable---or conversely, switch to high pricing when the market becomes more uncertain.

\subsection{Heterogeneous valuation}\label{het:val}

Consider a monopolistic service firm that operates as an $M/M/1$ queue. Potential customers arrive according to a Poisson process with arrival rate $\lambda$. Customers who decide to join the queue are served in order of arrival by a single server. Service times are exponentially distributed with rate $\theta$.  
Following the approach in 
e.g.~\cite{stidham1992pricing,cachon2011pricing,haviv2014pricing}, customers have individual service valuations that are modeled as i.i.d.~random
variables. In traditional models, the valuation distribution is assumed to be known to the system manager, an assumption that we will relax below. Customers know their valuation but cannot observe the queue length, and decide whether to join
or not only based on their assessment of the expected waiting times. Such queues are called unobservable queues and widely studied in the field of rational queueing theory; see \cite{hassin2003queue} for an extensive overview.

For notational consistency, denote the random valuation of service by $X$, so that a customer only joins the queue when this $X$ exceeds the sum of price $p$ and holding costs $h W(\gamma)$ with $W(\gamma)$ some continuous function that increases with $\gamma$, where a standard choice would be $W(\gamma)=\frac{\gamma}{\theta(\theta-\gamma)}$, the expected waiting time in an $M/M/1$ queue  with effective arrival rate $\gamma$ and service rate $\theta$.  The parameter  $h$ is the cost of waiting per time unit and measures the customer's delay sensitivity. Let $\bP$ denote the distribution of $X$.
%and $ \bar{F}(x)=1-F(x)$. 
An arbitrary customer then joins with probability $\bP(X > p+h W(\gamma))$.~(For mathematical convenience we use a strict inequality $>$ instead of $\geq$. Note that by Lemma \ref{pdinf} this will not affect the worst-case tail bound for the ambiguity sets considered in this paper.) Since customers arrive according to a Poisson process with rate $\lambda$, the effective arrival rate $\gamma=\gamma(p,\bP)$ solves
\begin{flalign}\label{eq:gamma}
&&\gamma=\lambda \bP(X > p+h W(\gamma)).&&
\end{flalign}
It is a priori not clear if \eqref{eq:gamma} has  a solution for given $p$ and $\Prob$; we will return to this matter later on. Note that $\lambda$ and $h$ are considered to be given constants at this point.

The above idea that individual utility leads to a global equilibrium is one of the key concepts in rational queueing theory and explained in detail in the book
\cite{hassin2003queue}. See e.g.~\cite{naor1969regulation,dewan1990user,chen2004monopoly,taylor2018demand,baron2023omnichannel,chen2022food} for examples of pricing in this area, and in particular   \cite{cachon2011pricing,haviv2014pricing} consider monopoly pricing for the same model as in this section. However, all the above references assume that the valuation distribution is known, and satisfies certain regularity conditions, whereas we shall instead pursue a robust approach for when only limited information on the valuation distribution is available. Regularity conditions in \cite{haviv2014pricing} include 
that the valuation distribution has a continuously differentiable density and a nondecreasing hazard rate, which warrants existence and uniqueness of the solution to \eqref{eq:gamma}.

A revenue-maximizing seller seeks to solve $\sup_p p \cdot \gamma(p,\bP)$ and we consider the robust setting where the seller does not know the distribution $\bP$ of $X$, but instead only knows that the valuation distribution lies in the ambiguity set $\mathcal{P}$. Taking a pessimistic approach, the seller then sets a deterministic price that maximizes the expected revenue under the worst-case scenario.

It is not guaranteed that \eqref{eq:gamma} has  a solution for given $p$ and $\Prob \in \cP$. We therefore set $\gamma(p,\bP) = \infty$ if no solution to \eqref{eq:gamma} exists. We assume that for every $p$, there is at least one $\Prob \in \cP$ for which $\gamma(p,\Prob) < \infty$. This is in particular true for ambiguity sets satisfying the two conditions in Theorem \ref{th:queue}, as becomes clear from the proof. We define \vspace*{-0.3cm}
\begin{flalign}
    &&\gamma^*(p,\cP) =  \inf_{\Prob \in \cP} \gamma(p,\Prob),&&
\label{eq:gamma_star}
\end{flalign}
and then the objective of a robust seller is to solve $\sup_p\, p \cdot \gamma^*(p,\cP)$.

\begin{theorem}[Robust Endogenous Demand I]\label{th:queue}
Let $p >0 $ and ambiguity set $\mathcal{P}$ be given, and assume that:
\begin{itemize}
    \item[{\rm (i)}] $f(x) \eqdef \inf_{\Prob \in \cP} \Prob(X > p + hW(x))$ is continuous.
    \item[{\rm (ii)}] $\inf_{\Prob \in \cP} \Prob(X > p + hW(x))$ is attained by some distribution $\Prob_x \in \cP$ for every $x$.
\end{itemize}
Then $\gamma^* = \gamma^*(p)$ defined in \eqref{eq:gamma_star} is the unique solution to the equation
\begin{flalign}\label{thmrq}
    &&\gamma^* = \lambda \cdot \inf_{\Prob \in \cP} \Prob(X > p + hW(\gamma^*)).&&
\end{flalign}
\end{theorem}

\iffalse

\fi

The proof can be found in Appendix \ref{ap:queue1}. Lemma~\ref{pdinf} shows that    
Theorem~\ref{th:queue} applies to all ambiguity sets considered in this paper: The function \eqref{ttb} is always continuous in $p$, and worst-case distributions are identified in the proof of Lemma \ref{pdinf} in Appendix \ref{pdp}. 
Take as example the ambiguity set  $\cP(\mu,\beta)$ of all distributions with given mean  $\mu$ and maximal valuation $\beta$.

\begin{figure}[ht]
\centering
\subfigure[Worst-case and upper bound on arrival rate.]{
\label{fig:A}
\begin{tikzpicture}
    \node at (0,0) {\includegraphics[width=0.45\linewidth]{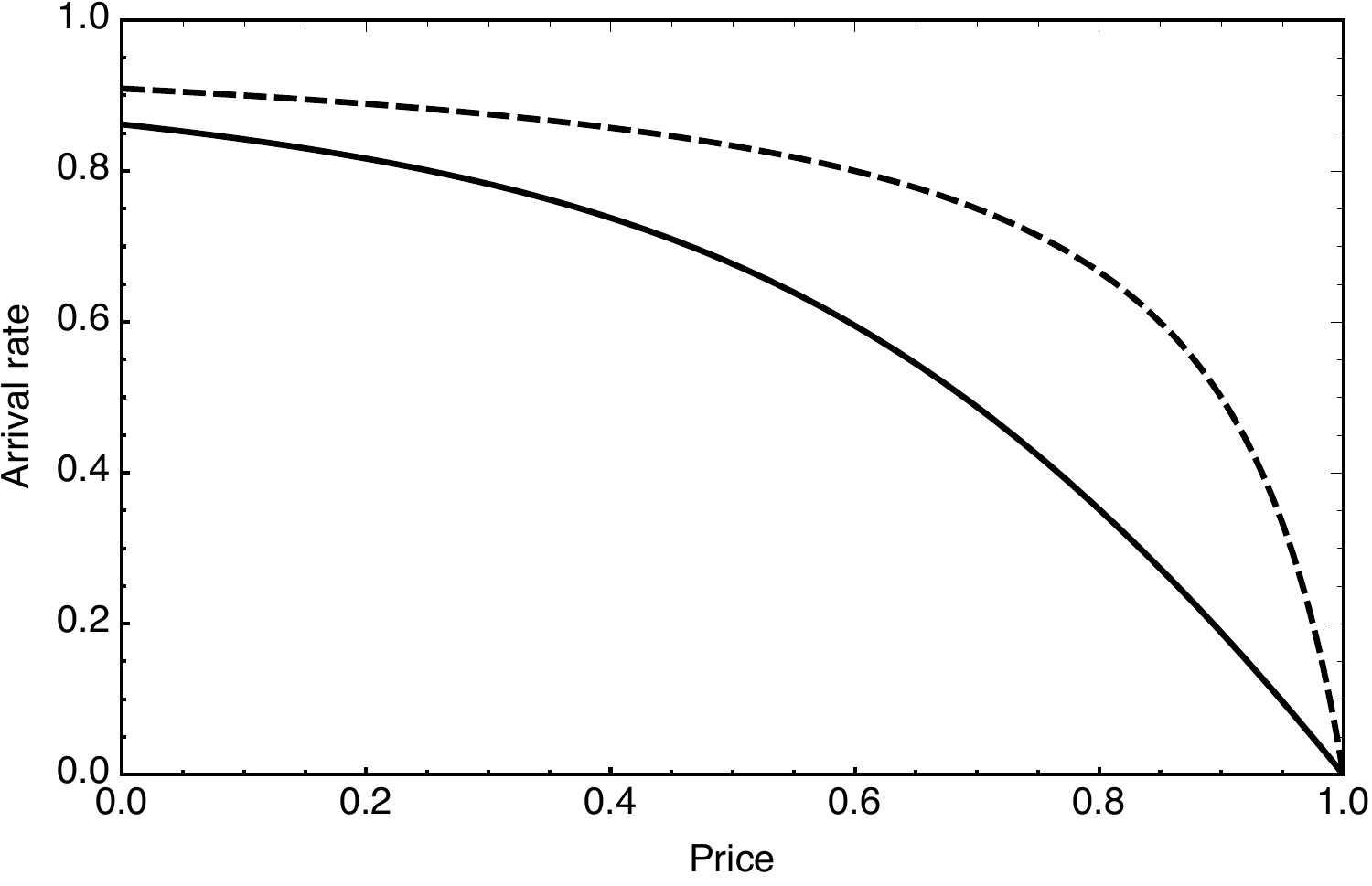}};
      \node at (.5,.5) {\begin{scriptsize}$\gamma^*$\end{scriptsize}};
       \node at (1.6,1.5) {\begin{scriptsize}$\gamma_{\max}$\end{scriptsize}};
    \end{tikzpicture}
}
\subfigure[Prices maximizing expected revenue.]{
\label{figmain1hoptQ}
\begin{tikzpicture}
    \node at (0,0) {\includegraphics[width=0.45\linewidth]{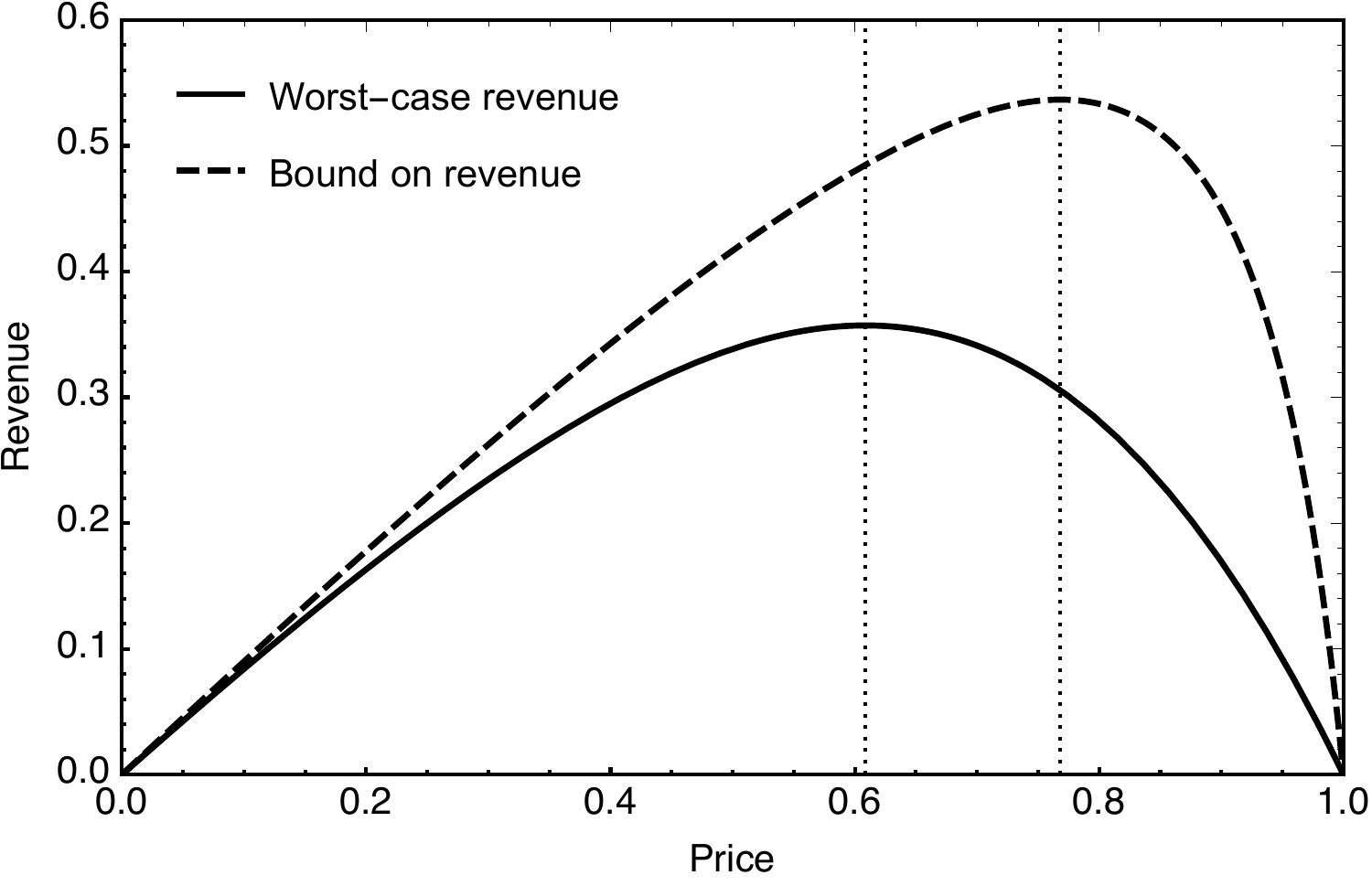}};
       \node at (1.2,-1.45) {\begin{scriptsize}$p^*$\end{scriptsize}};
       \node at (2.4,-1.5) {\begin{scriptsize}$p_{\max}$\end{scriptsize}};
    \end{tikzpicture}
}
\caption{Pricing in delay-prone services when only the mean and maximal valuation are known, with $\theta=1$, $\lambda=10$, $\mu=1$,  $\beta = 5$ and  $h = 0.1$.}\label{plotAB}
\end{figure}

Lemma~\ref{pdinf} with 
$\lsigma=0$ and $\usigma=\sqrt{\mu(\beta-\mu)}$ then yields
$\inf_{\bP\in\cP(\mu,\beta)}\bP(X> x) = \frac{\mu-x}{\beta-x}$ for $x \in (0,\mu)$, and hence, by taking $x = p + hW(\gamma^*)$, we have that $\gamma^*$
solves
   $\frac{\mu-p-hW(\gamma^*)}{\beta-p-hW(\gamma^*)}=\frac{\gamma^*}{\lambda}$.
 For  $W(\gamma)=\frac{\gamma}{\theta(\theta-\gamma)}$  the worst-case arrival rate $\gamma^*$ can be obtained in closed-form. Moreover, for this case, we can also derive an explicit upper bound on $\gamma^*$ using the fact that 
$\mu-p-hW(\gamma^*)$ needs to be nonnegative. This implies that
$\gamma^*\leq \gamma_{\max}\eqdef \theta -\frac{\theta  h}{\theta  \mu +h-\theta  p}$ and this
bound becomes exact in the large market size limit $\lambda\to\infty$. 
This also gives an upper bound  $p \gamma_{\max}$ on the worst-case revenue, which is maximized by
$p_{\max}\eqdef\mu-\frac{\sqrt{h^2+\theta  h \mu }-h}{\theta}$. Pricing in this service system depends on the delay function $W(\cdot)$, the service rate $\theta$, and the additional parameter $h$. Observe that for $h=0$ we retrieve the classical monopoly pricing problem. See Figure~\ref{plotAB} for an example.
Knowing the mean and maximal value thus leads to a unique optimal price that decreases with $\theta$. As in the original monopoly pricing problem with information on the mean and maximal value, there is no low-high pricing phenomenon. The result, however, is interesting in its own right, as distribution-free pricing for these service systems has not been covered in the literature. 

\begin{figure}[ht]
\centering
\subfigure[Function $p\mapsto\gamma^*(p,\lambda)$.]{
\label{figmain1hQ}
\begin{tikzpicture}
    \node at (0,0) {\includegraphics[width=0.46\linewidth]{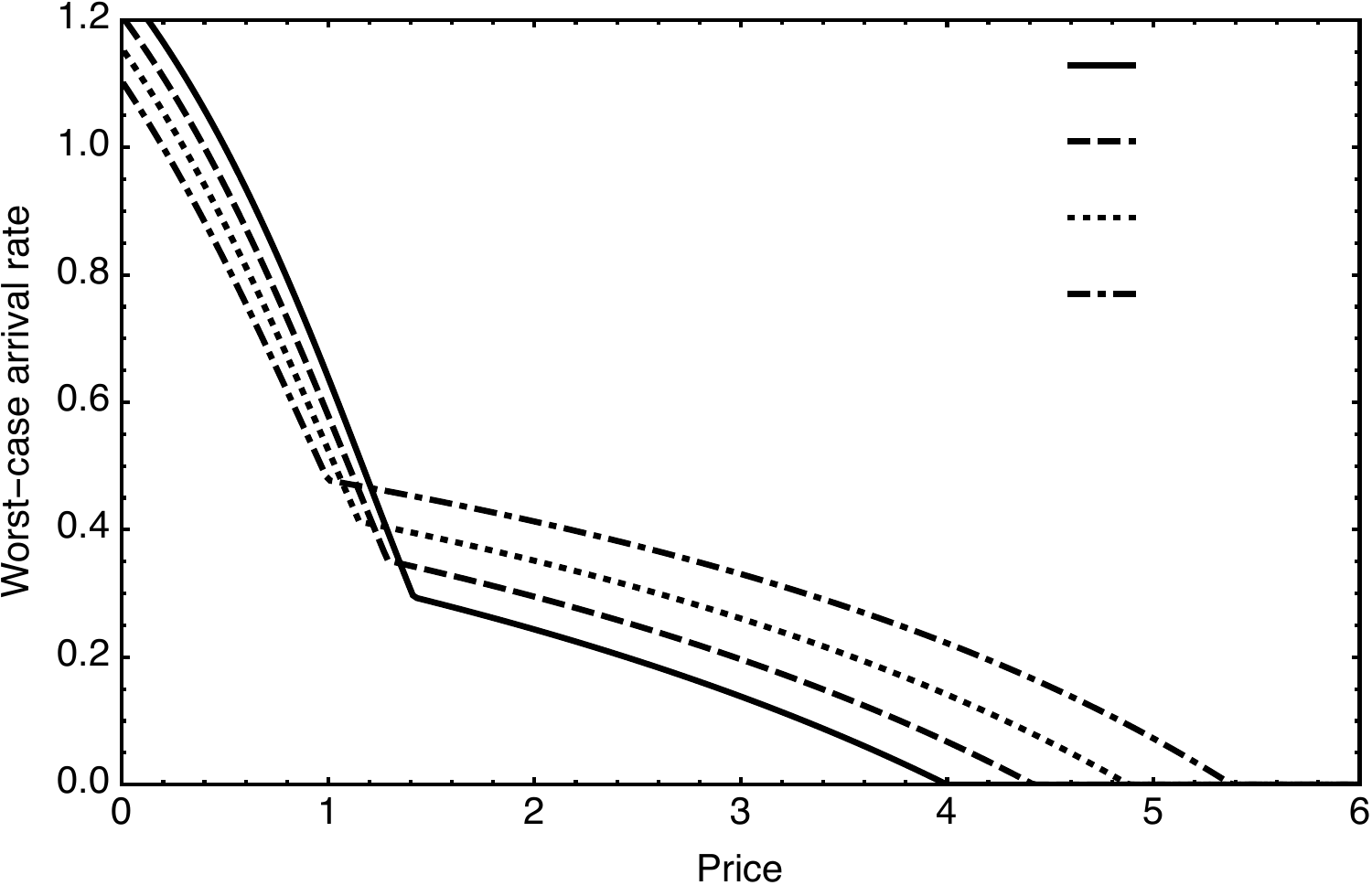}};
    \node at (3.1, 2.13)   {\begin{scriptsize}$ \sigma=2.0$\end{scriptsize}};
     \node at (3.1, 1.73)   {\begin{scriptsize}$ \sigma=2.2$\end{scriptsize}};
     \node at (3.1, 1.32)   {\begin{scriptsize}$ \sigma=2.4$\end{scriptsize}};
     \node at (3.1, 0.91)   {\begin{scriptsize}$ \sigma=2.6$\end{scriptsize}};
    \end{tikzpicture}
  }
\subfigure[Function $p\mapsto p\gamma^*(p,\lambda)$.]{
\label{figmain1hoptQ}
\begin{tikzpicture}
    \node at (0,0) {\includegraphics[width=0.46\linewidth]{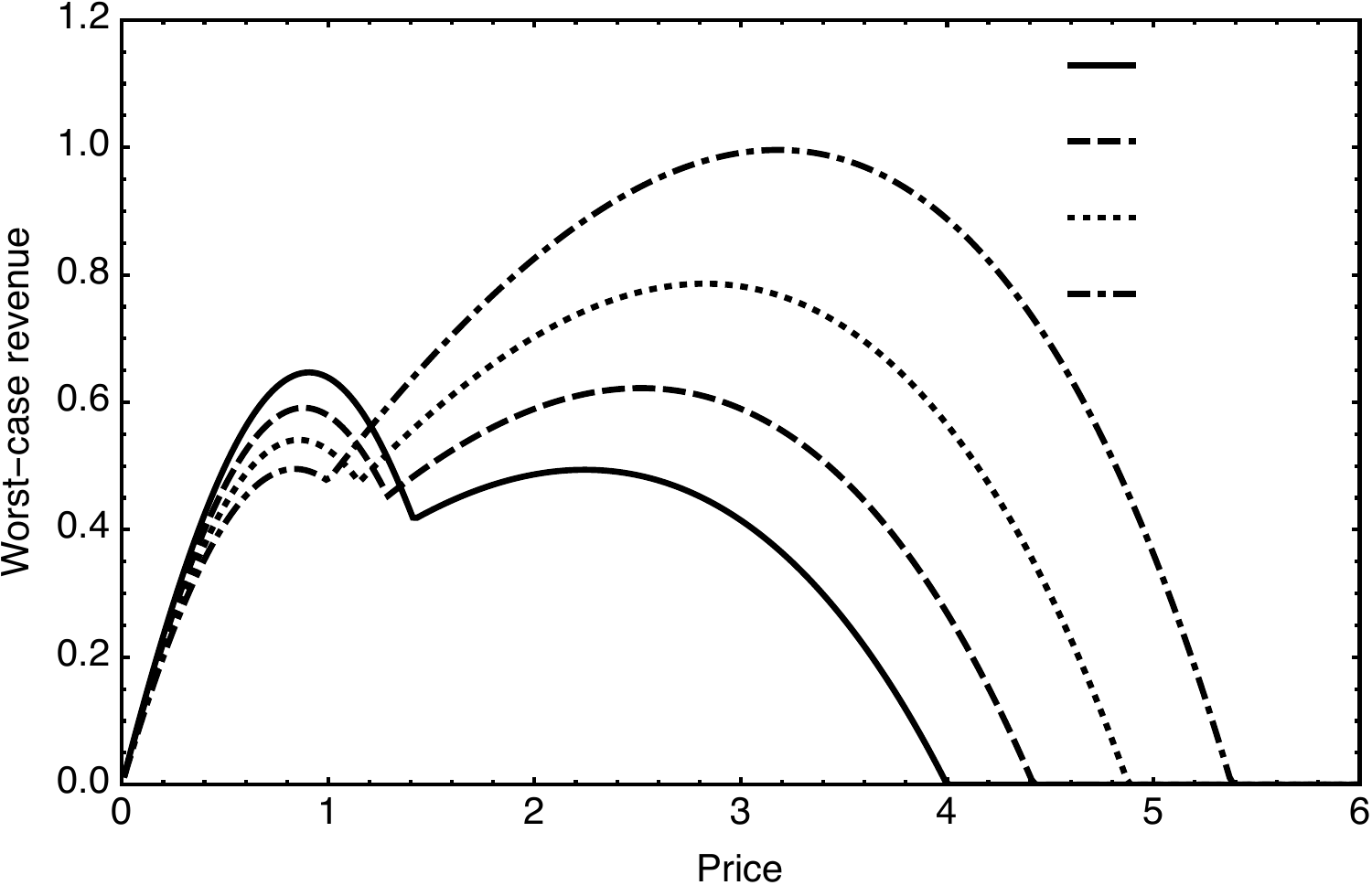}};
    \node at (3.1, 2.13)   {\begin{scriptsize}$ \sigma=2.0$\end{scriptsize}};
     \node at (3.1, 1.73)   {\begin{scriptsize}$ \sigma=2.2$\end{scriptsize}};
     \node at (3.1, 1.32)   {\begin{scriptsize}$ \sigma=2.4$\end{scriptsize}};
     \node at (3.1, 0.91)   {\begin{scriptsize}$ \sigma=2.6$\end{scriptsize}};
    \end{tikzpicture}
}
\caption{Illustration of queue with $\lambda=5$, $\mu= 2$, $\theta = 2$, $\beta = 10$ and $h = 1$.}\label{FigmainhQ3}
\end{figure}

We now return to the main theme of this paper and examine the pricing strategy a seller should adopt when there is information about the variation of the valuation distribution, and thus, the demand uncertainty. Specifically, we assume the seller has knowledge of the mean $\mu$, the support bound $\beta$, and the standard deviation $\sigma$. While the framework extends to cases where the seller only knows that $\sigma$ lies within a certain interval, for clarity of presentation, we focus on the scenario where $\sigma$ is precisely known.
Theorem~\ref{th:queue} implies that the worst-case tail bound, as established in Lemma~\ref{pdinf}, can be utilized to solve the max-min pricing problem given by $\sup_p p\gamma^*$, assuming $\lsigma = \usigma$. Figure~\ref{FigmainhQ3} illustrates how the optimal arrival rate $\gamma^*$ varies with $\sigma$, leading to a bimodal revenue structure reminiscent of the original monopoly pricing problem. Notably, as standard deviation increases beyond a certain threshold, a transition occurs from low to high pricing, highlighting a fundamental low-to-high pricing phenomenon. Furthermore, we show another switching effect in terms of the support upper bound. As shown in Figures~\ref{figmain1hQ} and \ref{figmain1hoptQ}, a decrease in $\beta$ also triggers a switch from low to high pricing. We capture these observations in the following theorem:
\begin{theorem}[Price Jumps in Queues I]\label{THMswithQ}
Assume that $W$ depends on $\sigma$ and $\beta$ solely through $\gamma^*$. Consider the optimal price solving $\sup_p p\gamma^*$ with the distribution of the valuation $X$ contained in $\cP(\mu,\sigma,\beta)$. There exist unique thresholds $\sigma^*$ and $\beta^*$ such that the optimal pricing switches from low pricing for $\sigma\leq \sigma^*$ and $\beta\geq \beta^*$ to high pricing for $\sigma\geq \sigma^*$ and $\beta\leq \beta^*$.
\end{theorem}
Theorem~\ref{THMswithQ} is proved in Appendix~\ref{proofQQs}.
Let us elaborate on its implications. Why does restricting valuations through $\beta$ lead the seller to post higher prices, and vice versa? Allowing larger valuations broadens the set of possible valuation distributions, enabling the adversary to concentrate mass on low valuations while still satisfying the mean constraint. This forces the seller to choose a low price. Conversely, restricting the adversary's ability to assign high valuations compels the seller to adopt a high-price strategy. These insights are universal and extend beyond the specific settings of monopoly pricing for products and delay-prone services discussed in this paper.

\begin{figure}[ht!]
\centering
\subfigure[Function $p\mapsto p\gamma^*(p,\lambda)$ with $\beta=10$.]{
\label{figmain1hQ}
\begin{tikzpicture}
    \node at (0,0) {\includegraphics[width=0.46\linewidth]{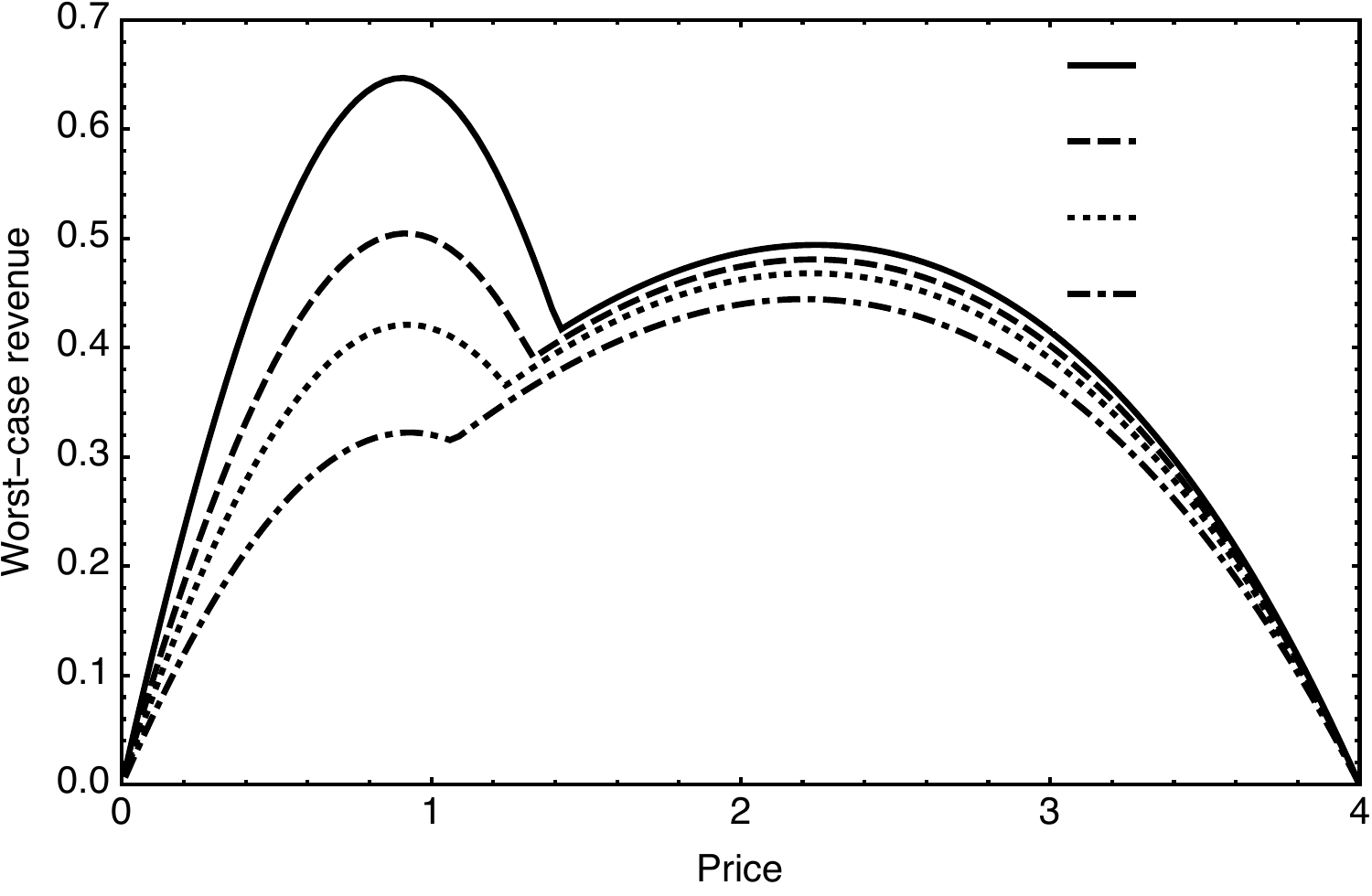}};
    \node at (3.1, 2.13)   {\begin{scriptsize}$ h=1$\end{scriptsize}};
     \node at (3.1, 1.73)   {\begin{scriptsize}$ h=2$\end{scriptsize}};
     \node at (3.1, 1.32)   {\begin{scriptsize}$ h=3$\end{scriptsize}};
     \node at (3.1, 0.89)   {\begin{scriptsize}$ h=5$\end{scriptsize}};
    \end{tikzpicture}
}
\subfigure[Function $p\mapsto p\gamma^*(p,\lambda)$ with $h=1$.]{
\label{figmain1hoptQ}
\begin{tikzpicture}
    \node at (0,0) {\includegraphics[width=0.46\linewidth]{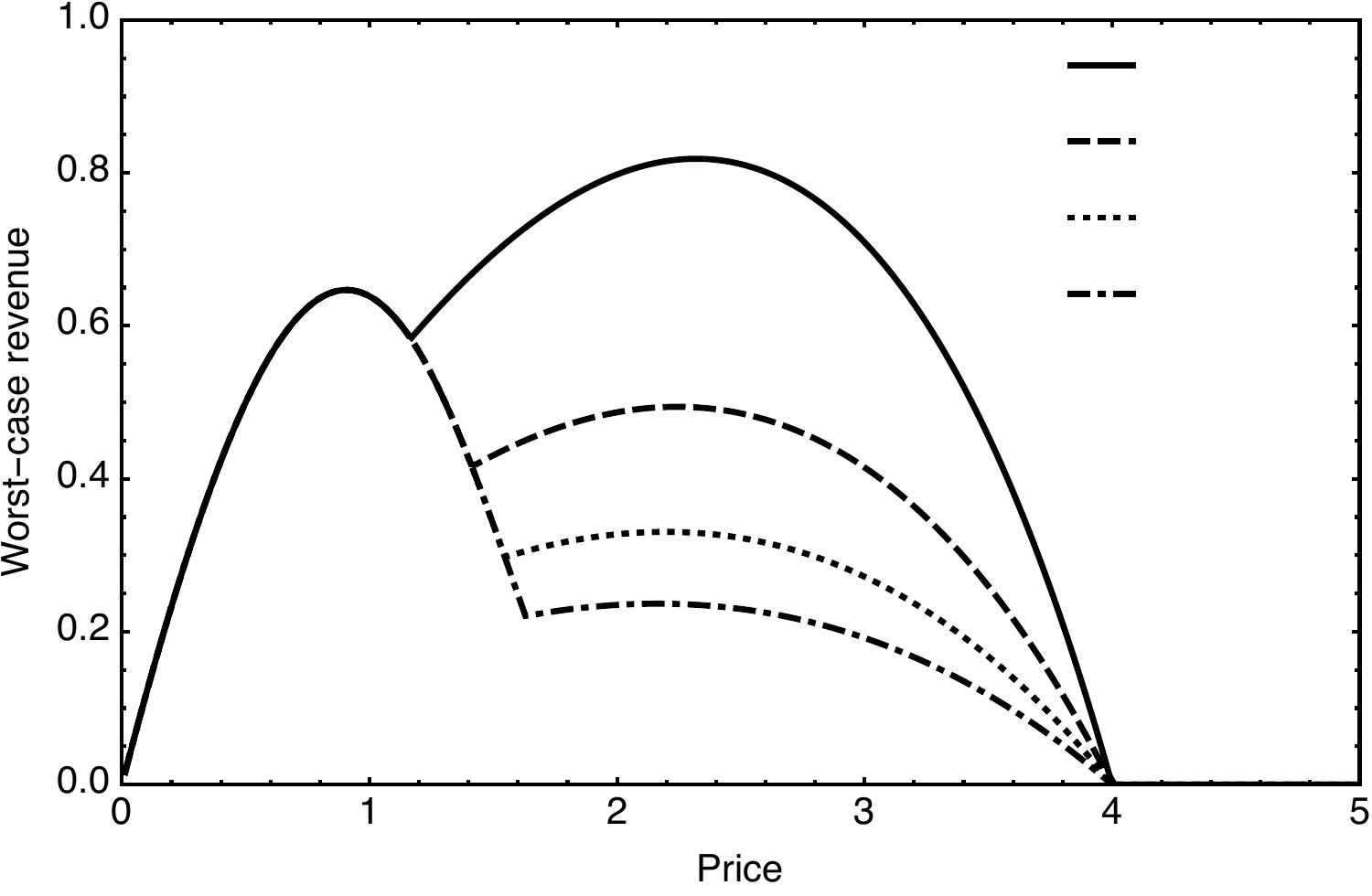}};
    \node at (3.02, 2.13)   {\begin{scriptsize}$ \beta=8$\end{scriptsize}};
     \node at (3.1, 1.73)   {\begin{scriptsize}$ \beta=10$\end{scriptsize}};
     \node at (3.1, 1.31)   {\begin{scriptsize}$ \beta=12$\end{scriptsize}};
     \node at (3.1, 0.89)   {\begin{scriptsize}$ \beta=14$\end{scriptsize}};
    \end{tikzpicture}
}
\caption{Illustration of queue with $\lambda=5$, $\mu= 2$, $\sigma=2$, $\theta = 2$ and different values of $h$ and $\beta$.}\label{FigmainhQ}
\end{figure}

\subsection{Heterogeneous delay cost}\label{het:cos}
In rational queueing theory, customer behavior is often modeled based on differences in how much they value the service itself, influencing their willingness to pay for priority or join the queue. This approach was adopted in Section~\ref{het:val} and works well in scenarios where customers derive varying benefits from the service. For example, diners at a luxury restaurant may pay more for faster seating to enhance their experience, while concertgoers might purchase VIP tickets for shorter wait times or expedited entry. An alternative perspective focuses on differences in delay costs rather than service valuation. Here, the service is assumed to have equal value for all customers, but individuals vary in their sensitivity to waiting. Those with urgent tasks or high opportunity costs experience greater discomfort from delays, often opting for faster, more expensive options like express lanes or priority services. Meanwhile, customers with lower delay costs are more willing to wait. Examples include airline check-ins, where passengers with tight schedules prioritize speed, or healthcare, where critical patients receive faster treatment due to higher delay costs.

The delay cost rate thus reflects heterogeneity among customers, and accurately characterizing the distribution of delay-cost rates is essential for devising effective pricing strategies.

Much of the existing literature assumes specific forms for the distribution of delay costs. Several works focus on discrete distributions that model multiple types of customers \citep{Ha2001, Plambeck2004, Afeche2013, Ata2018}. Others investigate continuous distributions, such as uniform delay-cost rates \citep{Gilland2009, Yu2013}. Empirical distributions are derived in \citet{Lu2013} from supermarket sales data and \citet{Aksin2013} use a log-normal distribution based on call center data.
Some studies explore broader classes of delay-cost distributions under general regularity conditions \citep{cao2019priority, Nazerzadeh2018}.  In contrast to all above works, we adopt a distributionally robust approach, avoiding assumptions about a specific distribution. Instead, we consider all potential distributions within an ambiguity set and explore the delay-cost approach for the robust max-min setting. 

Assume that the delay cost rate $H$, measured in dollars per unit time, is a positive random variable. Arriving customers cannot observe queue lengths, but customers do receive an expected delay announcement $W(\gamma)$. Given this information, the customer decides to join by comparing the
utility of receiving service $R$ with the cost of waiting. A customer whose delay cost rate is $h$ experiences a disutility of joining of $p+h W(\gamma)$ and joins if $p+h W(\gamma)<R$. The fraction of customers joining will thus be
$\bP\left(H\leq\frac{R-p}{W(\gamma)}\right)$. Since customers arrive according to a Poisson process with rate $\lambda$, the effective arrival rate $\gamma_H=\gamma_H(p,\bP)$ solves
$\gamma_H=\lambda \Prob\Big(H\leq\frac{R-p}{W(\gamma_H)}\Big)$.
Equivalent to \eqref{eq:gamma_star}, we define \vspace*{-0.3cm}
\begin{flalign}
    &&\gamma^*_H(p,\cP) :=  \inf_{\Prob \in \cP} \gamma_H(p,\Prob),&&
\label{eq:gamma_star_H}
\end{flalign}
and the seller's objective $\sup_p\, p \cdot \gamma^*_H(p,\cP)$. The following theorem then gives a robust bound on the demand:
\begin{theorem}[Robust Endogenous Demand II]\label{th:queue2}
Let $p >0 $ and ambiguity set $\mathcal{P}$ be given, and assume that:
\begin{itemize}
    \item[{\rm (i)}] $g(x) \eqdef \inf_{\Prob \in \cP} \Prob(H<\frac{R-p}{W(x)})$ is continuous.
    \item[{\rm (ii)}] $\inf_{\Prob \in \cP} \Prob(H<\frac{R-p}{W(x)})$ is attained by some distribution $\Prob_x \in \cP$ for every $x$.
\end{itemize}
Then $\gamma^*_H = \gamma^*_H(p)$ defined in \eqref{eq:gamma_star_H} is the unique solution to the equation
\begin{flalign}\label{thmrq2}
    &&\gamma^*_H = \lambda \cdot \inf_{\Prob \in \cP} \Prob\Big(H<\frac{R-p}{W(\gamma^*_H)}\Big).&&
\end{flalign}
\end{theorem}
The proof is provided in Appendix \ref{ap:queue2}. A notable difference between \eqref{thmrq} and \eqref{thmrq2} is the requirement of the supremum rather than the infimum of the tail bound. This follows from the relation
$$\inf_{\Prob \in \cP} \Prob\Big(H<\frac{R-p}{W(\gamma^*_H)}\Big) = 1 - \sup_{\Prob \in \cP} \Prob\Big(H\geq\frac{R-p}{W(\gamma^*_H)}\Big).$$ We derive this tail bound in Appendix \ref{ap:tutb} for $\cP = \cP(\mu,\lsigma,\usigma,\beta)$. Both assumptions in Theorem \ref{th:queue2} hold for all ambiguity sets considered in this paper. We also observe a price jump phenomenon.
\begin{theorem}[Price Jumps in Queues II]\label{THMswithQ2}
Assume that $W$ depends on $\sigma$ and $\beta$ solely through $\gamma^*_H$. Consider the optimal price solving $\sup_p p\gamma^*_H$ with the distribution of the valuation $H$ contained in $\cP(\mu,\sigma,\beta)$. There exist unique thresholds $\sigma^*$ and $\beta^*$ such that the optimal pricing switches from low pricing for $\sigma\leq \sigma^*$ and $\beta\geq \beta^*$ to high pricing for $\sigma\geq \sigma^*$ and $\beta\leq \beta^*$.
\end{theorem}
The proof is provided in Appendix \ref{ap:proofTHMQ2}. Observe that the assumption in Theorem \ref{THMswithQ2} is satisfied for the expected waiting time in an $M/M/1$ queue.

\section{Conclusions and outlook}\label{sec:con}
The transition from low to high pricing in response to increased market volatility has not been previously documented in the robust pricing literature. When the seller knows only the mean and variance, prices tend to be lower and decrease with risk. Likewise, relying solely on valuation bounds leads to risk-insensitive pricing. The distinctive shift from low to high pricing occurs only when the seller, in addition to knowing the mean and valuation bounds, has sufficiently precise information about the variance. This specific combination of factors may explain why this phenomenon has been overlooked. Another contributing factor could be the dominance of low-pricing strategies in prior research, which has largely focused on the intuition that prices should decrease with risk or volatility.  

The primary goal of this paper was to establish the low-high pricing effect as a consequence of increased segmentation in the worst-case market. To formally prove these results, we adopted a distributionally robust approach, where the seller solves a max-min optimization problem to identify the worst-case market among all distributions with a known mean, maximum valuation, and dispersion (the ambiguity set). Ambiguity sets beyond the scope of this paper include those incorporating the first three or more moments or those defined using distance metrics, such as \(\phi\)-divergence, Kullback-Leibler distance, or Wasserstein distance. Empirically driven ambiguity sets, which integrate an increasing number of percentiles (based on past pricing responses), also offer an interesting avenue for future research. 
While requiring a different mathematical framework, such sets might in some cases also lead to low-high pricing, particularly under high levels of ambiguity, though this is far from guaranteed. Indeed, not all ambiguity sets generate segmentation: for example, ambiguity sets defined only by mean and variance lead to unimodal worst-case distributions and hence more conservative, gradually adjusting prices. By contrast, the mean–variance–range ambiguity set considered in this paper does admit worst-case distributions with genuine segmentation, which in turn creates the observed low–high pricing effect.

The key insight is that the low–high pricing phenomenon arises only when the available information admits genuine market segmentation. In such cases, the seller does not adjust prices gradually with increasing risk, but instead switches discontinuously between low and high prices. Segmentation of this kind may also emerge directly from empirical data when distinct consumer groups are present, or from irregular but stylized price–response functions. What is distinctive in our setting is that segmentation is not imposed a priori through assumptions or fitted distributions, but rather emerges endogenously as the important worst-case market under the robust mean–variance–range model. Put differently, when an ambiguity set allows for segmentation, the worst-case distribution it produces is not an artifact of modeling but the critical scenario in the max–min sense. Robustness then guides the seller to prepare for exactly this possibility: a segmented market with discontinuous low–high pricing.

Building on these insights, we also detected the low–high pricing phenomenon in queueing models. Theorems~\ref{th:queue} and \ref{th:queue2} are more general than strictly necessary for the objectives of this paper. While our focus is on the M/M/1 queue with mean-variance range ambiguity, these theorems extend to more general queueing models. For instance, in an M/G/1 queue, the delay function can be expressed as \(W(\gamma) = \frac{1+c_s^2}{2}\frac{\gamma}{\theta(\theta-\gamma)}\), where \(c_s\) is the coefficient of variation of the service times. Since this delay function is continuous, non-decreasing in the effective arrival rate, and independent of the ambiguity set parameters, all insights and regime-switching phenomena identified for M/M/1 carry over directly to the M/G/1 setting. 
More broadly, exploring other queueing models, such as G/G/1 or M/G/s systems, would present additional technical challenges, but pursuing this direction is worthwhile. Similarly, studying alternative ambiguity sets beyond mean-variance range, such as mean-support, higher moments, or distance-based sets, offers another interesting avenue for future research. As such, this paper introduces the first robust analysis of rational queues and highlights several promising directions for extending these results, including heterogeneous valuations and delay sensitivities.

\vspace{.2cm}
\noindent {\bf Acknowledgments}
The authors thank Alex Suzdaltsev for exchanging ideas about some unpublished results in \cite{suzdaltsev2018distributionally}.

\bibliographystyle{informs2014}
\bibliography{bibbook}
\newpage
\begin{appendices}

\section{Tail bound results}
\subsection{Proof of Lemma~\ref{pdinf}}\label{pdp}
Let $\mathcal{M}$ be the set of probability measures defined on the measurable space $([0,\beta],\mathcal{B})$ with $\mathcal{B}$ the Borel sigma-algebra on $\mathbb{R_+}$. We solve $\inf_{\bP \in \cP(\mu,\lsigma,\usigma,\beta)}\bP(X > p)$, which can be rewritten as
\begin{equation}\label{primal1}
\begin{aligned}
&\! \inf_{\mathbb{P}\in \mathcal{M}} &  &\int_x \1{\{x> p\}}{\rm d} \mathbb{P}(x),\\
&\text{s.t.} &      &  \int_x {\rm d}\mathbb{P}(x)=1, \ \int_x x{\rm d}\mathbb{P}(x)=\mu,\ \int_x x^2{\rm d}\mathbb{P}(x)\geq\lsigma^2+\mu^2,\ \int_x x^2{\rm d}\mathbb{P}(x)\leq\usigma^2+\mu^2.
%&               &      & \int p(x){\rm d}x=1\\
%&               &      & p(x)\geq 0.\\
\end{aligned}
\end{equation}
Consider the corresponding dual, 
\begin{equation}\label{dual1}
\begin{aligned}
&\sup_{\lambda_0,\lambda_1,\lambda_2,\lambda_3} &  &\lambda_0 + \lambda_1 \mu+\lambda_2 (\lsigma^2+\mu^2) + \lambda_3 (\usigma^2+\mu^2),\\
&\text{s.t.} &      & \1{\{x> p\}}\geq \lambda_0  +\lambda_1 x+(\lambda_2+\lambda_3) x^2 =: F(x), \ \forall x\in[0,\beta], \lambda_2 \geq 0, \lambda_3 \leq 0.
\end{aligned}
\end{equation}

The proof is structured in the following manner: Initially, four scenarios are considered, each corresponding to a feasible solution for the dual problem. Subsequently, we construct four primal solutions corresponding to these scenarios. The final step is the demonstration that strong duality is maintained across all possible values of $p$, with each scenario uniquely aligned with a specific interval of $p$. We now consider these scenarios.

\begin{figure}[ht!]
    \centering
    \begin{tikzpicture}
    \node at (0,0) {\includegraphics[width=\linewidth]{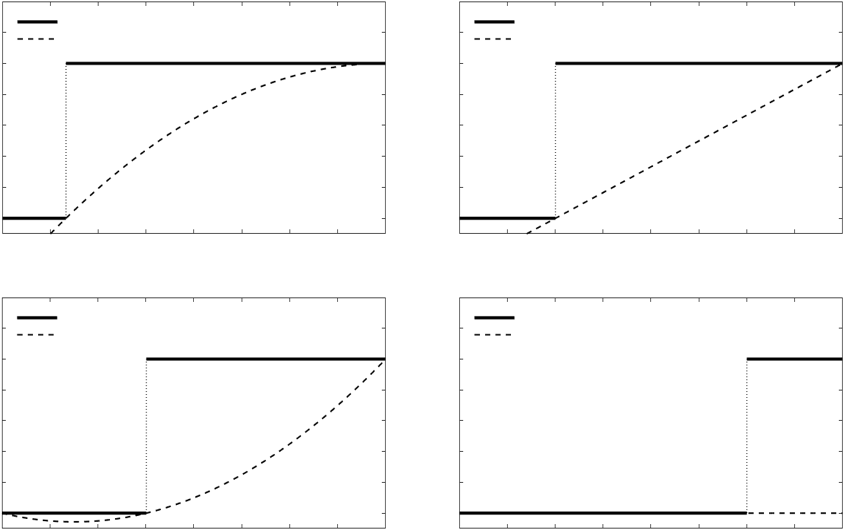}};
    \node at (-6, 4.52)   {\begin{scriptsize}$ \1\{X> p\}$\end{scriptsize}};
    \node at (2.5, 4.52)   {\begin{scriptsize}$ \1\{X> p\}$\end{scriptsize}};
    \node at (-6, -0.96)   {\begin{scriptsize}$ \1\{X> p\}$\end{scriptsize}};
    \node at (2.5, -0.96)   {\begin{scriptsize}$ \1\{X> p\}$\end{scriptsize}};
    \node at (-6.3, 4.18)   {\begin{scriptsize}$ F(x)$\end{scriptsize}};
    \node at (2.17, 4.19)   {\begin{scriptsize}$ F(x)$\end{scriptsize}};
    \node at (-6.32, -1.31)   {\begin{scriptsize}$ F(x)$\end{scriptsize}};
    \node at (2.17, -1.31)   {\begin{scriptsize}$ F(x)$\end{scriptsize}};
    \node at (-7.76, 0.35)   {\begin{scriptsize}$ 0$\end{scriptsize}};
    \node at (-6.6, 0.35)   {\begin{scriptsize}$ p$\end{scriptsize}};
    \node at (-0.76, 0.35)   {\begin{scriptsize}$ \beta$\end{scriptsize}};
    \node at (0.7, 0.35)   {\begin{scriptsize}$ 0$\end{scriptsize}};
    \node at (2.45, 0.35)   {\begin{scriptsize}$ p$\end{scriptsize}};
    \node at (7.75, 0.35)   {\begin{scriptsize}$ \beta$\end{scriptsize}};
    \node at (-7.75, -5.12)   {\begin{scriptsize}$ 0$\end{scriptsize}};
    \node at (-5.1, -5.12)   {\begin{scriptsize}$ p$\end{scriptsize}};
    \node at (-0.74, -5.12)   {\begin{scriptsize}$ \beta$\end{scriptsize}};
    \node at (0.75, -5.12)   {\begin{scriptsize}$ 0$\end{scriptsize}};
    \node at (6.05, -5.12)   {\begin{scriptsize}$ p$\end{scriptsize}};
    \node at (7.75, -5.12)   {\begin{scriptsize}$ \beta$\end{scriptsize}};
    \node at (-8, 0.86)   {\begin{scriptsize}$ 0$\end{scriptsize}};
    \node at (-8, 3.75)   {\begin{scriptsize}$ 1$\end{scriptsize}};
    \node at (0.48, 0.87)   {\begin{scriptsize}$ 0$\end{scriptsize}};
    \node at (0.48, 3.75)   {\begin{scriptsize}$ 1$\end{scriptsize}};
    \node at (-8, -4.60)   {\begin{scriptsize}$ 0$\end{scriptsize}};
    \node at (-8, -1.75)   {\begin{scriptsize}$ 1$\end{scriptsize}};
    \node at (0.48, -4.60)   {\begin{scriptsize}$ 0$\end{scriptsize}};
    \node at (0.48, -1.72)   {\begin{scriptsize}$ 1$\end{scriptsize}};
    \end{tikzpicture}
    \caption{Visualization of scenario 1 (top-left), scenario 2 (top-right), scenario 3 (bottom-left), and scenario 4 (bottom-right) with $\mu=1$, $\underset{\bar{}}{\sigma} = 0.1$, $\usigma = 1$, and $\beta = 4$.}
    \label{fig:inf_pd}
\end{figure}

In the first scenario, $F(x)$ intersects $\1\{x> p\}$ in $\{p, \alpha\}$ with $\alpha\in[p,\beta]$. In the second scenario, the intersection occurs in the points $\{p, \beta\}$. For the third scenario, the points of intersection are identified at $\{0,p,\beta\}$. Finally, in the fourth scenario, the intersection occurs at the interval $[0,p]$. The four scenarios are depicted in Figure \ref{fig:inf_pd}.
\color{black}

We now present four viable dual solutions that correspond to these scenarios. For scenario 1, the conditions $F(p)=0$, $F(\alpha)=1$ and $F'(\alpha)=0$ are observed. Moreover, this scenario has $F$ as a downward-sloping parabola. To align with these conditions, we set $\lambda_2 = 0$, as $\lambda_3 \leq 0$. This gives
\begin{flalign}
&&\lambda_0= -\frac{2 \alpha  p-p^2}{(\alpha -p)^2},\ \lambda_1= \frac{2 \alpha }{(\alpha -p)^2},\ \ \lambda_2 = 0,\ \lambda_3=-\frac{1}{(\alpha -p)^2}.&&\
\end{flalign}
The point $\alpha$ can be determined by considering a two-point distribution with support $\{p,\alpha\}$, mean $\mu$, and variance $\usigma^2$, yielding $\alpha= \mu +\usigma^2/(\mu-p)$.
Consequently,
\begin{flalign}
    && \lambda_0 + \lambda_1 \mu+\lambda_2 (\lsigma^2+\mu^2) + \lambda_3 (\usigma^2+\mu^2) = \frac{(\mu -p)^2}{(\mu-p )^2+\usigma^2}.&&
\end{flalign}

\noindent In scenario 2, we observe that $F(p)=0$ and $F(\beta)=1$. In addition, $F$ is characterized as a line, leading $\lambda_2=\lambda_3=0$. This gives dual solution
\begin{flalign}
&&\lambda_0= -\frac{p}{\beta-p},\ \lambda_1=\frac{1}{(\beta-p)},\ \lambda_2=0,\ \lambda_3=0,&&
\end{flalign}
and \begin{flalign}
    &&\lambda_0 + \lambda_1 \mu+\lambda_2 (\lsigma^2+\mu^2) + \lambda_3 (\usigma^2+\mu^2) = \frac{\mu-p}{\beta-p}.&&
\end{flalign}

\noindent In scenario 3 we observe that $F(0)=0$, $F(p)=0$, and $F(\beta)=1$. Moreover, $F$ is an upward-sloping parabola. Hence, in line with the functional form, we set $\lambda_3=0$, as $\lambda_2 \geq 0$. This results in
\begin{flalign}
&&\lambda_0= 0,\ \lambda_1= -\frac{p}{\beta(\beta-p)},\ \lambda_2=\frac{1}{\beta(\beta-p)}, \ \lambda_3 = 0,&&
\end{flalign}
and \begin{flalign}
    &&\lambda_0 + \lambda_1 \mu+\lambda_2 (\lsigma^2+\mu^2) + \lambda_3 (\usigma^2+\mu^2) = \frac{\mu ^2+\lsigma^2-\mu  p}{\beta (\beta-p)}.&&
\end{flalign}

\noindent Finally, in scenario 4, we note that $F(x)=0$, for all $x$ in the interval $[0,\beta]$. This gives $\lambda_0 = 0$, $\lambda_1 = 0$, $\lambda_2 = 0$, and $\lambda_3 = 0$. Hence,
\begin{flalign}
    &&\lambda_0 + \lambda_1 \mu+\lambda_2 (\lsigma^2+\mu^2) + \lambda_3 (\usigma^2+\mu^2) = 0.&&
\end{flalign}

Subsequently, we formulate four primal solutions, each with the same objective value as those obtained from the dual problem solutions. This equality in objective values asserts the claim of strong duality. A constructive method for deriving these primal solutions is now presented. Assuming strong duality holds, we have 
\begin{flalign}\label{compslack1}
    &&\int _x\1\{x > p\}{\rm d}\mathbb{P}^*(x) = \int _x(\lambda^*_0+\lambda^*_1\mu + \lambda^*_2   (\lsigma^2+\mu^2) +\lambda^*_3   (\usigma^2+\mu^2)  ){\rm d}\mathbb{P}^*(x).&&
\end{flalign}
Furthermore, due to dual feasibility, we obtain
\begin{flalign}\label{compslack2}
    &&\lambda^*_0+\lambda^*_1\mu + \lambda^*_2 (\lsigma^2+\mu^2) +\lambda^*_3 (\usigma^2+\mu^2) - \1\{x > p\}\geq 0,&&
\end{flalign} holding point-wise for each $x\in[0,\beta]$. Equation \eqref{compslack1} in conjunction with equation \eqref{compslack2} is known as the complementary slackness relationship in semi-infinite linear programming. A direct implication of this principle is that the atoms of the optimal primal solution coincide with the intersection of $\lambda^*_0+\lambda^*_1x + (\lambda^*_2+\lambda^*_3)x^2$ and the indicator function $\1\{x > p\}$. Following this, we will present four viable dual solutions corresponding to these scenarios.

In scenario 1, we construct a two-point distribution with atoms $\{p,\alpha\}$. The corresponding probabilities are
\begin{flalign}
    &&p_p=\frac{\usigma^2}{\mu ^2+\usigma^2+p^2-2 \mu  p}, \quad p_\alpha=1-p_p,&&
    \end{flalign}
leading to
\begin{flalign}
    &&\int_x\1\{x > p\}{\rm d}\mathbb{P}(x)=p_\alpha=\frac{(\mu -p)^2}{(\mu-p )^2+\usigma^2}.&&
    \end{flalign}
However, this distribution's feasibility requires $\alpha \leq \beta$, necessitating $p \in [0,\uupsilon_1]$.

\noindent For scenario 2, the construction involves a two-point distribution with atoms $\{p,\beta\}$. The probabilities are
\begin{flalign}
    &&p_p=\frac{\beta-\mu}{\beta-p}, \quad p_{\beta}=\frac{\mu-p}{\beta-p},&&
\end{flalign}
resulting in
\begin{flalign}
    &&\int_x\1\{x> p\}{\rm d}\mathbb{P}(x)=p_{\beta} = \frac{\mu-p}{\beta-p}.&&
\end{flalign}
However, the feasibility criterion necessitates
$\lsigma^2 \leq (p-\mu)^2p_p + (\beta-\mu)^2p_{\beta} \leq \usigma^2$,
requiring $p \in [\uupsilon_1,\lupsilon_1]$.

\noindent For scenario 3, we consider a three-point distribution with probability masses on the points $\{0,p,\beta\}$ and corresponding probabilities
\begin{flalign}
    &&p_0 = \frac{(\beta-\mu)(p-\mu)+\lsigma^2}{\beta p}, \quad p_p = 1-p_0-p_{\beta}, \quad p_{\beta}=\frac{\mu ^2+\lsigma^2-\mu  p}{\beta (\beta-p)}.&&
\end{flalign}
This yields
\begin{flalign}
    &&\int_x\1\{x> p\}{\rm d}\mathbb{P}(x)=p_{\beta} = \frac{\mu ^2+\lsigma^2-\mu  p}{\beta (\beta-p)}.&&
\end{flalign}
However, this distribution is only feasible when $p_0 \geq 0$ and $p_{\beta} \geq 0$, which requires $p \in [\lupsilon_1,\lupsilon_2]$.

\noindent Scenario 4 sees a two-point distribution with variance $\lsigma^2$ and atoms $\{\bar{\alpha},p\}$ with $\bar{\alpha} = \mu - \frac{\lsigma^2}{p-\mu}$. The corresponding probabilities are
\begin{flalign}
    &&p_{\bar{\alpha}} = 1 - p_p, \quad p_p = \frac{(p-\mu)^2}{(p-\mu)^2+\lsigma^2},&&
\end{flalign}
and hence
\begin{flalign}
    &&\int_x\1\{x > p\}{\rm d}\mathbb{P}(x) = 0.&&
\end{flalign}
  However, this distribution is only feasible when $\beta \geq 0$, requiring $p \in [\lupsilon_2,\beta]$. Since the primal and dual objective values are the same, we know that strong duality holds. Therefore, these objective values are optimal. In addition, the objective value is left-continuous in $p$. Therefore, we can apply Lemma \ref{lemma_pd}, stated below, and conclude that $\inf_{\bP \in \cP(\mu,\lsigma, \usigma,\beta)}\bP(X \geq p)=\inf_{\bP \in \cP(\mu,\lsigma, \usigma,\beta)}\bP(X > p).$
This completes the proof.

\begin{lemma}\label{lemma_pd}
Let $p \in (0,\beta]$ and $\cP$ any ambiguity set.\begin{enumerate}[label=(\roman*)]
  \item\label{claim1} If $\inf_{\bP\in\cP}\bP(X>p)$ is left-continuous in $p$, then $\inf_{\bP\in\cP}\bP(X \geq p) = \inf_{\bP\in\cP}\bP(X > p).$
  \item\label{claim2} If $\sup_{\bP\in\cP}\bP(X \geq p)$ is right-continuous in $p$, then $\sup_{\bP\in\cP}\bP(X \geq p) = \sup_{\bP\in\cP}\bP(X > p).$
\end{enumerate}
\end{lemma}
\begin{customproof}
Let $p\in(0,\beta]$ and $\epsilon>0$. The following identity trivially holds:
$\inf_{\bP\in\cP}\bP(X>p) \leq \inf_{\bP\in\cP}\bP(X\geq p) \leq \inf_{\bP\in\cP}\bP(X>p-\epsilon).$
However, since $\epsilon$ is chosen arbitrarily, we also have that
$\inf_{\bP\in\cP}\bP(X\geq p) \leq \lim_{\epsilon\xrightarrow{}0^+}\inf_{\bP\in\cP}\bP(X>p-\epsilon) = \lim_{\epsilon\xrightarrow{}p^-}\inf_{\bP\in\cP}\bP(X>\epsilon) = \inf_{\bP\in\cP}\bP(X>p).$
Hence, we obtain $\inf_{\bP\in\cP}\bP(X>p) \leq \inf_{\bP\in\cP}\bP(X\geq p)\leq \inf_{\bP\in\cP}\bP(X> p),$
which proves \ref{claim1}.

Next, observe that $\sup_{\bP\in\cP}\bP(X\geq p) \geq \sup_{\bP\in\cP}\bP(X > p) \geq \sup_{\bP\in\cP}\bP(X\geq p+\epsilon).$
However, since $\epsilon$ is chosen arbitrarily, we also have that
$\sup_{\bP\in\cP}\bP(X > p) \geq \lim_{\epsilon\xrightarrow{}0^+}\sup_{\bP\in\cP}\bP(X\geq p+\epsilon) = \lim_{\epsilon\xrightarrow{}p^+}\sup_{\bP\in\cP}\bP(X\geq\epsilon) = \sup_{\bP\in\cP}\bP(X\geq p).$
Hence, we obtain $\sup_{\bP\in\cP}\bP(X\geq p) \leq \sup_{\bP\in\cP}\bP(X > p)\leq \sup_{\bP\in\cP}\bP(X \geq p),$
which proves \ref{claim2}.
\end{customproof}

\subsection{Tight upper tail bounds}\label{ap:tutb}
\begin{lemma}\label{th:tutb}
Consider ambiguity set $\cP_{(\mu,\lsigma,\usigma, \beta)}$ and denote $\lupsilon_1 := \mu - \frac{\lsigma^2}{\beta-\mu}$, $\lupsilon_2 := \mu + \frac{\lsigma^2}{\mu}$ and $\uupsilon_2 := \mu + \frac{\usigma^2}{\mu}$. Then, 
\begin{align}
    &&\sup_{\bP \in \cP_{(\mu,\lsigma,\usigma, \beta)}}\bP(X \geq p)=\sup_{\bP \in \cP_{(\mu,\lsigma,\usigma, \beta)}}\bP(X > p)=
    \begin{cases}
     1, \quad & p\in(0,\lupsilon_1], \\
     \frac{(\beta+p)\mu-\mu^2-\lsigma^2}{\beta p}, \quad & p\in[\lupsilon_1,\lupsilon_2], \\
     \frac{\mu}{p},\quad & p\in[\lupsilon_2,\uupsilon_2], \\
     \frac{\usigma^2}{(\mu-p)^2+\usigma^2}, & p\in[\uupsilon_2,\beta].
    \end{cases}&&
\end{align}
\end{lemma}
\begin{customproof}
Let $\mathcal{M}$ be the set of probability measures defined on the measurable space $([0,\beta],\mathcal{B})$. We will solve $\sup_{\bP \in \cP(\mu,\lsigma,\usigma,\beta)}\bP(X \geq p)$, which can be rewritten as
\begin{equation}\label{primalsup}
\begin{aligned}
&\! \sup_{\mathbb{P}\in \mathcal{M}} &  &\int_x \1_{\{x\geq p\}}{\rm d} \mathbb{P}(x),\\
&\text{s.t.} &      &  \int_x {\rm d}\mathbb{P}(x)=1, \ \int_x x{\rm d}\mathbb{P}(x)=\mu,\ \int_x x^2{\rm d}\mathbb{P}(x)\geq\lsigma^2+\mu^2,\ \int_x x^2{\rm d}\mathbb{P}(x)\leq\usigma^2+\mu^2.
\end{aligned}
\end{equation}
Consider the corresponding dual, 
\begin{equation}\label{dualsup}
\begin{aligned}
&\inf_{\lambda_0,\lambda_1,\lambda_2,\lambda_3} &  &\lambda_0 + \lambda_1 \mu+\lambda_2 (\lsigma^2+\mu^2) + \lambda_3 (\usigma^2+\mu^2),\\
&\text{s.t.} &      & \1_{\{x\geq p\}}\leq \lambda_0  +\lambda_1 x+(\lambda_2+\lambda_3) x^2 =: F(x), \ \forall x\in[0,\beta], \lambda_2 \leq 0, \lambda_3 \geq 0.
\end{aligned}
\end{equation}
The proof is structured analogously to the proof of Lemma \ref{pdinf}. Hence, we derive both feasible dual and primal solutions and show that strong duality holds through matching objective values. Four scenarios are considered. In the first scenario, $F(x)$ intersects with $\1\{x\geq p\}$ within the interval $[p,\beta]$. In the second scenario, the intersection occurs in the points $\{0,p,\beta\}$. For the third scenario, the points of intersection are identified at $\{0,p\}$. Finally, in the fourth scenario, the intersections are at $\{\alpha,p\}$, where $\alpha \in [0,p]$. We now present four viable dual solutions corresponding to these scenarios.

For scenario 1 we observe that $F(x) = 1$ for all $x$ in the interval $[0,\beta]$. This gives $\lambda_0=1, \lambda_1=0, \lambda_2=0$, and $\lambda_3=0$. Hence,
\begin{flalign}
    &&\lambda_0 + \lambda_1 \mu+\lambda_2 (\lsigma^2+\mu^2) + \lambda_3 (\usigma^2+\mu^2) = 1.&&
\end{flalign}

\noindent In scenario 2, the conditions $F(0)=0$, $F(p)=1$, and $F(\beta)=1$ are observed. Furthermore, this scenario is characterized by $F$ taking the form of a downward-sloping parabola. To align with these conditions, we set $\lambda_3 = 0$, as $\lambda_2 \leq 0$. This gives
\begin{flalign}
&&\lambda_0= 0,\ \lambda_1= \frac{1}{\beta}+\frac{1}{p},\ \lambda_2=-\frac{1}{\beta p},\ \lambda_3=0,\ &&
\end{flalign}
and \begin{flalign}
    &&\lambda_0 + \lambda_1 \mu+\lambda_2 (\lsigma^2+\mu^2) + \lambda_3 (\usigma^2+\mu^2) = \frac{(\beta+p)\mu-\mu^2-\lsigma^2}{\beta p}.&&
\end{flalign}

\noindent In scenario 3 we observe that $F(0)=0$ and $F(p)=1$. Moreover, $F$ is linear. Hence, we set $\lambda_0 = \lambda_2 = \lambda_3 = 0$, resulting in
\begin{flalign}
&&\lambda_0= 0,\ \lambda_1= \frac{1}{p},\ \lambda_2=0,\ \lambda_3=0,\ &&
\end{flalign}
and \begin{flalign}
   &&\lambda_0 + \lambda_1 \mu+\lambda_2 (\lsigma^2+\mu^2) + \lambda_3 (\usigma^2+\mu^2) = \frac{\mu}{p}.&&
\end{flalign}

\noindent Finally, in scenario 4, we note that $F(p)=1$, $F(\alpha)=0$, and $F'(\alpha)=0$. In this scenario, the functional form of $F$ is an upward-sloping parabola. Hence, following these specified conditions, we set $\lambda_2 = 0$, as $\lambda_3 \geq 0$. This gives dual solution
\begin{flalign}
&&\lambda_0= \frac{\alpha^2}{(\alpha-p)^2},\ \lambda_1= -\frac{2\alpha}{(\alpha-p)^2},\ \lambda_2=0,\ \lambda_3=\frac{1}{(\alpha-p)^2}.&&
\end{flalign}
The point $\alpha$ can be determined by considering a two-point distribution with support $\{\alpha,p\}$, mean $\mu$, and maximal variance $\usigma^2$, giving 
$$
\alpha= \mu-\frac{\usigma^2}{p-\mu},
$$
and
\begin{flalign}
    &&\lambda_0 + \lambda_1 \mu+\lambda_2 (\lsigma^2+\mu^2) + \lambda_3 (\usigma^2+\mu^2) = \frac{\usigma^2}{(\mu-p)^2+\usigma^2}.&&
\end{flalign}

Subsequently, we formulate four primal solutions, each with the same objective value as those obtained from the dual problem solutions. This equality in objective values asserts the claim of strong duality. A constructive method for deriving these primal solutions is now presented. Assuming strong duality holds, we have 
\begin{flalign}\label{compslack3}
    &&\int _x\1\{x \geq p\}{\rm d}\mathbb{P}^*(x) = \int _x(\lambda^*_0+\lambda^*_1\mu + \lambda^*_2   (\lsigma^2+\mu^2) +\lambda^*_3   (\usigma^2+\mu^2)  ){\rm d}\mathbb{P}^*(x).&&
\end{flalign}
Furthermore, due to dual feasibility, we obtain
\begin{flalign}\label{compslack4}
    &&\lambda^*_0+\lambda^*_1\mu + \lambda^*_2 (\lsigma^2+\mu^2) +\lambda^*_3 (\usigma^2+\mu^2) - \1\{x\geq p\}\geq 0,&&
\end{flalign} holding point-wise for each $x\in[0,\beta]$. Equation \eqref{compslack3} in conjunction with equation \eqref{compslack4} is known as the complementary slackness relationship in semi-infinite linear programming. A direct implication of this principle is that the atoms of the optimal primal solution coincide with the intersection of $\lambda^*_0+\lambda^*_1x + (\lambda^*_2+\lambda^*_3)x^2$ and the indicator function $\1\{x \geq p\}$.

In scenario 1, the construction involves a two-point distribution with atoms $\{p,\bar{\alpha}\}$, where $\bar{\alpha} := \mu + \frac{\lsigma^2}{\mu-p}$. The corresponding probabilities are
\begin{flalign}
    &&p_p = \frac{\lsigma^2}{(p-\mu)^2+\lsigma^2}, \quad p_{\bar{\alpha}} = 1 - p_p,&&
\end{flalign}
leading to
\begin{flalign}
    &&\int_x\1\{x \geq p\}{\rm d}\mathbb{P}(x) = 1.&&
\end{flalign}
However, this distribution's feasibility requires $\bar{\alpha} \leq \beta$, necessitating $p \in [0,\lupsilon_1]$.\\

\noindent For scenario 2, a unique three-point distribution is considered with probability masses at $\{0,p,\beta\}$. The probabilities are
\begin{flalign}
    &&p_0 = 1-\frac{(\beta+p)\mu-\mu^2-\lsigma^2}{\beta p}, \quad p_p = 1-p_0-p_\beta, \quad p_\beta=\frac{\lsigma^2+\mu^2-p\mu}{\beta(\beta-p)},&&
\end{flalign}
resulting in
\begin{flalign}
    &&\int_x\1\{x\geq p\}{\rm d}\mathbb{P}(x)=p_p + p_\beta = \frac{(\beta+p)\mu-\mu^2-\lsigma^2}{L p}.&&
\end{flalign}
However, this distribution is only feasible when $p_0 \geq 0$ and $p_\beta \geq 0$, requiring $p \in [\lupsilon_1,\lupsilon_2]$.\\

\noindent For scenario 3, we consider a two-point distribution on $\{0,p\}$ with probabilities 
\begin{flalign}
    &&p_0=1-p_p, \quad p_p=\frac{\mu}{p}.&&
\end{flalign}
This yields
\begin{flalign}
    &&\int_x\1\{x \geq p\}{\rm d}\mathbb{P}(x)=p_p = \frac{\mu}{p}.&&
\end{flalign}
The feasibility criterion necessitates
\begin{flalign}
    &&\lsigma^2 \leq \mu^2p_0 + (p-\mu)^2p_p \leq \usigma^2,&&
\end{flalign}
requiring $p \in [\lupsilon_2,\uupsilon_2]$.\\

\noindent Scenario 4 sees a two-point distribution with probability masses on the points $\{\alpha,p\}$ with probabilities
\begin{flalign}
    &&p_\alpha=1-p_p, \quad p_p=\frac{\usigma^2}{(\mu-p)^2+\usigma^2},&&
    \end{flalign}
leading to
\begin{flalign}
    &&\int_x\1\{x \geq p\}{\rm d}\mathbb{P}(x)=p_p=\frac{\usigma^2}{(\mu-p )^2+\usigma^2}.&&
    \end{flalign}
However, this distribution is only feasible when $\alpha \leq \beta$, requiring $p \in [\uupsilon_2,\beta]$. Since the primal and dual objective values are the same, we know that strong duality holds. Therefore, these objective values are optimal. Furthermore, the objective value is right-continuous in $p$. Therefore, we can apply Lemma \ref{lemma_pd} and conclude that $\sup_{\bP \in \cP(\mu,\lsigma, \usigma,\beta)}\bP(X \geq p)=\sup_{\bP \in \cP(\mu,\lsigma, \usigma,\beta)}\bP(X > p).$
\end{customproof}

\section{Comparison with optimal pricing}
As is common in the pricing literature, we will now evaluate the performance of the robust pricing strategy with the optimal Bayesian strategy; see e.g.~\cite{azar2012optimal,chen2022distribution}. The optimal Bayesian strategy is known to be fixed pricing \citep{riley1983optimal} and achieves the optimal expected revenue
$\textup{OPT}(\bP) \eqdef \sup_p p\bP(X\geq p)$. Markov's inequality yields $p\bP(X \geq p) \leq \mu$, and hence $\textup{OPT}(\bP)\leq \mu$. This is in fact a tight bound. To see this, let $p \in [\mu+{\lsigma^2}/{\mu},\mu+{\usigma^2}/{\mu}]$ and consider the two-point distribution $\bP_0\in \cP(\mu,\lsigma,\usigma,\beta)$ with atoms $\{0,p\}$ and associated probabilities $1-\mu/p$ and ${\mu}/{p}$, respectively, for which $p\bP_0(X\geq p) = \mu$. Combining this bound with Theorem~\ref{3pthm} yields, upon some rewriting, the following performance guarantee for robust pricing:
\begin{corollary}\label{cor:cpx}
    For any $\bP \in \cP(\mu,\lsigma,\usigma,\beta)$ consider the optimal price $p^*$ in \eqref{p^*_gmain}.
    Then,
    \begin{flalign}\label{cpx}
    &&\frac{p^*\bP(X\geq p^*)}{\textup{OPT}(\bP)} \geq \begin{cases}
    \rho_l(\mu,\usigma) \eqdef \frac{1}{2}\frac{\usigma}{\mu}\kappa^3(\mu,\usigma), \quad & (\lsigma,\usigma) \in \Sigma_l, \\
    \rho_m(\mu,\beta) \eqdef \frac{2}{\mu}(\beta-\sqrt{\beta(\beta-\mu)})-1, \quad & (\lsigma,\usigma) \in \Sigma_m, \\
    \rho_h(\mu,\lsigma,\beta) \eqdef 2(1-\sqrt{1-\frac{\mu^2+\lsigma^2}{\mu\beta}})-\tfrac{\mu^2+\lsigma^2}{\mu\beta}, \quad & (\lsigma,\usigma) \in \Sigma_h.
    \end{cases}&&
    \end{flalign}
\end{corollary}
\begin{figure}[h]
    \centering
    \begin{tikzpicture}
    \node at (0,0) {\includegraphics[width=0.53\linewidth]{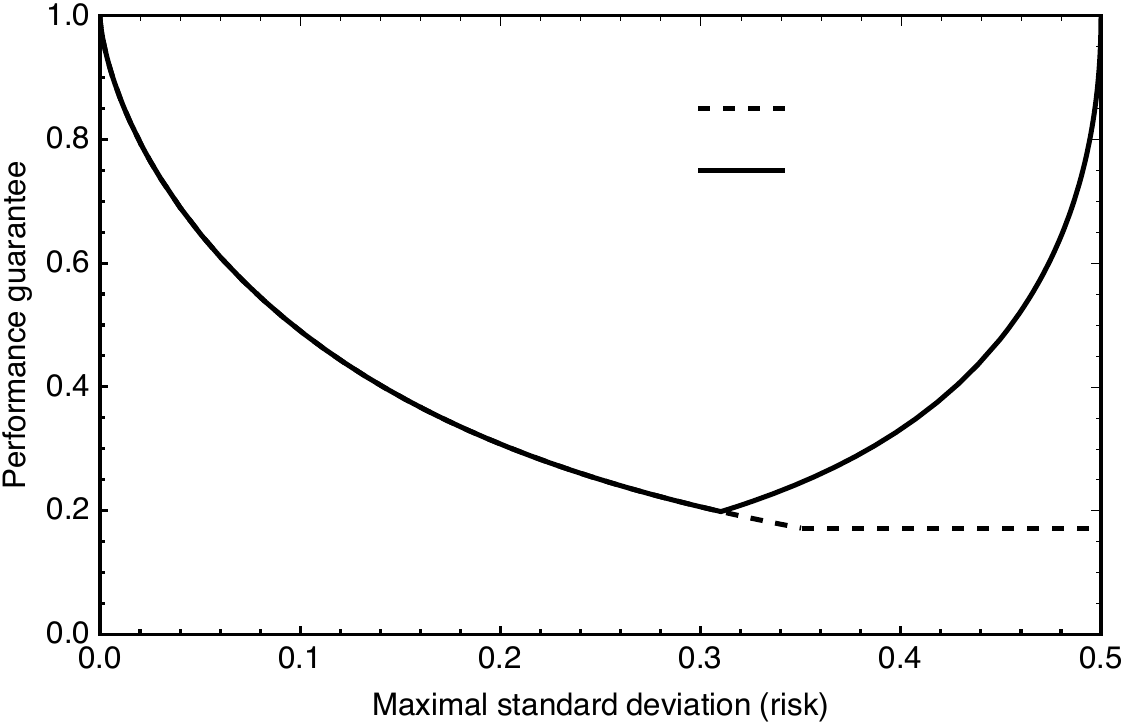}};
    \node at (2.15,1.42)   {\begin{scriptsize}$\lsigma = \usigma$\end{scriptsize}};
    \node at (2.15,1.85)   {\begin{scriptsize}$\lsigma = 0$\end{scriptsize}};
    \end{tikzpicture}
    \caption{Corollary~\ref{cor:cpx} with $\mu=0.5$ and $\beta=1$. }
    \label{fig:cpx}
\end{figure}    

\noindent{\bf U-shaped performance guarantee.}
Figure~\ref{fig:cpx} illustrates Corollary \ref{cor:cpx}. When the seller has precise information about the variance, there occurs a U-shaped performance guarantee, aligned with the U-shaped revenue curves found in monopolistic markets \citep{johnson2006simple,sun2012does} discussed in the introduction. In our case, this U-shape occurs as the seller switches from low to high pricing, changing the curve from being decreasing to increasing. More generally, we have found that the U-shape occurs when the seller has precise enough (not perfect per se) information on the variance. When the seller only knows an upper bound, the U-shape is replaced with a function that first decreases and then remains constant, see Figure~\ref{fig:cpx}. We substantiate these numerical insights with the following result:
\begin{corollary}[Properties of price and revenue]\label{cor:mon} 
    {\rm (i)} $p^*_l(\usigma)$ and  $\rho_l(\mu,\usigma)$ increase in $\mu$ and decrease in $\usigma$;  {\rm (ii)} $p^*_m(\beta)$  and $\rho_m(\mu,\beta)$ increase in $\mu$ and decrease in $\beta$;
    {\rm (iii)} $p^*_h(\lsigma,\beta)$ and $\rho_h(\mu,\lsigma,\beta)$ increase in $\lsigma$ and decrease in $\beta$.
\end{corollary}

\vspace{.2cm}
\noindent{\bf Risk-aligned pricing.} 
The optimal price and performance guarantee adjust in the same direction in response to changes in mean, support, and variance. In uncertain markets, robust pricing involves higher prices to reflect risk levels. Higher volatility leads to a high-price strategy, targeting a smaller, high-value customer segment. Corollary \ref{cor:mon} shows that only the high price is positively correlated with variance, resulting in a performance guarantee that increases with variance. Managers in volatile environments should focus on value alignment rather than broad appeal, catering to customers willing to pay for exclusivity, as seen in luxury fashion and high-end tech. Conversely, Corollary \ref{cor:mon} shows that low prices are negatively correlated with variance. For example, electric vehicles initially targeted a niche market with high prices due to production and demand uncertainties. As technology matured and demand increased, prices dropped, attracting mainstream consumers and a larger customer base.

\section{Remaining proofs of  pricing results}\label{sec:proofs}

\subsection{Proof of  Lemma~\ref{middlepointnot}}
We first set $\mu = 1$ without loss of generality. Recall $\hat{R}(p)$ from \eqref{rev_h} and define $\hat{R}(p;\mu,\sigma,\beta) \eqdef \hat{R}(p)$. In addition, consider the scaled parameters $\nu = \frac{\sigma}{\mu}$, $\tau = \frac{\beta}{\mu}$, and $\psi = \frac{p}{\mu}$. From \eqref{rev_h} it follows that $\hat{R}(p;\mu,\sigma,\beta) = \hat{R}(\psi\mu;\mu,\nu\mu,\tau\mu) = \mu \cdot \hat{R}(\psi;1,\nu,\tau).$
Hence, using the substitution $p = \psi\mu$, we have $\arg\sup_p \hat{R}(p;\mu,\sigma,\beta) = \mu \arg\sup_\psi \hat{R}(\psi;1,\nu,\tau).$ 
We proceed with $\mu = 1$. 

Consider that for $\lupsilon_1$ to be the unique global optimum, the following ordering needs to hold: 
\begin{flalign}\label{ordering}
    &&p^*_{h}(\sigma,\beta) < \lupsilon_1 < p^*_{l}(\sigma).&&
\end{flalign}
We will show that this ordering never occurs by considering two cases that exhaust all possibilities. The first case considers $\sigma \geq 1$. We claim the following:
\begin{flalign}
    &&p^*_{h}(\sigma,\beta) \geq 1 \geq \lupsilon_1.&&
\end{flalign}
\noindent The second inequality is trivial. To see why the first inequality holds true, note that
\begin{flalign}
    &&\beta - \sqrt{\beta(\beta-(1+\sigma^2))} \geq 1
    &\Leftrightarrow \beta - 1 \geq \sqrt{\beta^2-\beta-\beta\sigma^2}&& \nonumber \\
    && &\Leftrightarrow \beta^2 -2\beta+1 \geq \beta^2-\beta-\beta\sigma^2&& \nonumber \\
    && &\Leftrightarrow 1+\beta\sigma^2 \geq \beta.&&
\end{flalign}
This is true since both $\sigma^2 \geq 1$ and $\beta \geq 1$. 

The second case considers $\sigma \in (0,1)$. Denote
\begin{flalign}
    &&\beta_l \eqdef \frac{1}{2} \left(2+\sigma^2+A^{\frac{1}{3}}+\frac{\sigma^4}{A^{\frac{1}{3}}} \right) \text{with } A\eqdef\sigma^4(\sqrt{1+\sigma^2}+1)^2,&&
\end{flalign}
as the unique value of $\beta$ for which $p^*_{l}(\sigma) = \lupsilon_1$,
and
\begin{flalign}
    &&\beta_h \eqdef \frac{2-2^{\frac{4}{3}}\sigma^4 B^{-\frac{1}{3}}-2^\frac{2}{3}B^{\frac{1}{3}}}{2(1-\sigma^2)} \text{with } B\eqdef-2\sigma^8,&&
\end{flalign}
as the unique value for which $p^*_{h}(\sigma,\beta) = \lupsilon_1$. For \eqref{ordering} to hold, it has to be the case that $\beta_h < \beta_l$. Hence, we will show that the opposite, namely $\beta_h \geq \beta_l$ is always true. It is sufficient to demonstrate that $\beta_h \geq (1-\sigma^2)\beta_h \geq \beta_l$. This results in the following inequality:
\begin{flalign}\label{in1}
    &&2\sigma^{\frac{4}{3}}+2\sigma^{\frac{8}{3}} \geq \sigma^2 + \sigma^\frac{4}{3}(\sqrt{1+\sigma^2}+1)^{\frac{2}{3}} + \frac{\sigma^{\frac{8}{3}}}{(\sqrt{1+\sigma^2}+1)^{\frac{2}{3}}}.&&
\end{flalign}
\noindent Dividing by $\sigma^{\frac{4}{3}}$ and substituting $c$ for $\sigma^\frac{2}{3}$ results in
$2+2c^2 \geq c + (\sqrt{1+c^3}+1)^\frac{2}{3} + c^2/((\sqrt{1+c^3}+1)^\frac{2}{3})$,
which is the same as 
\begin{flalign}\label{eq1.8}
    &&2+2c^2 - c \geq (\sqrt{1+c^3}+1)^{-\frac{2}{3}}\left[(\sqrt{1+c^3}+1)^{\frac{4}{3}}+c^2\right].&&
\end{flalign}
Since $\sqrt{1+c^3}+1 \geq 2$ for every $c \in (0,1)$, we have
\begin{flalign}\label{eq1.9}
    &&\frac{2}{3} \geq 2^{-\frac{2}{3}} \geq (\sqrt{1+c^3}+1)^{-\frac{2}{3}}.&&
\end{flalign}
Hence, combining \eqref{eq1.8} with \eqref{eq1.9} results in $2+2c^2 - c \geq \frac{2}{3}\left[(\sqrt{1+c^3}+1)^{\frac{4}{3}}+c^2\right],$
which is equivalent to
\begin{flalign}\label{eq1.11}
    &&3+2c^2-\frac{3}{2}c \geq (\sqrt{1+c^3}+1)^\frac{4}{3}.&&
\end{flalign}
Denote the left-hand side as $r(c) \eqdef 3+2c^2-\frac{3}{2}c$ and the right-hand side as $t(c) \eqdef (\sqrt{1+c^3}+1)^\frac{4}{3}$. Note that the quadratic function $r(c)$ has one local minimum at $c^* = \frac{3}{8}$ and that $t(c)$ is strictly increasing in $c$. By inspection, one can verify that \eqref{eq1.11} is true for both $c = 0$ and $c = c^*$. Since $r(c)$ is decreasing on the interval $(0,c^*]$, and $t(c)$ increasing, it follows that \eqref{eq1.11} is true for all $c \in (0,c^*]$. Both functions are increasing on the interval $(c^*,1]$. We will argue that it suffices to check \eqref{eq1.11} for a finite number of points to conclude that it holds for all $c \in (c^*,1)$. 

Assume that for a given $\epsilon>0$ and $d,d+\epsilon \in (c^*,1]$, we have $r(d) > t(d+\epsilon)$. Then \eqref{eq1.11} holds for any $c \in [d,d+\epsilon]$, as
$r(c) \geq r(d) \geq t(d+\epsilon) \geq t(c),$
since $r$ and $t$ are both continuous and increasing. Therefore, we can make a discretization of the interval $(c^*,1]$ formed by $K$ points $c^* = c_0 < c_1 < ... < c_K = 1$ with $c_k = c_0 + \epsilon * k$, for all $k \in \{1,...,K\}$. As it turns out, $K=6$ is the smallest number for which this method works and has corresponding $\epsilon = \frac{1-c^*}{6} = \frac{5}{48}$. Table \ref{table_calc} shows the function evaluations that together prove \eqref{eq1.11}. This proves that \eqref{ordering} never occurs when $\sigma \in (0,1)$. Hence, $\lupsilon_1$ is never the unique global optimum.

\begin{table}[H]
\centering
\renewcommand{\arraystretch}{1.2} % Triple row height
\begin{tabular}{>{\centering\arraybackslash}p{2.1cm}|>{\centering\arraybackslash}p{2.1cm}|>{\centering\arraybackslash}p{2.1cm}}
\hline
\textbf{$c_k$} & \textbf{$r(c_k) \approx$} & \textbf{$t(c_k + \epsilon) \approx$} \\
\hline
$c_0 = \frac{3}{8}$ & 2.72 & 2.61 \\
\hline
$c_1 = \frac{23}{48}$ & 2.74 & 2.68 \\
\hline
$c_2 = \frac{7}{12}$ & 2.81 & 2.78 \\
\hline
$c_3=\frac{11}{16}$ & 2.91 & 2.90 \\
\hline
$c_4 = \frac{19}{24}$ & 3.07 & 3.06 \\
\hline
$c_5 = \frac{43}{48}$ & 3.26 & 3.24 \\
\hline
\end{tabular}
\caption{\textnormal{Numerical calculations required in the proof of \eqref{eq1.11}}.}
\label{table_calc}
\end{table}

\subsection{Proof of Corollary~\ref{cor:mon}}
With shorthand notation $r \eqdef \frac{\mu}{\sigma}$ and $r_+ \eqdef \sqrt{1+(\frac{\mu}{\sigma})^2}>r$, observe that
    \begin{flalign}
        &&\frac{\partial \rho_l(\mu,\usigma)}{\partial r} = \dfrac{\left(\sqrt[3]{r+r_+}+\sqrt[3]{r-r_+}\right)^2\left(\left(r+r_+\right)^\frac{2}{3}-\left(r-r_+\right)^\frac{2}{3}\right)}{2r^2r_+\left(r-r_+\right)^\frac{2}{3}\left(r+r_+\right)^\frac{2}{3}} \geq 0,&&
    \end{flalign}
    as $r>0$ and
    \begin{flalign}
        &&\left(r+r_+\right)^\frac{2}{3}-\left(r-r_+\right)^\frac{2}{3} \geq 0.&&
    \end{flalign}
    Hence, $\rho_l(\mu,\usigma)$ increases in $\mu$ and decreases in $\usigma$. 
    Next, denote $s \eqdef \frac{\beta}{\mu}$ and observe that
    \begin{flalign}
        &&\frac{\partial \rho_m(\mu,\beta)}{\partial s} = 2\left(1-\dfrac{2s-1}{2\sqrt{\left(s-1\right)s}}\right) \leq 0.&&
    \end{flalign}
    Therefore, $\rho_m(\mu,\beta)$ increases in  $\mu$ and decreases in $\beta$. 
    With $q \eqdef \frac{\mu^2+\lsigma^2}{\mu\beta}$ observe that 
    \begin{flalign}
        &&\frac{\partial \rho_h(\mu,\lsigma,\beta)}{\partial q} = \dfrac{1}{\sqrt{1-q}}-1 \geq 0,&&
    \end{flalign}
    since $q = \frac{\lupsilon_2}{\beta} < 1$. This shows that $\rho_h(\mu,\lsigma,\beta)$ increases in $\lsigma$ and decreases in $\beta$.
     Next, recall that $p^*_l(\usigma)$ is the unique real interior solution of $R_1(p)$, so solving
    \begin{flalign}\label{eq:m1}
        &&(p^*_l(\usigma))^3-3\mu (p^*_l(\usigma))^2+3(\mu^2+\usigma^2)p^*_l(\usigma) = \mu(\mu^2 + \usigma^2).&&
    \end{flalign}
    To simplify the analysis, we change variables by substituting $\mu_2 \eqdef \mu^2 + \usigma^2$ (the maximal second moment), and write \eqref{eq:m1} as
    \begin{flalign}\label{eq:m2}
        &&(p^*_l(\mu_2-\mu))^3-3\mu (p^*_l(\mu_2-\mu))^2+3\mu_2p^*_l(\mu_2-\mu) = \mu\mu_2.&&
    \end{flalign}
    By keeping $p^*_l(\mu_2-\mu)$ implicit while isolating $\mu$ and $\mu_2$ in \eqref{eq:m2}, we can write \eqref{eq:m1} as
    \begin{flalign}\label{eq:m3}
        &&L_1(p^*_l(\mu_2-\mu),\mu_2) \eqdef \frac{(p^*_l(\mu_2-\mu))^3+3p^*_l(\mu_2-\mu)\mu_2}{3(p^*_l(\mu_2-\mu))^2+\mu_2} = \mu,&&
    \end{flalign}
    or
    \begin{flalign}\label{eq:m4}
        &&L_2(p^*_l(\mu_2-\mu),\mu) \eqdef \frac{-(p^*_l(\mu_2-\mu))^2(p^*_l(\mu_2-\mu)-3\mu)}{3p^*_l(\mu_2-\mu)-\mu} = \mu_2.&&
    \end{flalign}
    Differentiating both sides of \eqref{eq:m3} with respect to $\mu$ yields
    \begin{flalign}\label{eq:m5}
         && \frac{\partial L_1(p^*_l(\mu_2-\mu),\mu_2)}{\partial p^*_l(\mu_2-\mu)} \cdot \frac{\partial p^*_l(\mu_2-\mu)}{\partial \mu} = \frac{3((p^*_l(\mu_2-\mu)))^2-\mu_2)^2}{(3p^*_l(\mu_2-\mu)^2+\mu_2)^2} \cdot \frac{\partial p^*_l(\mu_2-\mu)}{\partial \mu}= 1,&&
    \end{flalign}
    showing that $p^*_l(\usigma)$ is (strictly) increasing in $\mu$. Differentiating both sides of \eqref{eq:m4} with respect to $\mu_2$ gives
    \begin{flalign}\label{eq:m6}
         &&\frac{\partial L_2(p^*_l(\mu_2-\mu),\mu)}{\partial p^*_l(\mu_2-\mu)} \cdot \frac{\partial p^*_l(\mu_2-\mu)}{\partial \mu_2} = \frac{-6p^*_l(\mu_2-\mu)(p^*_l(\mu_2-\mu)-\mu)^2}{(3p^*_l(\mu_2-\mu)-\mu)^2} \cdot \frac{\partial p^*_l(\mu_2-\mu)}{\partial \mu_2}=1,&&
    \end{flalign}
    showing that $p^*_l(\usigma)$ is (strictly) decreasing in $\usigma$. To see why $p^*_m(\beta)$ is increasing in $\mu$ and decreasing in $\beta$, we examine $\rho_m(\mu,\beta) = {2p^*_m(\beta)}/{\mu}-1$, which is clearly increasing in $p^*_m(\beta)$. Hence, the same monotonicity properties directly apply to $p^*_m(\beta)$. Finally, $p^*_h(\lsigma,\beta)$ increases in $\lsigma$ and decreases in $\beta$, due to $p^*_h(\lsigma,\beta) = p^*_m(\mu+{\lsigma^2}/{\mu},\beta)$ and monotonicity of both $\mu+{\lsigma^2}/{\mu}$ and $p^*_m(\beta)$.

\section{Proofs related to queueing models}
\subsection{Proof of Theorem \ref{th:queue}}\label{ap:queue1}
Recall that $W$ is an increasing function and so if $x < x'$, then also $p + hW(x) < p + hW(x')$. This implies that $\lambda \cdot f(x)$ is non-increasing in $x$. Since $x$ is increasing in $x$, the equation $x = \lambda \cdot \inf_{\Prob \in \cP} \Prob(X > p + hW(x))$ has a unique solution $x^*$. 
Recall that $\gamma^* = \inf_{\Prob \in \bP} \gamma(p,\Prob)$ with $\gamma(p,\Prob)$ the unique solution to the equation $\gamma = \lambda \cdot \Prob(X > p + hW(\gamma))$
for the given distribution $\bP$. Due to the second assumption we know that there is a distribution $\Prob_{x^*}$ so that $x^* = \lambda \cdot \inf_{\Prob \in \cP} \Prob(X > p + hW(x^*)) = \lambda \cdot \Prob_{x^*}(X > p + hW(x^*)).$
This means that  $\gamma^* \leq x^* < \infty$. 
We will next show that $\gamma^* = x^*$ using a proof by contradiction. Suppose that equality does not hold, i.e., $x^* > \gamma^*$. Then there exists a distribution $\Prob' \in \cP$ where $\gamma' = \gamma'(p,\Prob')$  satisfies
$\gamma' = \lambda \cdot \Prob'(X > p + hW(\gamma'))$ 
and $\gamma^* \leq \gamma' < x^*$. Because the function $x \mapsto \Prob'(X> p + hW(x))$ is non-increasing, it follows that 
$\Prob'(X> p + hW(x^*)) \leq \Prob'(X> p + hW(\gamma')) = \frac{\gamma'}{\lambda} <  \frac{x^*}{\lambda} = \Prob_{x^*}(X> p + hW(x^*)).$
By definition of $\Prob_{x^*}$, i.e., the second assumption of the theorem statement, this then implies that
$\Prob'(X> p + hW(x^*))  < \inf_{\Prob \in \cP} \Prob(X > p + hW(x^*))$
which is a contradiction by definition of the infimum. This shows that $x^* = \gamma^*$ and completes the proof.

\subsection{Proof of Theorem \ref{THMswithQ}}\label{proofQQs}
We establish Theorem \ref{THMswithQ} by combining Lemma \ref{QM1_sigma_pricejump} and \ref{QM1_beta_pricejump}, which we will provide below. To begin, we define key quantities that will be used throughout the proofs. First, we recall the worst-case tail bound from Lemma \ref{pdinf}, where we set $\lsigma = \usigma = \sigma$:
\begin{flalign}\label{P_i}
    &&P(t,\sigma,\beta) := \begin{cases}
        P_1(t,\sigma,\beta) :=\frac{(\mu-t)^2}{(\mu-t)^2+\sigma^2}, \quad & t \in (0,\upsilon_1], \\
     P_2(t,\sigma,\beta):=\frac{\mu(\mu-t)+\sigma^2}{\beta(\beta-t)}, \quad & t \in [\upsilon_1,\upsilon_2].
    \end{cases}&&
    \end{flalign}

\noindent By utilizing \eqref{P_i}, we define worst-case expected revenue as
\begin{flalign}\label{R_i}
    &&R(p,\sigma,\beta) := \begin{cases}
        R_1(p,\sigma,\beta) := p\lambda P_1\left(p+hW(\gamma_1^*),\sigma,\beta \right), \quad & p+hW(\gamma_1^*) \in (0,\upsilon_1], \\
     R_2(p,\sigma,\beta):=p\lambda P_2\left(p+hW(\gamma_2^*),\sigma,\beta\right), \quad & p +hW(\gamma_2^*) \in [\upsilon_1,\upsilon_2],
    \end{cases}&&
    \end{flalign}
where $\gamma^*_i = \gamma^*_i(p,\sigma,\beta)$ is the unique solution to $$\gamma^*_i = \lambda P_i(p + hW(\gamma^*_i),\sigma,\beta), \text{ for } i = 1, 2.$$ 

\noindent Denote $p^* := \arg\sup_pR(p,\sigma,\beta)$ and the two potential optimal prices as $$p^*_l(\sigma,\beta):=\arg\sup_{p+hW(\gamma^*_1)\in(0,\upsilon_1]}R_1(p,\sigma,\beta), \ \ p^*_h(\sigma,\beta) := \arg\sup_{p+hW(\gamma^*_2)\in[\upsilon_1,\upsilon_2]}R_2(p,\sigma,\beta).$$
For clarity and brevity, we occasionally adopt a slight abuse of notation by omitting either $\sigma$ or $\beta$ as a function argument when the variable is not relevant in a given proof. This should be understood as a notational simplification with the omitted variable remaining fixed and the function itself unchanged.
\begin{lemma}\label{QM1_sigma_pricejump}
    A unique $\sigma^* \in (0,\sigma_{\textup{max}})$ exists such that $p^* = p^*_l(\sigma)$ when $\sigma \leq \sigma^*$ and $p^* = p^*_h(\sigma)$ when $\sigma \geq \sigma^*$.
\end{lemma}
\begin{customproof}
Let $p$ be fixed and consider $\sigma_2 > \sigma_1 \geq 0$. Now let $x$ be the unique solution to $x = \lambda P_1(p + h W(x),\sigma)$. Clearly, for fixed $x$, we have $\lambda P_1(p + h W(x),\sigma_1) > \lambda P_1(p + h W(x),\sigma_2)$. Hence, $\gamma_1^*(p,\sigma_1) > \gamma_1^*(p,\sigma_2)$. But then also $p\gamma^*_1(p,\sigma_1) > p\gamma^*_1(p,\sigma_2)$. Since $p$ was chosen arbitrarily, the following holds:
$$R_1(p^*_l(\sigma_1),\sigma_1) > R_1(p^*_l(\sigma_2),\sigma_1) > R_1(p^*_l(\sigma_2),\sigma_2).$$
Hence, $R_1(p^*_l,\sigma)$ is decreasing in $\sigma$.

Now let $y$ be the unique solution to $y = \lambda P_2(p + h W(y),\sigma)$. Clearly, for fixed $y$, we have $\lambda P_2(p + h W(y),\sigma_1) < \lambda P_2(p + h W(y),\sigma_2)$. Hence, $\gamma_2^*(p,\sigma_1) < \gamma_2^*(p,\sigma_2)$. But then also $p\gamma^*_2(p,\sigma_1) < p\gamma^*_2(p,\sigma_2)$. Since $p$ was chosen arbitrarily, the following holds:
$$R_2(p^*_h(\sigma_1),\sigma_1) < R_2(p^*_h(\sigma_1),\sigma_2) < R_2(p^*_h(\sigma_2),\sigma_2).$$
Hence, $R_2(p^*_h,\sigma)$ is increasing in $\sigma$.

Now consider $\sigma = 0$. Then $p^*_h(0) + hW(\gamma^*_2) = \mu$, as $\upsilon_1 = \upsilon_2 = \mu$, and consequently, $R_2(p^*_h(0),0) = 0$. Hence, since $\sup_{p+hW(\gamma^*_1) \in (0,\mu]}R_1(p,0)>0$, we have that $p^* = p^*_l(0)$. We now consider $\sigma = \sigma_{\textup{max}} = \sqrt{\mu(\beta-\mu)}$. Then $p^*_l(\sigma_{\textup{max}}) + hW(\gamma^*_1) = 0$, as $\upsilon_1 \xrightarrow{} 0^+$, and consequently, $R_1(p^*_l(\sigma_{\textup{max}}),\sigma_{\textup{max}}) = 0$. Hence, as $\sup_{p+hW(\gamma^*_2) \in (0,\beta]}R_2(p,\sigma_{\textup{max}})>0$, we have that $p^* = p^*_h(\sigma_{\textup{max}})$.
\end{customproof}

\begin{lemma}\label{QM1_beta_pricejump}
    A unique $\beta^* \in (\beta_{\textup{min}},\infty)$ exists such that $p^* = p^*_l(\beta)$ when $\beta \geq \beta^*$ and $p^* = p^*_h(\beta)$ when $\beta \leq \beta^*$.
\end{lemma}
\begin{customproof}
Let $p$ be fixed and consider $\beta_2 > \beta_1 \geq \upsilon_2$. Clearly, $R_1(p^*_l,\beta)$ is constant in $\beta$. Now let $x$ be the unique solution to $x = \lambda P_2(p + h W(x),\beta)$.  For fixed $x$, we have $\lambda P_2(p + h W(x),\beta_1) > \lambda P_2(p + h W(x),\beta_2)$. Hence, $\gamma_2^*(p,\beta_1) > \gamma_2^*(p,\beta_2)$. But then also $p\gamma^*_2(p,\beta_1) > p\gamma^*_2(p,\beta_2)$. Since $p$ was chosen arbitrarily, the following holds:
$$R_2(p^*_h(\beta_1),\beta_1) > R_2(p^*_h(\beta_2),\beta_1) > R_2(p^*_h(\beta_2),\beta_2).$$
Hence, $R_2(p^*_h,\beta)$ is decreasing in $\beta$.

Consider $\beta = \upsilon_2 = \mu + \frac{\sigma^2}{\mu}$. However, from Lemma \ref{QM1_sigma_pricejump} it then follows that $p^* = p^*_h(\upsilon_2)$. Now consider $\beta = \infty$. Then $R_2(p,\infty) = 0$ for any $p$, so also for $p^*_h(\infty)$. Hence, $p^* = p^*_l(\infty)$.
\end{customproof}

\subsection{Proof of Theorem \ref{th:queue2}}\label{ap:queue2}
Recall that $W$ is an increasing function and so if $x < x'$, then also $\frac{R-p}{W(x)} > \frac{R-p}{W(x')}$. This implies that $\lambda \cdot g(x)$ is non-increasing in $x$. Because $x$ is increasing in $x$, it follows that the equation $x = \lambda \cdot \inf_{\Prob \in \cP} \Prob(H \leq \frac{R-p}{W(x)}$ has a unique solution $x^*$. We will next argue that $x^* \geq \gamma^*_H$. Recall that $\gamma^*_H = \inf_{\Prob \in \bP} \gamma_H(p,\Prob)$ with $\gamma_H(p,\Prob)$ the unique solution to the equation $\gamma_H = \lambda \cdot \Prob(H \leq \frac{R-p}{W(\gamma_H)})$ for the given distribution $\bP$. Because of the second assumption we know that there is a distribution $\Prob_{x^*}$ so that $x^* = \lambda \cdot \inf_{\Prob \in \cP} \Prob(H<\frac{R-p}{W(x^*)}) = \lambda \cdot \Prob_{x^*}(H<\frac{R-p}{W(x^*)}).$ This means that  $\gamma_H^* \leq x^* < \infty$. We will next argue that $\gamma^* = x^*$ using a proof by contradiction. Suppose that equality does not hold, i.e., $x^* > \gamma^*_H$. Then there exists a distribution $\Prob' \in \cP$ where $\gamma_H' = \gamma_H'(p,\Prob')$  satisfies
$\gamma_H' = \lambda \cdot \Prob'(H<\frac{R-p}{W(\gamma'_H)})$ and $\gamma^*_H \leq \gamma_H' < x^*$. Because the function $x \mapsto \Prob'(H<\frac{R-p}{W(x)})$ is non-increasing, it follows that $\Prob'(H<\frac{R-p}{W(x^*)}) \leq \Prob'(H<\frac{R-p}{W(\gamma_H')}) = \frac{\gamma_H'}{\lambda} <  \frac{x^*}{\lambda} = \Prob_{x^*}(H<\frac{R-p}{W(x^*)}).$ By definition of $\Prob_{x^*}$, i.e., the second assumption of the theorem statement, this then implies that $\Prob'(H<\frac{R-p}{W(x^*)})  < \inf_{\Prob \in \cP} \Prob(H<\frac{R-p}{W(x^*)})$ which is a contradiction by definition of the infimum. This shows that $x^* = \gamma^*_H$, completing the proof.

\subsection{Proof of Theorem \ref{THMswithQ2}}\label{ap:proofTHMQ2}
We establish Theorem \ref{THMswithQ2} by combining Lemma \ref{QM2_sigma_pricejump} and \ref{QM2_beta_pricejump}, which we will provide below. To begin, we define key quantities that will be used throughout the proofs. We utilize the tail bound from Lemma \ref{th:tutb}, where we set $\lsigma = \usigma = \sigma$:
\begin{flalign}\label{P_i2}
    &&\bar{P}(t,\sigma,\beta) := \inf_{\bP \in \cP(\mu,\sigma,\beta)}\bP(H < t)=
    \begin{cases}
     \bar{P}_{2}(t,\sigma,\beta) := \frac{(\beta-\mu)(t-\mu)+\sigma^2}{\beta t},\quad& t\in[\upsilon_1,\upsilon_2], \\
     \bar{P}_{1}(t,\sigma,\beta) :=\frac{(\mu-t)^2}{(\mu-t)^2+\sigma^2}, \quad& t\in[\upsilon_2,\beta].
    \end{cases}&&
 \end{flalign}

\noindent By utilizing \eqref{P_i2}, we obtain expected revenue. However, since this quantity is decreasing in $p$, the ordering of the regimes will flip. This results in
\begin{flalign}\label{R_i2}
    &&\bar{R}(p,\sigma,\beta) = \begin{cases}
        \bar{R}_{1}(p,\sigma,\beta) := p\lambda \bar{P}_{1}\left(\frac{R-p}{W(\gamma^*_{H1})}),\sigma,\beta\right), \quad & p \in (0,\upsilon^H_1], \\
     \bar{R}_{2}(p,\sigma,\beta):=p\lambda \bar{P}_{2}\left(\frac{R-p}{W(\gamma^*_{H2})}),\sigma,\beta\right), \quad & p \in [\upsilon^H_1,\upsilon^H_2],
    \end{cases}&&
    \end{flalign}
where $\upsilon^H_1$ is the unique solution to $\upsilon^H_1 = R-\upsilon_2W(\gamma^*_{H1}(\upsilon^H_1))$, $\upsilon^H_2$ the unique solution to $\upsilon^H_2 = R-\upsilon_1W(\gamma^*_{H2}(\upsilon^H_2))$, and $\gamma^*_{Hi} = \gamma^*_{Hi}(p,\sigma,\beta)$ the unique solution to $$\gamma^*_{Hi} = \lambda \bar{P}_i\left(\frac{R-p}{W(\gamma^*_{Hi})},\sigma,\beta\right), \text{ for } i = 1, 2.$$ 

\noindent Denote $p^* := \arg\sup_p\bar{R}(p,\sigma,\beta)$ and the two potential optimal prices as $$p^*_l(\sigma,\beta):=\arg\sup_{p\in(0,\upsilon^H_1]}\bar{R}_1(p,\sigma,\beta), \ \ p^*_h(\sigma,\beta) := \arg\sup_{p\in[\upsilon^H_1,\upsilon^H_2]}\bar{R}_2(p,\sigma,\beta).$$
For brevity, we sometimes slightly abuse notation by omitting $\sigma$ or $\beta$ as a function argument when it is not relevant to a proof. This is purely a notational simplification with the omitted variable remaining fixed.

\begin{lemma}\label{QM2_sigma_pricejump}
    A unique $\sigma^*\in(0,\sigma_{\textup{max}})$ exists such that $p^* = p^*_{l}(\sigma)$ when $\sigma \leq \sigma^*$ and $p^* = p^*_{h}(\sigma)$ when $\sigma \geq \sigma^*$.
\end{lemma}
\begin{customproof}
Let $p$ be fixed and consider $\sigma_2 > \sigma_1 \geq 0$. Now let $x$ be the unique solution to $x = \lambda \bar{P}_{1}(\frac{R-p}{W(x)},\sigma)$. Clearly, for fixed $x$, we have $\lambda \bar{P}_{1}(\frac{R-p}{W(x)},\sigma_1) > \lambda \bar{P}_{1}(\frac{R-p}{W(x)},\sigma_2)$. Hence, $\gamma_{H1}^*(p,\sigma_1) > \gamma_{H1}^*(p,\sigma_2)$. But then also $p\gamma^*_{H1}(p,\sigma_1) > p\gamma^*_{H1}(p,\sigma_2)$. Since $p$ was chosen arbitrarily, the following holds:
$$\bar{R}_{1}(p^*_{l}(\sigma_1),\sigma_1) > \bar{R}_{1}(p^*_{l}(\sigma_2),\sigma_1) > \bar{R}_{1}(p^*_{l}(\sigma_2),\sigma_2).$$
Hence, $\bar{R}_{1}(p^*_{l},\sigma)$ is decreasing in $\sigma$.

Now let $y$ be the unique solution to $y = \lambda \bar{P}_{2}(\frac{R-p}{W(y)},\sigma)$. Clearly, for fixed $y$, we have $\lambda \bar{P}_{2}(\frac{R-p}{W(y)},\sigma_1) < \lambda \bar{P}_{2}(\frac{R-p}{W(y)},\sigma_2)$. Hence, $\gamma_{H2}^*(p,\sigma_1) < \gamma_{H2}^*(p,\sigma_2)$. But then also $p\gamma^*_{H2}(p,\sigma_1) < p\gamma^*_{H2}(p,\sigma_2)$. Since $p$ was chosen arbitrarily, the following holds:
$$\bar{R}_{2}(p^*_{h}(\sigma_1),\sigma_1) < \bar{R}_{2}(p^*_{h}(\sigma_1),\sigma_2) < \bar{R}_{2}(p^*_{h}(\sigma_2),\sigma_2).$$
Hence, $\bar{R}_{2}(p^*_{h},\sigma)$ is increasing in $\sigma$.

Now consider $\sigma = 0$. Then $\frac{R-p^*_{h}(0)}{W(\gamma^*_{H2})} = \mu$, as $\upsilon_1 = \upsilon_2 = \mu$, and consequently, $\bar{R}_{2}(p^*_{h}(0),0) = 0$. Hence, since $\sup_{p \in (0,\upsilon_1^H]}\bar{R}_{1}(p,0)>0$, we have that $p^* = p^*_{l}(0)$. We now consider $\sigma = \sigma_{\textup{max}} = \sqrt{\mu(\beta-\mu)}$. Then $\frac{R-p^*_{l}(\sigma_{\textup{max}})}{W(\gamma^*_{H1})} = 0$ and consequently, $\bar{R}_{1}(p^*_{l}(\sigma_{\textup{max}}),\sigma_{\textup{max}}) = 0$. Hence, as $\sup_{p \in (\upsilon_1^H,\upsilon_2^H]}\bar{R}_{2}(p,\sigma_{\textup{max}})>0$, we have that $p^* = p^*_{h}(\sigma_{\textup{max}})$.
\end{customproof}

\begin{lemma}\label{QM2_beta_pricejump}
    A unique $\beta^* \in (\beta_{\textup{min}},\infty)$ exists such that $p^* = p^*_l(\beta)$ when $\beta \geq \beta^*$ and $p^* = p^*_h(\beta)$ when $\beta \leq \beta^*$.
\end{lemma}
\begin{customproof}
Let $p$ be fixed and consider $\beta_2 > \beta_1 \geq \upsilon_2$. Clearly, $\bar{R}_1(p^*_l,\beta)$ is constant in $\beta$. Now let $x$ be the unique solution to $x = \lambda \bar{P}_2(\frac{R-p}{W(x)},\beta)$. Notice that for fixed $x$ we have that $t = \frac{R-p}{W(x)} \leq \upsilon_2 = \mu + \frac{\sigma^2}{\mu}$, which means
$$\frac{\partial \bar{P_2}(t,\beta)}{\partial \beta} =\frac{t\mu - \mu^2 - \sigma^2}{\beta^2 t} \leq 0.$$

Therefore, $\lambda \bar{P}_2(\frac{R-p}{W(x)},\beta_1) > \lambda \bar{P}_2(\frac{R-p}{W(x)},\beta_2)$, and hence, $\gamma_{H2}^*(p,\beta_1) > \gamma_{H2}^*(p,\beta_2)$. But then also $p\gamma^*_{H2}(p,\beta_1) > p\gamma^*_{H2}(p,\beta_2)$. Since $p$ was chosen arbitrarily, the following holds:
$$\bar{R}_2(p^*_h(\beta_1),\beta_1) > \bar{R}_2(p^*_h(\beta_2),\beta_1) > \bar{R}_2(p^*_h(\beta_2),\beta_2).$$
Hence, $\bar{R}_2(p^*_h,\beta)$ is decreasing in $\beta$.

Consider $\beta = \upsilon_2 = \mu + \frac{\sigma^2}{\mu}$. This is equivalent to $\sigma^2 = \mu(\beta-\mu)$. However, from Lemma \ref{QM2_sigma_pricejump} it then follows that $p^* = p^*_h(\upsilon_2)$. Hence, $p^* = p^*_h(\upsilon_2)$. Now consider $\beta = \infty$. Then $\bar{R}_2(p,\infty) = 0$ for any $p$, so also for $p^*_h(\infty)$. Hence, $p^* = p^*_l(\infty)$.
\end{customproof}

\section{Worst-case demand and revenue}\label{9plotsexamples}
In Section~\ref{sec:pre} of the main paper, we proved and illustrated in Figure~\ref{fig:1demrev} that as the standard deviation increases, the local maximum of the higher segment overtakes that of the lower segment. Here, we provide additional results for the worst-case demand function \eqref{ttb} and the corresponding worst-case revenue \eqref{rev_g}, further demonstrating this overtaking effect visually.

\newpage
\begin{figure}[h!]
    \centering
{\includegraphics[width=0.75\linewidth]{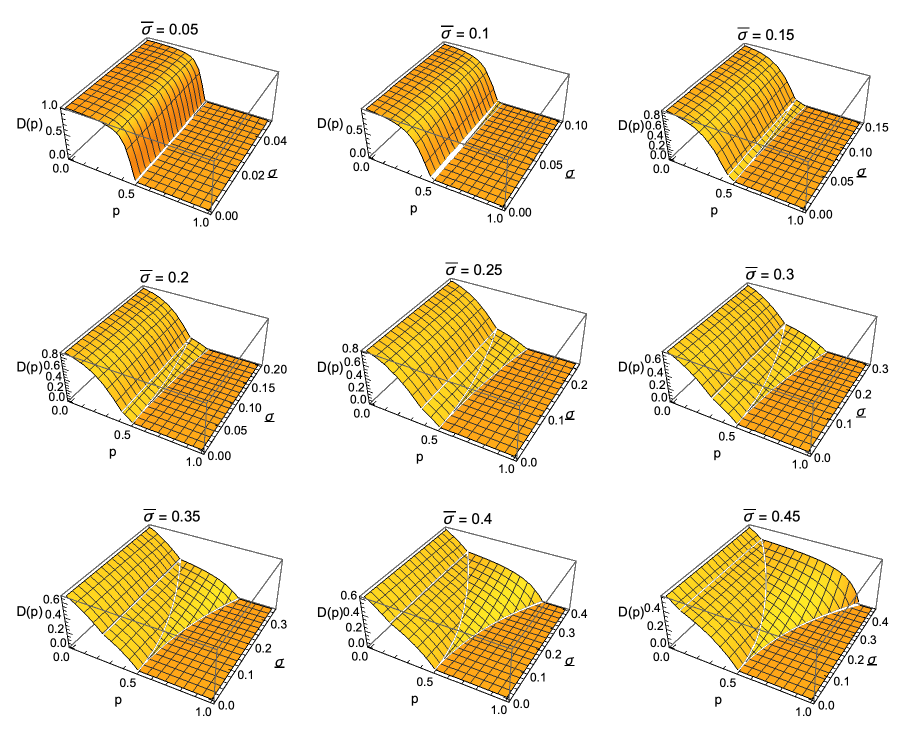}}
\caption{Worst-case demand function  \eqref{ttb} for $\mu=0.5$, $\beta=1$ and $\usigma$ ranging from 0.05 to 0.45.}
    \label{fig:9demand}
\end{figure}

\begin{figure}[h!]
    \centering
{\includegraphics[width=0.75\linewidth]{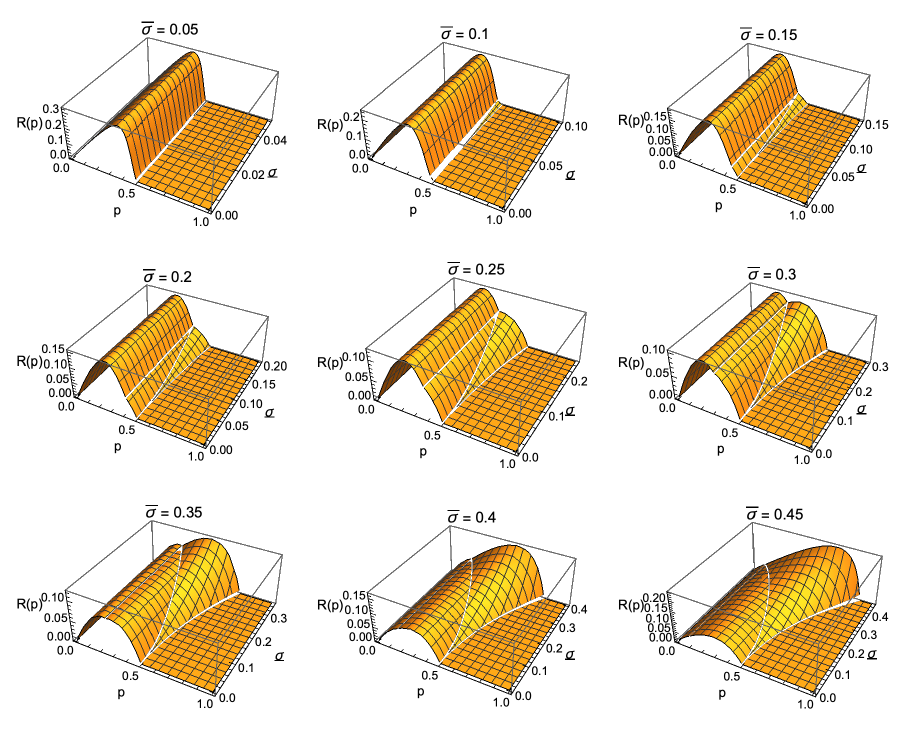}}
\caption{Worst-case revenue \eqref{rev_g} for $\mu=0.5$, $\beta=1$ and $\usigma$ ranging from 0.05 to 0.45.}
    \label{fig:9revenue}
\end{figure}
\color{black}

\newpage
\section{Comparison with ground truth markets}\label{app:gtms}

\begin{figure}[H]
\centering
\subfigure[$\mu=0.4$]{
\begin{tikzpicture}
    \node at (0,0) {\includegraphics[width=0.25\linewidth]{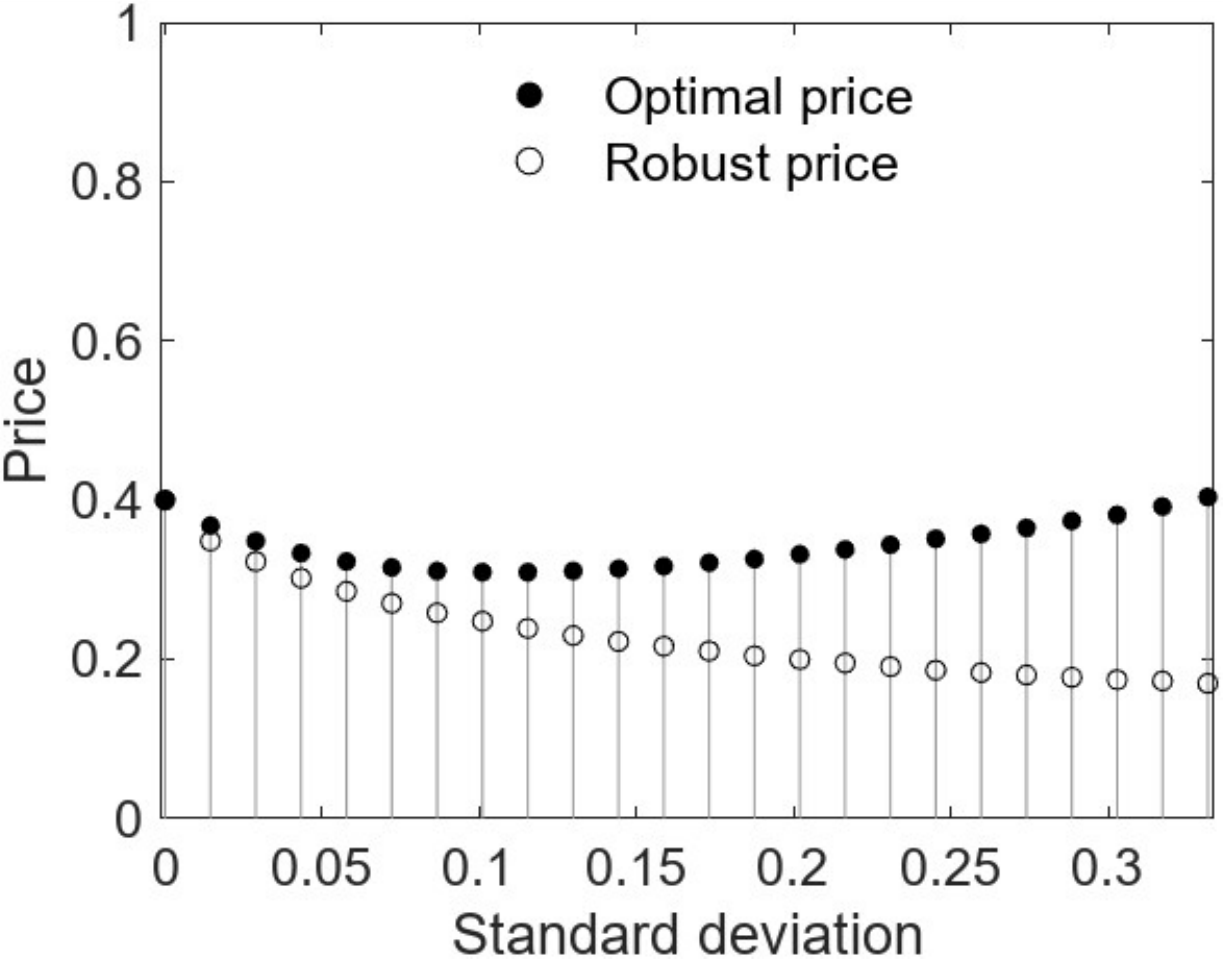}};
\end{tikzpicture}
}
\subfigure[$\mu=0.5$]{
\begin{tikzpicture}
    \node at (0,0) {\includegraphics[width=0.25\linewidth]{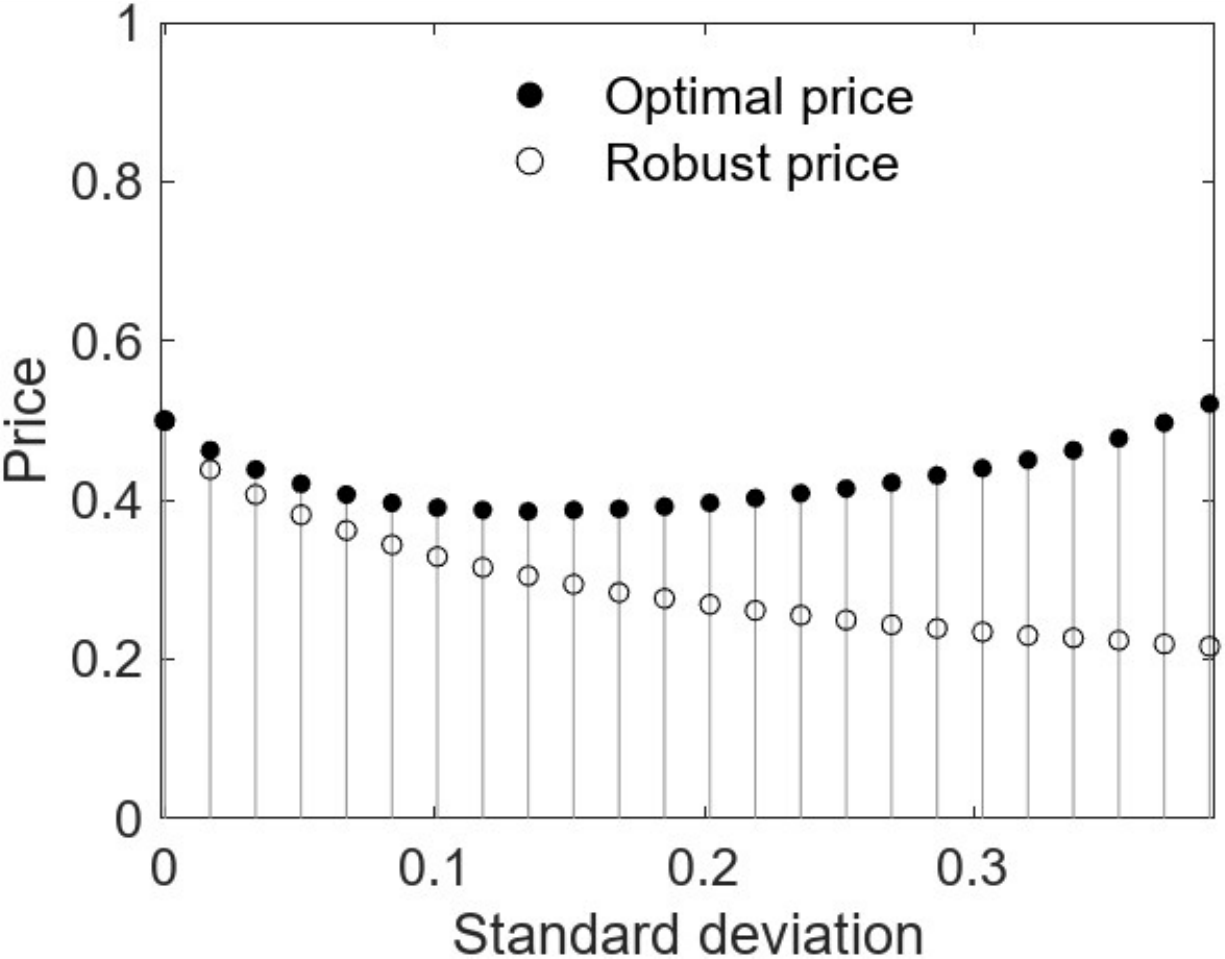}};
\end{tikzpicture}
}
\subfigure[$\mu=0.8$]{
\begin{tikzpicture}
    \node at (0,0) {\includegraphics[width=0.25\linewidth]{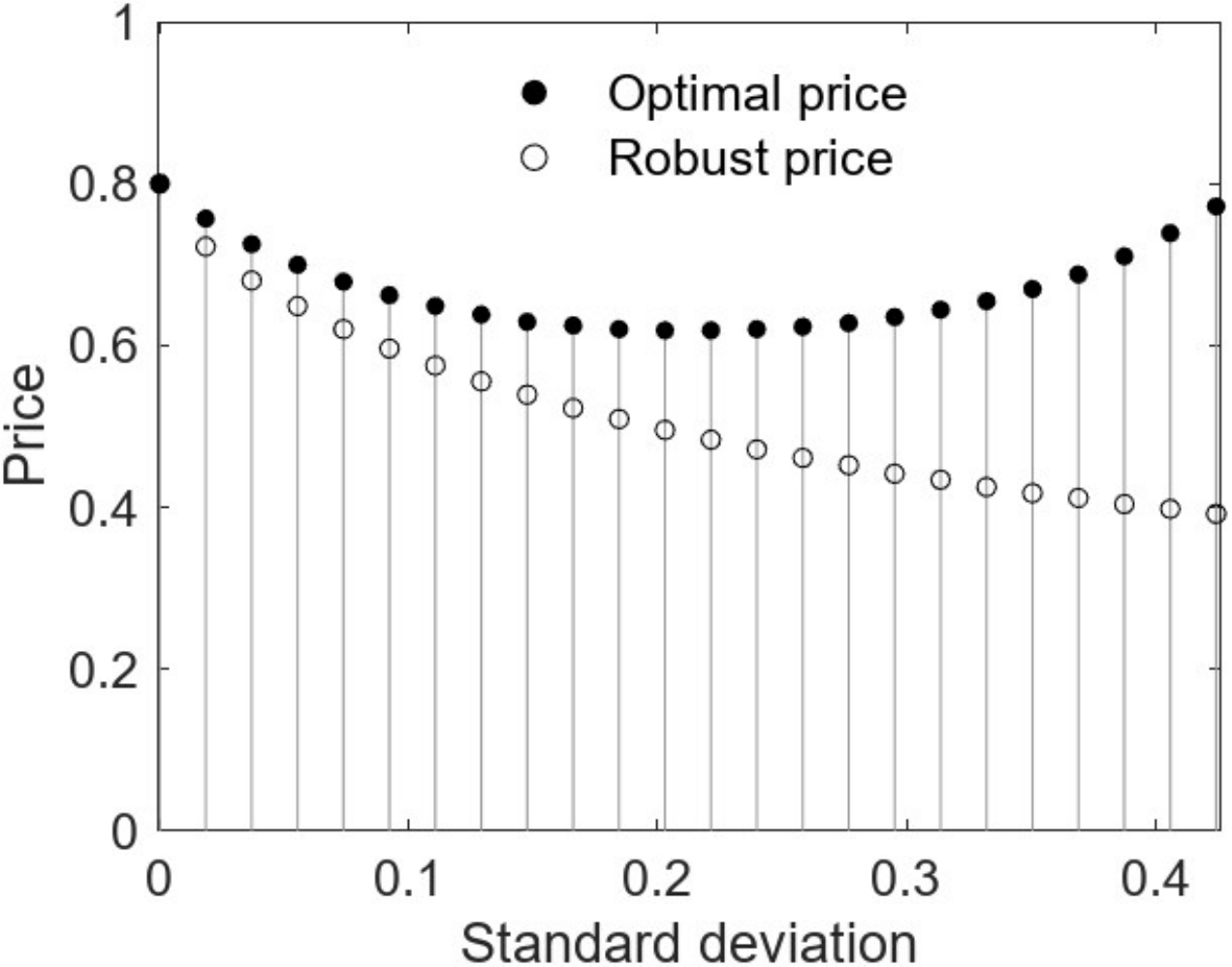}};
\end{tikzpicture}
}

\subfigure[$\mu=0.4$]{
\begin{tikzpicture}
    \node at (0,0) {\includegraphics[width=0.25\linewidth]{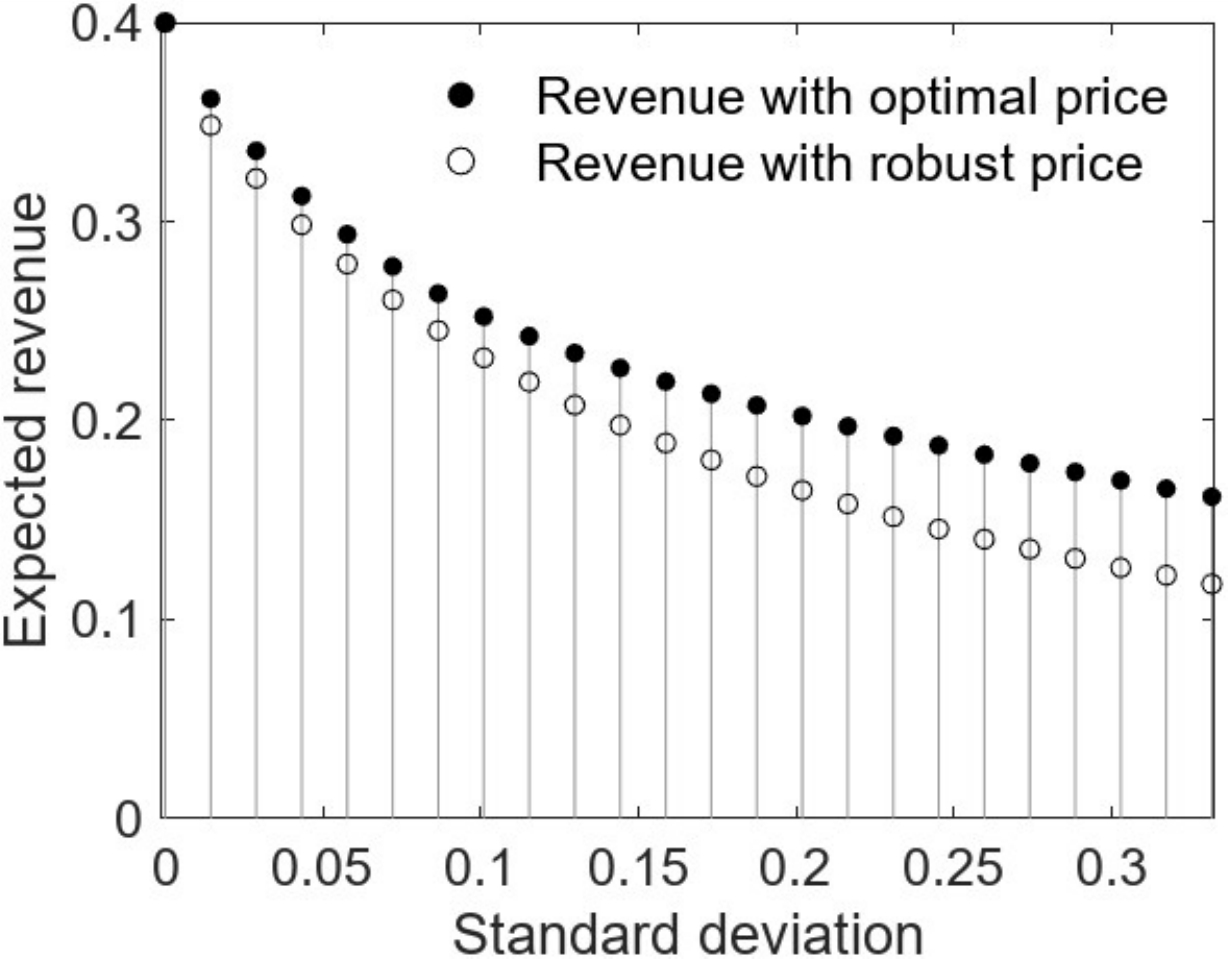}};
\end{tikzpicture}
}
\subfigure[$\mu=0.5$]{
\begin{tikzpicture}
    \node at (0,0) {\includegraphics[width=0.25\linewidth]{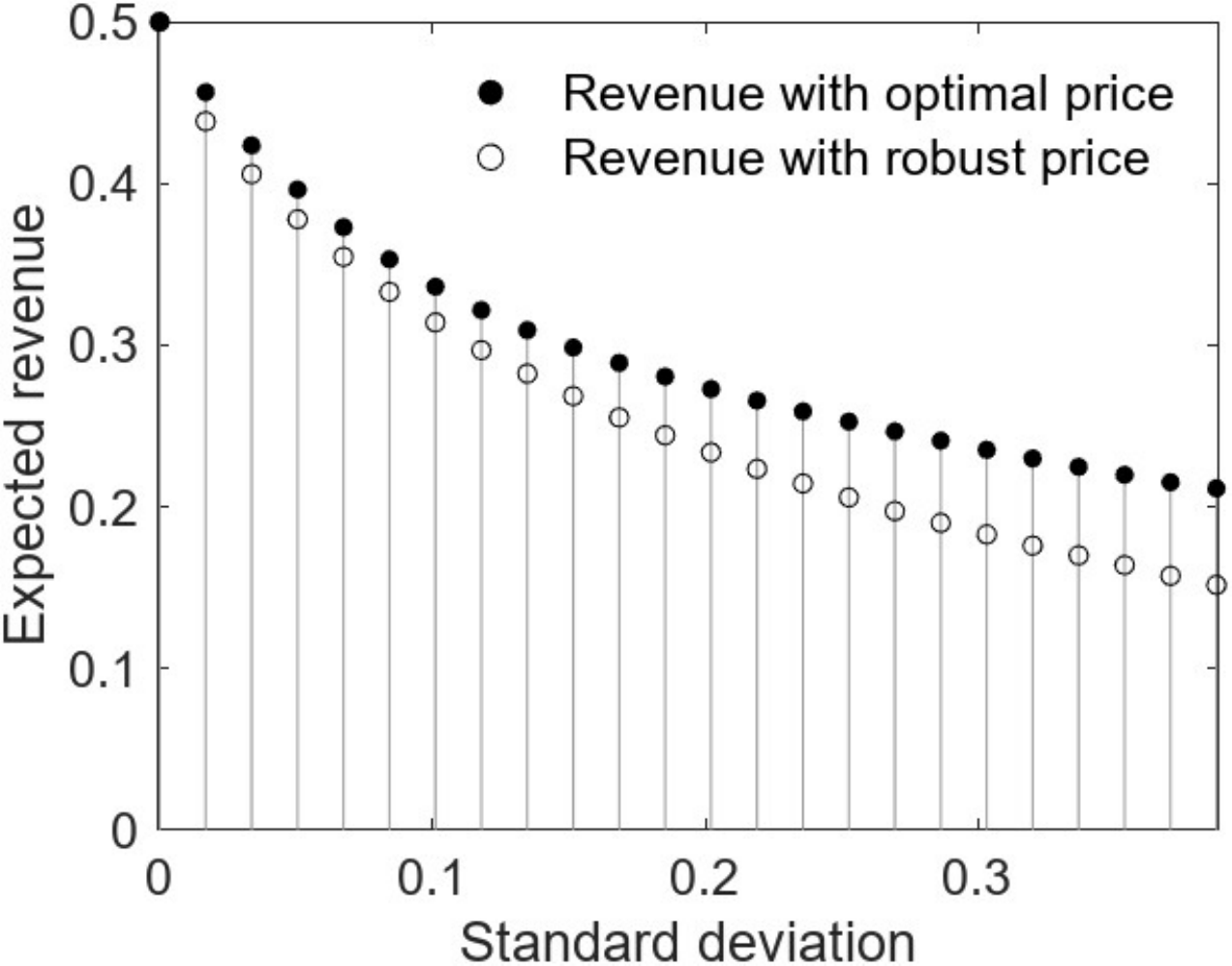}};
\end{tikzpicture}
}
\subfigure[$\mu=0.8$]{
\begin{tikzpicture}
    \node at (0,0) {\includegraphics[width=0.25\linewidth]{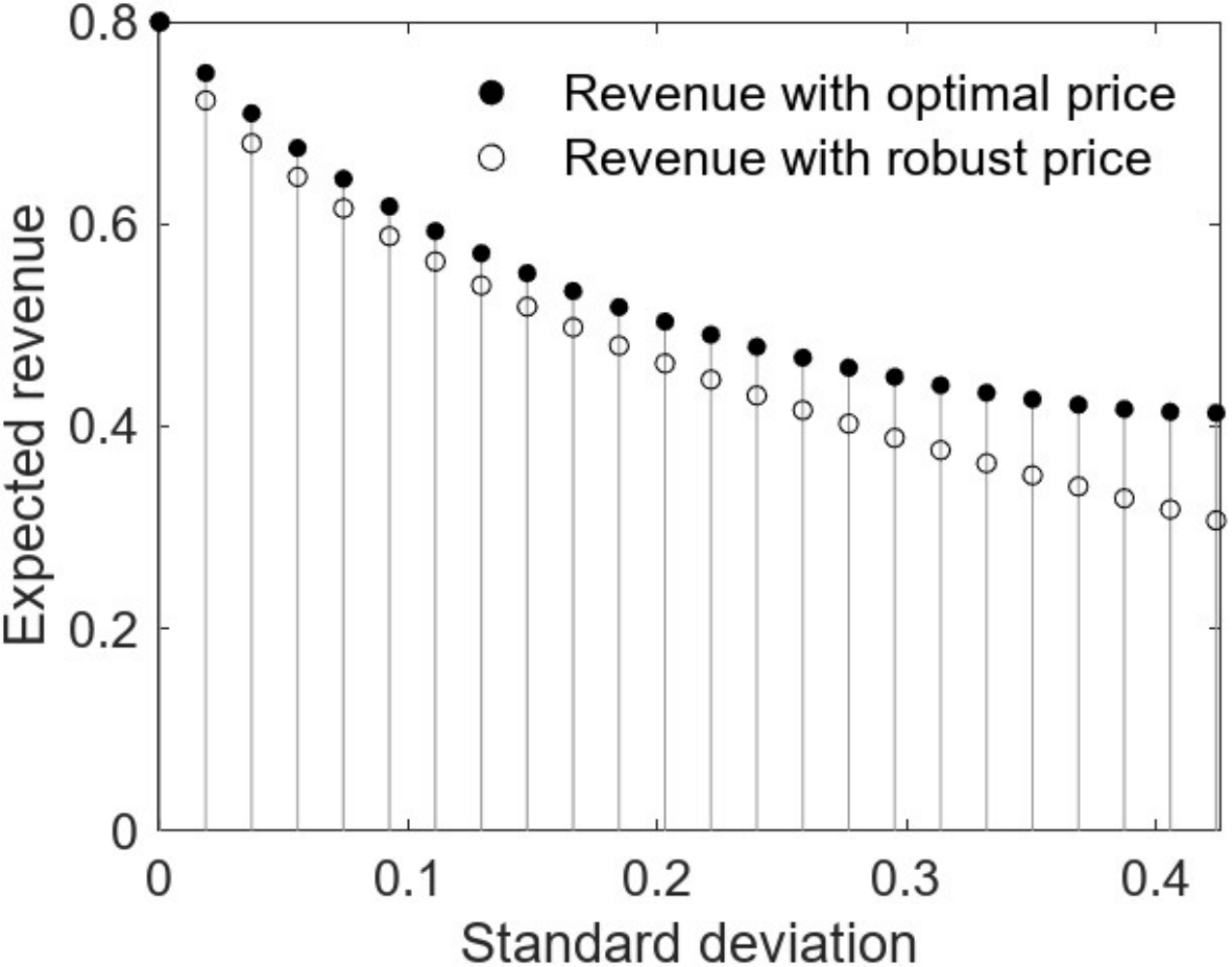}};
\end{tikzpicture}
}
\caption{Optimal and robust prices (top row) and expected revenue (bottom row) for truncated normal ground truth distribution with $\beta = 1.5$, several $\mu$ values and all possible $\sigma$ values.}
\label{fig:gt_tn}
\end{figure}

\begin{figure}[H]
\centering
\subfigure[$\mu=0.4$]{
\begin{tikzpicture}
    \node at (0,0) {\includegraphics[width=0.25\linewidth]{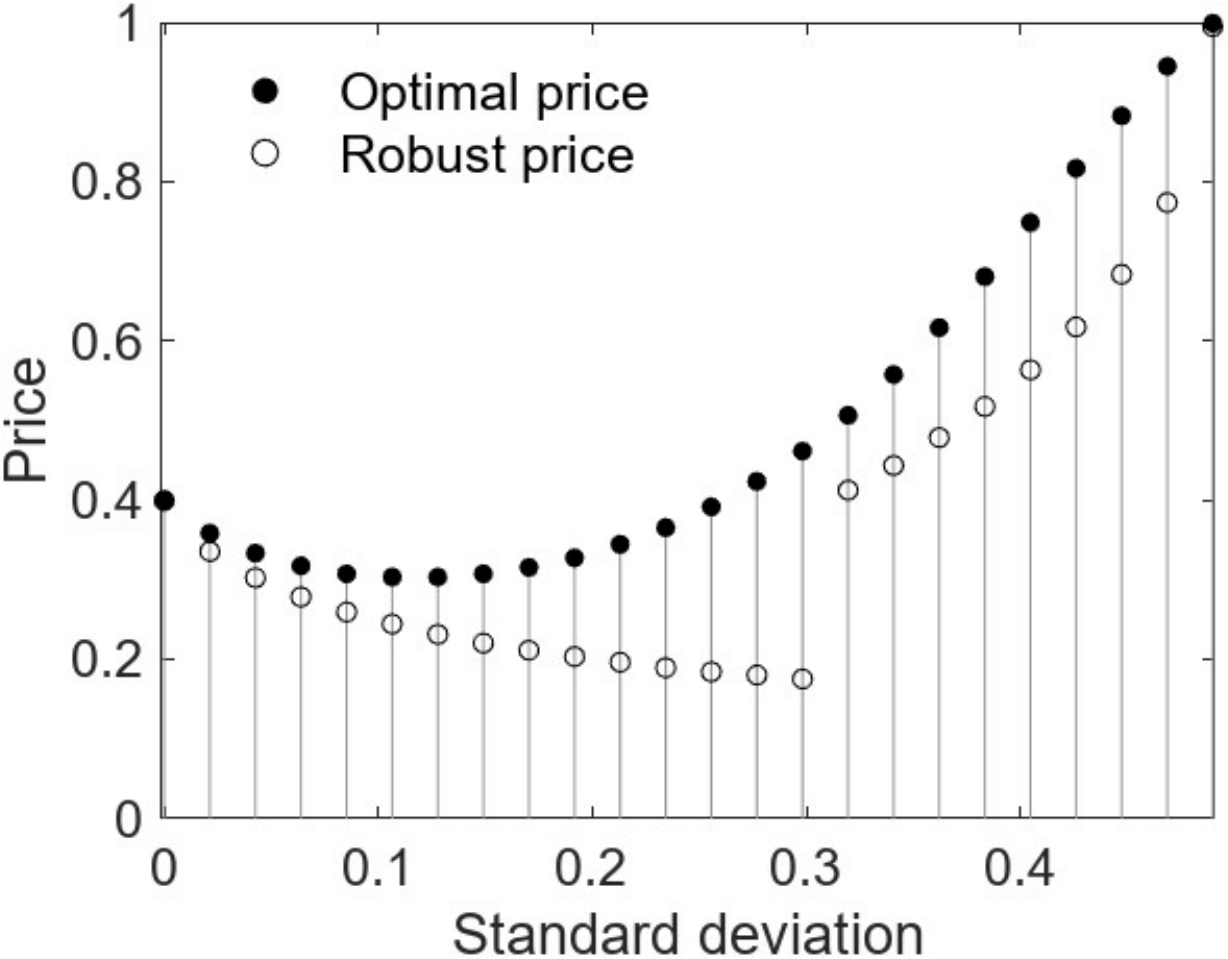}};
\end{tikzpicture}
}
\subfigure[$\mu=0.5$]{
\begin{tikzpicture}
    \node at (0,0) {\includegraphics[width=0.25\linewidth]{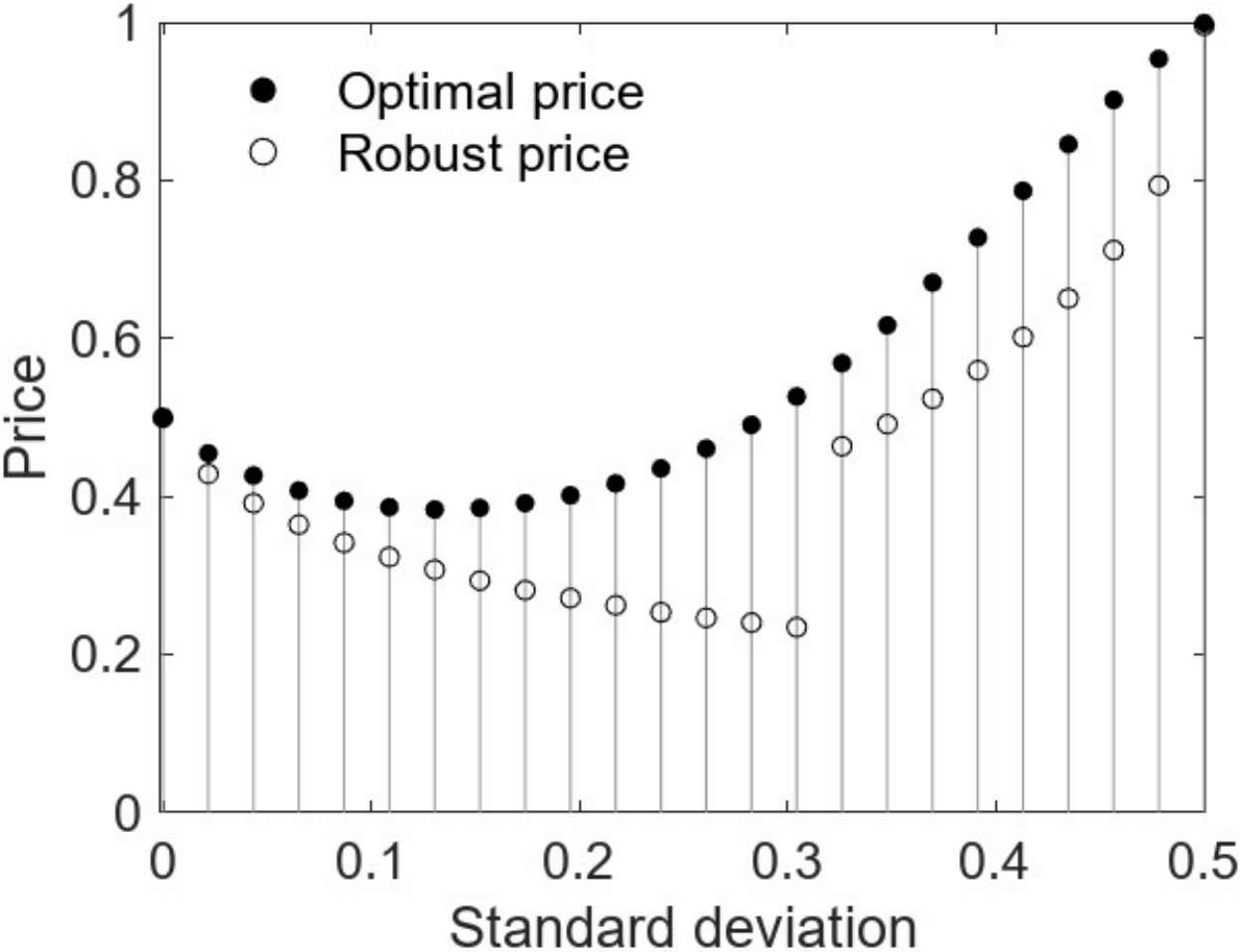}};
\end{tikzpicture}
}
\subfigure[$\mu=0.8$]{
\begin{tikzpicture}
    \node at (0,0) {\includegraphics[width=0.25\linewidth]{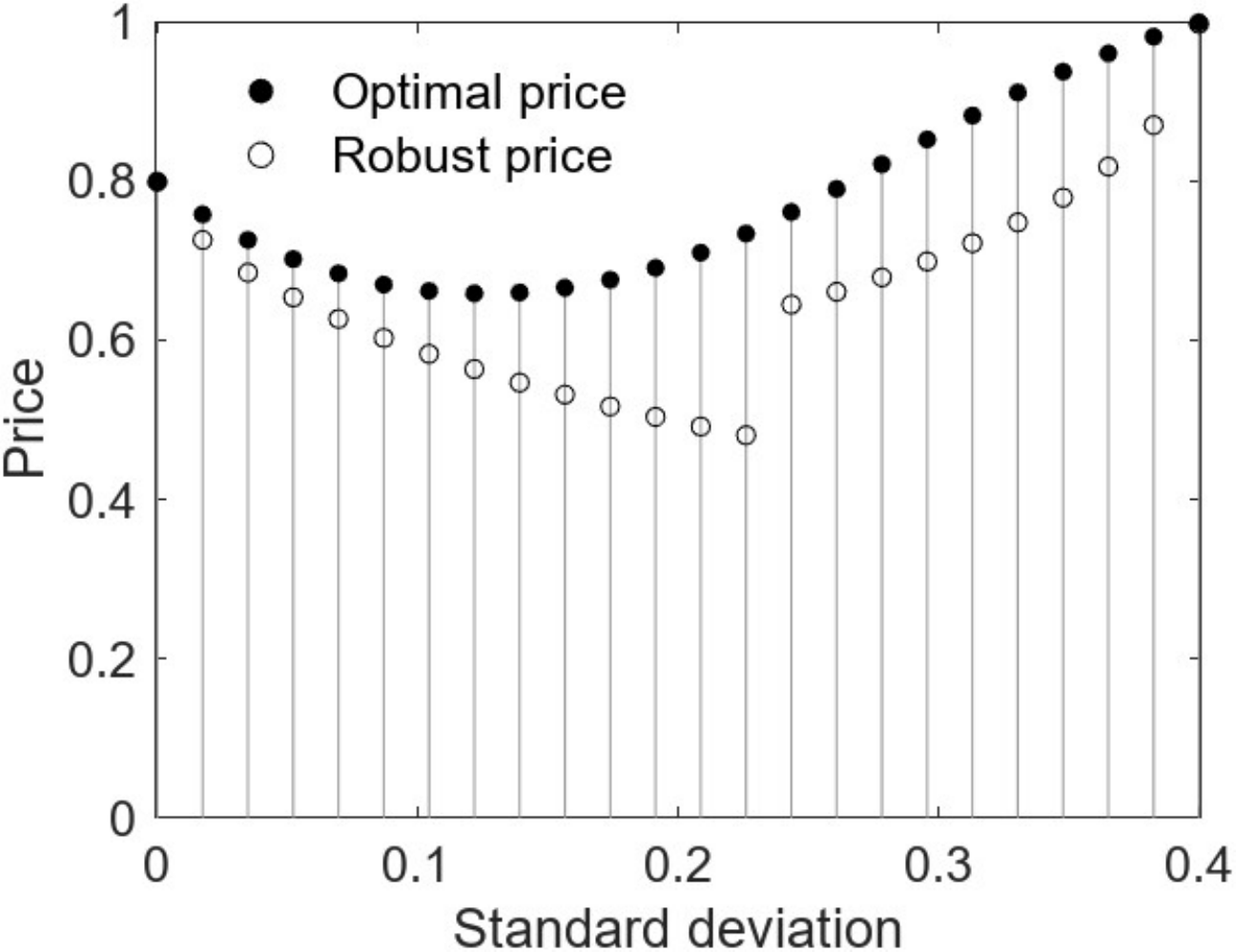}};
\end{tikzpicture}
}

\subfigure[$\mu=0.4$]{
\begin{tikzpicture}
    \node at (0,0) {\includegraphics[width=0.25\linewidth]{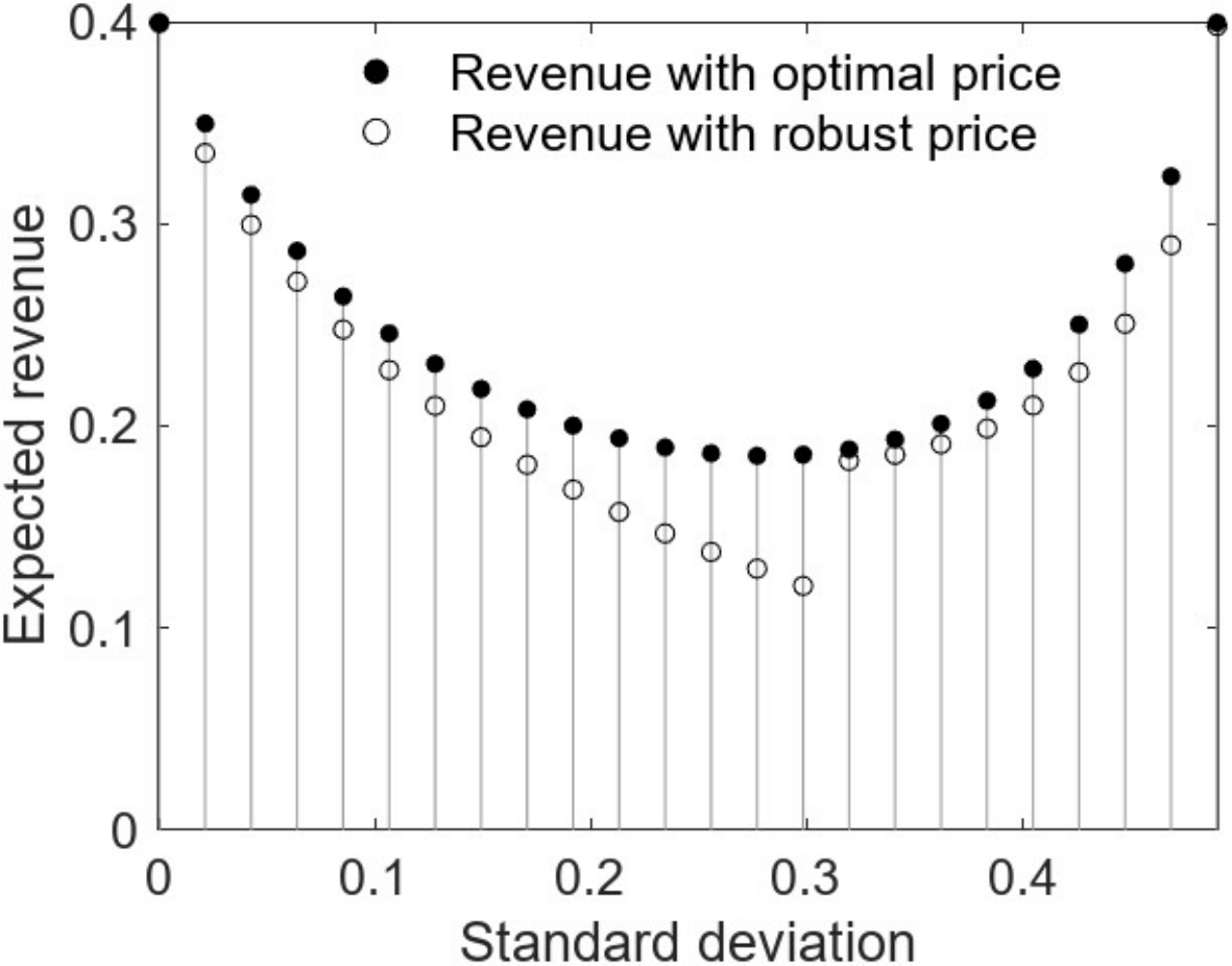}};
\end{tikzpicture}
}
\subfigure[$\mu=0.5$]{
\begin{tikzpicture}
    \node at (0,0) {\includegraphics[width=0.25\linewidth]{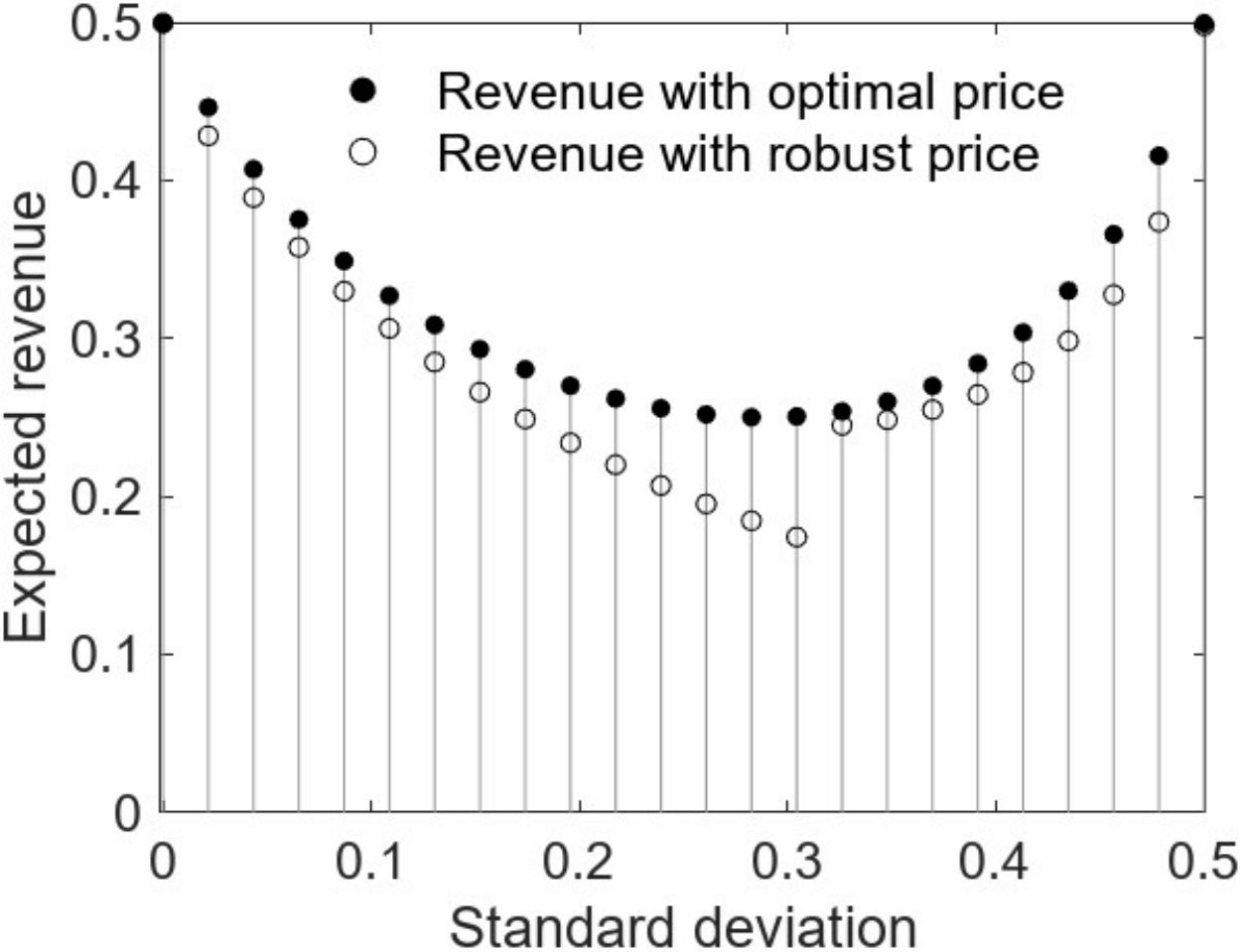}};
\end{tikzpicture}
}
\subfigure[$\mu=0.8$]{
\begin{tikzpicture}
    \node at (0,0) {\includegraphics[width=0.25\linewidth]{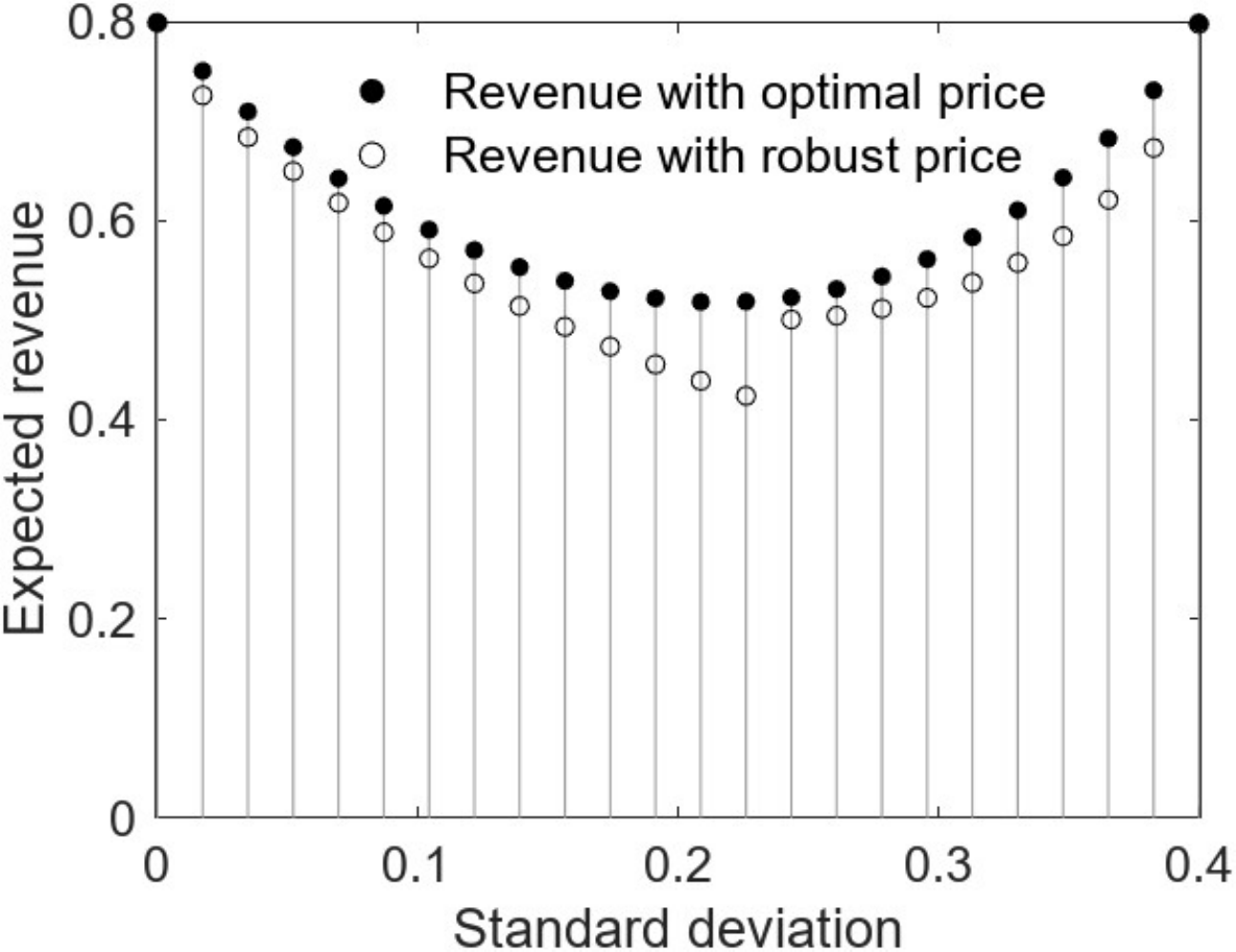}};
\end{tikzpicture}
}
\caption{Optimal and robust prices (top row) and expected revenue (bottom row) under a beta ground truth distribution for several $\mu$ values and all possible $\sigma$ values.}
\label{fig:gt_beta}
\end{figure}

\end{appendices}
\end{document}